\documentclass[12pt,a4paper,reqno]{amsart}
\usepackage{amssymb}
\usepackage{graphicx}
\usepackage{enumitem}

\input epsf.tex

\def\myitemmargin{
\leftmargini=25pt 
\leftmarginii=15pt 
\leftmarginiii=10pt 
\leftmarginiv=10pt  
}

\newdimen\xsize
\newdimen\oldbaselineskip
\newdimen\oldlineskiplimit
\def\restorelineskip{\baselineskip=\oldbaselineskip%
\lineskiplimit=\oldlineskiplimit}
\def\putm[#1][#2]#3{
\hbox{\vbox to 0pt{\parindent=0pt%
\vskip#2\xsize\hbox to0pt{\hskip#1\xsize $#3$\hss}\vss}}}%
\def\putt[#1][#2]#3{
\vbox to 0pt{\noindent\hskip#1\xsize\lower#2\xsize%
\vtop{\restorelineskip#3}\vss}}

\DeclareFontFamily{U}{rsf}{\skewchar\font'177}%
\DeclareFontShape{U}{rsf}{m}{n}{<-6>rsfs5<6-8>rsfs7<8->rsfs10}{}%
\DeclareFontShape{U}{rsf}{b}{n}{<-6>rsfs5<6-8>rsfs7<8->rsfs10}{}%
\DeclareMathAlphabet\RSFS{U}{rsf}{m}{n}
\SetMathAlphabet\RSFS{bold}{U}{rsf}{b}{n}
\let\scr=\rfs
\def\mib#1{\boldsymbol{#1}} 
\def\sf#1{{\mathsf{#1}}} 

\def\slsf#1{{\slshape \sffamily #1\/}}

\def\fr#1{{\mathfrak{#1}}}

\def\msmall#1{\mathchoice{\hbox{\small$\displaystyle {#1}$}}{#1}{#1}{#1}}

\let\3=\ss
\def\cc{{\mathbb C}}

\def\rr{{\mathbb R}}
\def\nn{{\mathbb N}}
\def\pp{{\mathbb P}}
  \def\cp{{\cc\pp}}
  \def\rp{{\rr\pp}}
\def\zz{{\mathbb Z}}

\def\aut{{\sf{Aut}}}
   \def\inn{{\sf{Inn}}}
   
   \def\out{{\sf{Out}}}

\def\st{{\mathsf{st}}}

\def\br{\sf{Br}}
\def\pbr{\sf{PBr}}

\def\cos{\sf{cos}\,}

\def\dim{\sf{dim}\,}

\def\ind{\sf{ind}}

\def\sym{\sf{Sym}}

\def\frg{{\fr{g}}}

\def\sfh{{\sf{H}}}

\def\id{\sf{Id}}

\def\sfk{{\sf{K}}}

\def\ker{\sf{Ker}\,}
\def\lim{\mathop{\sf{lim}}}

\def\mod{\;\sf{mod}\,}
\def\sfm{{\sf{M}}}  

\def\pr{\sf{pr}}
\def\pt{\sf{pt}}
\def\rank{\sf{rank}}

\def\sin{\sf{sin}}
\def\sin{\sf{sin}\,}

\def\SO{{\mib{SO}}}

\def\Gl{\mib{Gl}}

\def\sfw{\sf{W}}

\def\mbfe{{\mib{e}}}

\def\mbfs{{\mib{s}}}

\def\mbfv{{\mib{v}}}

\def\mbfx{{\mib{x}}}
\def\mbfy{{\mib{y}}}

\def\bfe{{\mib{E}}}
\def\bff{{\mib{F}}}
\def\bfl{{\mib{L}}}

\def\bfo{{\mib{O}}}
\def\bfp{{\mib{P}}}

\def\bfs{{\mib{S}}}

\def\bfu{{\mib{U}}}

\def\bfy{{\mib{Y}}}

\def\<{\langle}\let\lan=\<
\def\>{\rangle}\let\ran=\>
 \let\bs=\bss
\def\bsq{{\mathchoice{\rule{3.5pt}{3.5pt}}%
{\rule{3.5pt}{3.5pt}}%
{\rule{3.pt}{3.pt}}%
{\rule{2.2pt}{2.2pt}}%
}}

\def\dbar{{\barr\partial}}
\def\barz{{\bar z}}

\def\ddef{\mathrel{{=}\raise0.3pt\hbox{:}}}
\def\deff{\mathrel{\raise0.3pt\hbox{\rm:}{=}}}
\def\inv{^{-1}}

\def\fraction#1/#2{\mathchoice{{\msmall{ \frac{#1}{#2} }}}%
{{ \frac{#1}{#2} }}{{#1/#2}}{{#1/#2}}}

\def\half{{\fraction1/2}}
\def\leq{\leqslant}
\def\leq{\leqslant}
\def\vph{^{\mathstrut}}

\def\longpoints{\leaders\hbox to 0.5em{\hss.\hss}\hfill \hskip0pt}
\def\stateskip{\smallskip}
\def\state#1. {\stateskip\noindent{\bf#1. }} 
\def\statep#1. {\stateskip\noindent{\bf#1 }} 

\def\proof{\state Proof. \,}

\def\Chi{\raise 2pt\hbox{$\chi$}}
\def\eg{\,{\sl{e.g.}}}
\def\ie{\,{\sl i.e.}}

\def\isl{{\text{\slsf{i}}\mskip1mu}} 
\let\cpi=\isl
\def\sli{{\sl i)\/} \,}         \def\slip{{\sl i)}}
\def\slii{{\sl i$\!$i)\/} \,}   \def\sliip{{\sl i$\!$i)}}
\def\sliii{{\sl i$\!$i$\!$i)\/} \,}    \def\sliiip{{\sl i$\!$i$\!$i)}}
\def\sliv{{\sl i$\!$v)\/} \,}   \def\slivp{{\sl i$\!$v)}}
\def\slv{{\sl v)\/} \,}         \def\slvp{{\sl v)}}
\def\slvi{{\sl v$\!$i)\/} \,}   \def\slvip{{\sl v$\!$i)}}
\def\slvii{{\sl v$\!$i$\!$i)\/} \,}   \def\slviip{{\sl v$\!$i$\!$i)}}

\def\diff{{{\scr D}\mskip -2.8mu i\mskip -5mu f\mskip-6.3mu f\mskip -4mu}}

\def\map{{\mathsf{Map}}}

\def\barr#1{\mskip1mu\overline{\mskip-1mu{#1}\mskip-1mu}\mskip1mu}
\def\symp{{{\scr S}\mskip -3.5mu y\mskip -1.3mu m\mskip-1.3mu p\mskip -.6mu}}

\def\tic{{\ti c}}
\def\tiw{{\ti w}}
\def\tix{{\,\wt{\times}\,}}
\def\Chi{\raise 2pt\hbox{$\chi$}}

\def\iom{{\mathchar"010A}} 

\def\bfSigma{{\boldsymbol{\Sigma}}}
    \def\bfSig{\bfSigma}%

  \def\rma{{\mathrm A}}

   \def\rmb{{\mathrm{B}}}
   \def\rmbd{{\mathrm{BD}}}
   \def\rmbe{{\mathrm{BE}}}

   \def\rmd{{\mathrm{D}}}

   \def\rme{{\mathrm{E}}}

  \def\calf{{\mathcal F}}

  \def\rmi{{\mathrm{I}}}
\def\scrj{\scr{J}}

  \def\calk{{\mathcal K}}
\def\scrl{\scr{L}}
  
  \def\rml{{\mathrm{L}}}
\def\scrm{\scr{M}}

    \def\calm{{\mathcal M}}

\def\scro{{\scr O}}

  \def\cals{{\mathcal S}}

\def\scrz{{\scr Z}}

\let\xrar=\xrightarrow

\def\bs{\backslash}

\def\dbar{{\barr\partial}}
\def\1{{1\mkern-5mu{\rom l}}}

\def\id{\sf{id}}

\def\geq{\geqslant}\def\geq{\geqslant}
\def\inv{^{-1}}
\let\wh=\widehat
\let\wt=\widetilde
\def\half{{\fraction1/2}}
\def\leq{\leqslant}
\def\vph{^{\mathstrut}}

\def\ti#1{{\tilde{#1}}}

\def\comment#1\endcomment{}

 \belowdisplayskip=\abovedisplayskip

\def\lineeqqno(#1){\hfill\llap{\vbox to 10pt%
{\vss\begin{align} \eqqno(#1)\end{align}\vss}}\vskip1pt}

\def\Matrix{\begin{matrix}}
\def\Endmatrix{\end{matrix}}

\def\Cases{\begin{cases}}
\def\Endcases{\end{cases}}

\advance\headheight 1.2pt

\def\newsection[#1]#2{\section{#2} \vskip-15pt\label{#1}\vskip15pt}
   
\def\newsubsection[#1]#2{\subsection{#2}\vskip10pt\label{#1}\vskip-10pt}
   \def\refsubsection#1{\slsf{\S\ref{#1}}}

\newtheorem{thma}{Theorem}
   \def\newthma#1{\begin{thma} \label{#1}}
\newtheorem{lem}{Lemma}[section]
   \def\newlemma#1{\begin{lem} \label{#1}}
   \def\lemma#1{{\it \slsf{Lemma \ref{#1}}}}

\newtheorem{thm}[lem]{Theorem}
   \def\newthm#1{\begin{thm} \label{#1}}
   \def\refthm#1{{\it \slsf{Theorem \ref{#1}$\,$}}}
\newtheorem{prop}[lem]{Proposition}
   \def\newprop#1{\begin{prop} \label{#1}}
   \def\propo#1{{\it \slsf{Proposition \ref{#1}}}}
\newtheorem{corol}[lem]{Corollary}
   \def\newcorol#1{\begin{corol} \label{#1}}
   \def\refcorol#1{{\it \slsf{Corollary \ref{#1}}}}
\newtheorem{defi}{Definition}[section]
   \def\newdefi#1{\begin{defi} \label{#1}\sl }
   \def\refdefi#1{{\it \slsf{Definition \ref{#1}}}}

\newcounter{stepnr}

{\begin{list}%
{ \slsf{Step \arabic{stepnr}.} }%
{\usecounter{stepnr}  \setlength{\leftmargin}{0pt} 
\setlength{\itemindent}{35pt} \setlength{\listparindent}{0pt} 
\setlength{\itemsep}{7pt plus 3pt}
}}%
{\end{list}}

\def\reffig#1{{\it \slsf{Figure \ref{#1}}}}

\def\eqqno(#1){\label{(#1)}}
\def\eqqref(#1){\eqref{(#1)}}

\usepackage{amscd}
\usepackage{hyperref}

\usepackage[all]{xy}

\numberwithin{equation}{section}

\begin{document}

\myitemmargin 
\baselineskip =14.5pt plus 1.5pt

\title[Secondary Stiefel-Whitney class of\/ symplectomorphisms]%
{Secondary Stiefel-Whitney class\\[3pt] and diffeomorphisms\\[3pt]
of rational and ruled symplectic $4$-manifolds}
\author[V.~Shevchishin]{Vsevolod V.~Shevchishin}
\address{Universit\"at Hamburg\\
Department Mathematik\\
Bundesstra\ss e 55\\
 20146 Hamburg\\
Germany
} \email{shevchishin@googlemail.com}
\dedicatory{} \subjclass{} \thanks{This research was carried out with the
 financial support of the Deutsche
 Forschungsgemeinschaft and University of Hamburg.}
\begin{abstract} 
 We introduce the secondary Stiefel-Whitney class $\tiw_2$ of homotopically
 tri\-vial diffeomorphisms and show that a homotopically
 trivial symplectomorphism of a~ruled $4$-manifold is isotopic to identity if
 and only if the class  $\tiw_2$ vanishes. 

 Using this, we describe the combinatorial structure of the diffeotopy
 group of ruled symplectic $4$-manifolds $X$, either minimal or blown-up, and
 its action on the homology and homotopy groups $\sfh_2(X,\zz)$, $\pi_1(X)$,
 $\pi_2(X)$.
\end{abstract}
\maketitle
\setcounter{section}{-1} \pagebreak[1]

\ \\[-60pt]

\newsection[intro]{Introduction}

One of the problems of symplectic geometry is to classify, or at least
distinguish, possible symplectic structures on a given manifold $X$. Except
for the easy case of dimension 2, the only class of manifolds for which such a
classification is known consists of rational and ruled $4$-manifolds.  Recall
that any ruled symplectic $4$-manifold $(X,\omega)$ is either a product $Y{\times}S^2$
of the two-sphere $S^2$ with the surface $Y$ of genus $g(Y)$ or can be
obtained from such a product by a sequence of symplectic blow-ups and
contractions, see \refsubsection{ruled} for more details.  This construction
does not change the fundamental group, so that $\pi_1(X)\cong\pi_1(Y)$, and we call
$g(Y)=\half\rank\sfh_1(Y)=\half\rank\sfh_1(X)$ the \slsf{genus} of $X$. Such a
symplectic ruled $4$-manifold is \slsf{minimal} if it contains no symplectic
exceptional sphere. A ruled manifold is called \slsf{rational} if $g(X)=0$ and
\slsf{irrational} otherwise.

The classification of symplectic structures on ruled $4$-manifolds is made up
to diffeomorphisms which could be topologically non-trivial.  In this paper we
introduce an invariant of homotopically trivial diffeomorphisms $F:X\to X$,
called \slsf{secondary Stiefel-Whitney class} $\tiw_2(F)$, and show that
homotopically trivial symplectomorphisms $F$ with $\tiw_2(F)=0$ are smoothly
isotopically trivial. In general, the class $\tiw_2(F)\in \sfh^1(Z,\zz_2)$ is
defined for any diffeomorphism $F:Z\to Z$ of an arbitrary smooth manifold $Z$
such that $F$ acts trivially on all $\zz_2$-homology groups
$\sfh_\bullet(Z,\zz_2)$. The first main result of the paper is as follows.

\state Theorem 1. {\it Let $(X,\omega)$ be a 
 ruled or rational symplectic manifold.
\begin{itemize}
\item[\sli] For every $w\in\sfh^1(X,\zz_2)$ there exists a symplectomorphism
 $F:(X,\omega)\to(X,\omega)$ with $\tiw_2(F)=w$ and such that $F$ acts trivially on
 $\pi_1(X)$, $\pi_2(X)$, and\/ $\sfh_2(X,\zz)$.
\item[\slii] Let $F:(X,\omega)\to(X,\omega)$ be a homotopically trivial 
symplectomorphism such that $\tiw_2(F)=0$. Then $F$
 is smoothly isotopic to the identity.
\end{itemize} 
}

The part \sli is due to Denis Auroux, and I am indebted to him for
pointing me out this fact. Notice that $\sfh^1(X,\zz_2)=0$ and $\tiw_2(F)=0$ for
rational $X$.

\state Corollary 2. {\it \sli The smooth isotopy class of a symplectomorphism
 of a \emph{rational} $4$-manifold is determined by its action on the homology
 group $\sfh_2(X,\zz)$.

\slii For a rational or ruled $4$-manifold  $X$ there exists a natural
 bijection between the set of deformation classes of symplectic structures on
 $X$ with fixed cohomology class and the group $\Gamma_\bsq$ of isotopy classes of
 homotopically trivial diffeomorphisms $F:X\to X$ with  $\tiw_2(F)=0$.}

The corollary is a part of \slsf{Theorem 3} below giving a description of the
structure of the diffeotopy group of rational and ruled symplectic
$4$-manifolds. 

\smallskip%
In view of \slsf{Corollary 2} the complete classification of deformation
classes of symplectic structures $\omega$ on a given rational or ruled $4$-manifold
$X$ representing a given cohomology class reduces to an intriguing question
about the structure of the group $\Gamma_\bsq$ (``black box'') which is defined in
purely differential-topological terms. In this connection we notice that in
the case of \emph{homeomorphisms} of $S^2{\times}S^2$ the corresponding group is
trivial, see \cite{Qui} and \cite{Per}.

\tableofcontents

\newsubsection[diffeo-grp]{Diffeotopy groups of rational and ruled
 $4$-manifolds.}  Using \slsf{Theorem 1} we give a complete description of the
structure of the diffeotopy group ($\Leftrightarrow$ mapping class group) of rational and
ruled symplectic $4$-manifolds. For this description we need to introduce
various (sub)groups related to the problem.

\newdefi{def-Gamma} For a smooth manifold $X$, denote by $\diff(X)$ the group
of diffeomorphisms, by $\diff_0(X)$ its component of the identity, and by
$[X,X]$ the monoid of homotopy classes of continuous maps $F:X\to X$.

For an orientable manifold $X$ let  $\diff_+(X)<\diff(X)$ be the
subgroup which preserves the orientation of $X$.

Notice that the intersection form of\/ $\sfh_2(X,\rr)$ of a rational or ruled
$4$-manifold has Lorentz signature, and thus the set of classes
$A\in\sfh_2(X,\rr)$ with $A^2>0$ consists of two convex cones. Given a symplectic
form $\omega$ on such $X$, we call \slsf{positive} the cone containing the class
$[\omega]$ and \slsf{negative} the opposite one. Denote by $\diff_\lor(X)<\diff_+(X)$
be the subgroup which preserves the positive cone.   

The group $\diff_\bullet(X)$ consists of diffeomorphisms $f\in\diff(X)$ which are
\slsf{homotopically trivial}.  The latter means the following $2$ conditions:
\begin{itemize}
\item[(i)] $f_*:\sfh_2(X,\zz)\to\sfh_2(X,\zz)$ is the identical homomorphism;
\item[(ii)] $f$ is homotopic to a map $f'$ which preserves some base point
 $x_0$ and such that the induced homomorphisms $f'_*:\pi_1(X,x_0)\to\pi_1(X,x_0)$
 and $f'_*:\pi_2(X,x_0)\to\pi_2(X,x_0)$ are identical.
\end{itemize} 
Further, set $\diff_\bsq :=\ker\big(\tiw_2:\diff_\bullet\to\sfh^1(X,\zz_2)\big)$.

Now we define the diffeotopy group and its quotient- and subgroups. The full
\slsf{diffeotopy group} is group of isotopy classes of diffeomorphisms,
$\Gamma=\Gamma(X):=\diff(X)/\diff_0(X)$. Further, we set\/
$\Gamma_\bullet:=\diff_\bullet(X)\big/\diff_0(X)$ --- these are homotopically trivial
diffeomorphisms, and\/ $\Gamma_\bsq:= \diff_\bsq(X)\big/\diff_0(X)$, the ``black
box'' group introduced above. In particular,
$\Gamma_\bsq=\ker\big(\tiw_2:\Gamma_\bullet\to\sfh^1(X,\zz_2)\big)$.  Denote by $\wt\Gamma_H$ and
respectively by $\wt\Gamma_\sfw$ the images of $\diff(X)$ in $[X,X]$ and in
$\aut(\sfh_2(X,\zz))$, similarly, by $\Gamma_H$ and respectively by $\Gamma_\sfw$ we
denote the image of $\diff_\lor(X)$ in $[X,X]$ and in $\aut(\sfh_2(X,\zz))$. \rm
Thus $\wt\Gamma_H$ is the group of diffeomorphisms modulo homotopy, $\Gamma_H<\wt\Gamma_H$
the subgroup preserving orientation and the positive cone, and
$\Gamma_\bullet=\ker\big(\Gamma\to\wt\Gamma_H\big)$.
\end{defi}

The group $\barr\map(Y,\ell)$ is introduced below in \refdefi{def-barMap}.

\smallskip%
For a closed surface $Y$ we denote by $Y{\tix}S^2$ a closed $4$-manifold $X$
with $w_2(X)\neq0$ which admits an $S^2$-fibre bundle structure over the base
$Y$, given by a projection $\pr:X\to Y$. We call
$Y{\tix}S^2$ the \slsf{skew-product} of the surface $Y$ with $S^2$.

\state Theorem 3. {\it Let $(X,\omega)$ be a rational or ruled symplectic manifold 
 and $\iom(X,\omega)$ the space of symplectic forms on $X$ homologous to $\omega$.

\begin{itemize}
\item[\sli] 
The group $\diff_\bsq(X)$ acts transitively on $\iom(X,\omega)$. The
 group $\Gamma_\bsq$ acts \emph{simply transitively} on the set $\pi_0\iom(X,\omega)$ of connected
 components of $\iom(X,\omega)$.

\item[\slii] The group $\Gamma_\bullet$ is a semi-direct product
 $\Gamma_\bsq{\rtimes}\sfh^1(X,\zz_2)$. In particular, the homomorphism
 $\tiw_2:\Gamma_\bullet\to\sfh^1(X,\zz_2)\cong\zz_2^{2g}$ is surjective and has kernel
 $\Gamma_\bsq$.  Moreover, for any symplectic form $\omega$ the group $\Gamma_\bullet(\omega)<\Gamma_\bullet$ of
 isotopy classes of homotopically trivial symplectomorphisms is mapped by
 $\tiw_2$ isomorphically onto $\sfh^1(X,\zz_2)$ and defines a splitting of the
 product $\Gamma_\bsq{\rtimes}\sfh^1(X,\zz_2)$.

\item[\sliii] Assume that $X$ is an $\ell$-fold blow-up of the product
 $X_0=Y{\times}S^2$ or of the skew-product $Y{\tix}S^2$ with $\ell\geq0$. Then there
 exists a subgroup $\Gamma_\sfm <\Gamma_\lor$ and a natural homo\-mor\-phism
 $\Gamma_\lor\to\barr\map(Y,\ell)$ which induces an isomorphism $\Gamma_\sfm\to\barr\map(Y,\ell)$.
 Moreover, every element of\/ $\Gamma_\sfm$ is represented by a symplectomorphism
 $F:(X,\omega)\to(X,\omega)$.

\item[\sliv] The image $\Gamma_H$ of $\diff_\lor(X)$ in $[X,X]$ is isomorphic to
 $\barr\map(Y,\ell)\rtimes\Gamma_\sfw$.  The group $\Gamma_\lor$ is a semidirect product
 $\big(\Gamma_\bullet\rtimes\barr\map(Y,\ell)\big)\rtimes\Gamma_\sfw$. 

\item[\slv] For any ruled $4$-manifold $X$ there exists a diffeomorphism
 $F_c:X\to X$ which preserves orientation and acts by $-\id$ on the 2-homology
 $\sfh_2(X,\zz)$. In particular, $F_c$ inverts the positive cone in
 $\sfh_2(X,\zz)$. Moreover, $F_c$ can be realised as a complex conjugation,
 \ie, an anti-holomorphic involution with respect to some complex structure
 $J$ induced by some structure of complex projective algebraic surface on $X$;

\item[\slvi] If a ruled $4$-manifold $X$ admits an orientation inverting
 diffeomorphism, then $X$ is diffeomorphic to $Y{\times}S^2$ or to $Y{\tix}S^2$.
\end{itemize} 
}

\smallskip%
Now we describe the group $\Gamma_\sfw$ which gives the action of $\diff_\lor(X)$ on
the second homo\-logy group.  We start with the ``stable case'' of ruled
$4$-manifolds blown-up several times. The group $\Gamma_\sfw$ has different
structure in the rational and irrational cases, and we treat them separately.

By the \slsf{twist along an exceptional sphere} $E$ in a $4$-manifold $X$ we
mean a diffeomorphism supported in some small neighbourhood of $E$ and
inverting the orientation this sphere, see \refsubsection{non-symp} and
\lemma{E-twist} for further details.

For the definition of \slsf{Coxeter systems, Coxeter-Dynkin diagrams/graphs},
and \slsf{Coxeter-Weyl-} and \slsf{(quasi)reflection groups} we refer to
\refsubsection{cox}. Here we only notice that the diagrams $\rmd_\ell,\rme_\ell$
have usual meaning and the groups $\sfw(\rmd_\ell),\sfw(\rme_\ell)$ are usual Weyl
groups of the corresponding type, whereas $\wt\rma_\ell$ and $\wt\rmb_\ell$ are the
diagrams of \slsf{affine systems} of types $\rma_\ell$ and $\rmb_\ell$,
respectively, and the corresponding groups $\sfw(\wt\rma_\ell)$ and
$\sfw(\wt\rmb_\ell)$ are the usual \slsf{affine Weyl groups} of type $\rma_\ell$ and $\rmb_\ell$.

\state Theorem 4'. {\it Let $X$ be $\cp^2$ blown up at $\ell\geq3$ points,
 $\bfl,\bfe_1,\ldots,\bfe_\ell\in\sfh_2(X,\zz)$ the homology classes of the line and
 respectively of the exceptional rational curves $E_1,\ldots,E_\ell$.
\begin{itemize}
\item[\sli] The group $\Gamma_\sfw<\aut(\sfh_2(X,\zz))$ is a Coxeter-Weyl group
 $\sfw(\cals(X))$ with the Coxeter graph $\cals(X)$\\[-29pt]\
\begin{align} 
&\xymatrix@C-10pt@R-13pt{
& & & s_0\ar@{-}[d]
\\
\rmbe_{\ell+1}\ (\ell\geq4): 
&  s_1 \ar@<2pt>@{-}[r] & s_2\ar@<2pt>@{-}[r] & s_3 \ar@<2pt>@{-}[r] &
\ldots \ar@<2pt>@{-}[r] & s_{\ell-2} \ar@<2pt>@{-}[r] & s_{\ell-1}\ar@<1pt>@{=}[r] & s_{\ell}
}
\\
\noalign{\ \ \;
in the case $\ell\geq4$ or respectively}
&\xymatrix@C-10pt@R-13pt{
\rml_4(3,4,4)\ (\ell=3): 
&  s_1 \ar@<2pt>@{-}[r] & s_2\ar@<1pt>@{=}[r] & 
s_{3}\ar@<1pt>@{=}[r] & s_{0}
}
\end{align}
\ \\[-17pt]in the case $\ell=3$. 

The generators $s_i$ are represented by the Dehn twists along $(-2)$-spheres in
the homology classes $\bfs_{i,i+1}:=\bfe_{i}-\bfe_{i+1}$ for the with
$i=1,\ldots,\ell-1$, by the Dehn twist along $(-2)$-sphere in the homology class
$\bfs'_{123}:=\bfl-(\bfe_1+\bfe_2+\bfe_3)$ in the case of the generator $s_0$,
and the twist along the exceptional sphere $E_\ell$ in the case of the generator
$s_\ell$.

\item[\slii] Let $\omega$ be a symplectic form on $X$ positive on the class $\bfl$.
 Then there exists a diffeomorphism $F\in\diff_+(X)$ which is a product of the
 generators as in Part \slip, such that the form $\omega^*:=F_*\omega$ has the following
 properties:
\begin{itemize}
\item the classes $\bfl,\bfe_1,\ldots,\bfe_\ell$ are represented by $\omega^*$-symplectic
 spheres;
\item the periods  $\lambda^*:=\int_\bfl\omega^*, \mu^*_i:=\int_{\bfe_i}\omega^*$ of the form
 $\omega^*$ satisfy the \slsf{period conditions}
\begin{equation}\eqqno(period-cond-E)
\lambda^*\geq\mu^*_1+\mu^*_2+\mu^*_3,\ \ \mu^*_1\geq\mu^*_2\geq\cdots\geq\mu^*_\ell>0;
\end{equation}
\item Let $\cals(\omega^*)$ be the set of those classes
 $\bfs_{123},\bfs_{1,2},\bfs_{2,3},\ldots,\bfs_{\ell-1,\ell}$ for which
 $\int_{\bfs_{\ldots}}\omega^*=0$. Then every class $\bfs_{\ldots}$ from $\cals(\omega^*)$ is
 represented by an $\omega^*$-Lagrangian sphere $\Sigma_{\ldots}$. Furthermore, the image
 $\Gamma_\sfw(\omega^*)$ of the symplectomorphism group $\symp(X,\omega^*)$ in $\Gamma_\sfw$ is a
 finite reflection group generated by symplectic Dehn twists along Lagrangian
 spheres $\Sigma_{\ldots}$ whose homology classes $[\Sigma_{\ldots}]=\bfs_{\ldots}$ lie in $\cals(\omega^*)$.
 In particular, $\Gamma_\sfw(\omega^*)$ is a reflection group with the Coxeter system
 $\cals(\omega^*)$.
\end{itemize} 

\item[\sliii] The system $\cals(\omega^*)$ is contained in one of the following
 \slsf{maximal systems} $\cals$ listed below. In particular, the group
 $\Gamma_\sfw(\omega^*)$ is a subgroup the a Coxeter-Weyl group $\sfw(\cals)$. The
 maximal system $\cals$ is
\begin{itemize}
\item $\rme_{\ell}=
 \{\bfs'_{123},\bfs_{1,2},\bfs_{2,3},\bfs_{3,4},\ldots,\bfs_{\ell-1,\ell}\}$ \ in the case $3\leq\ell\leq8$, 
\end{itemize}
and in the case $\ell\geq9$ one of the following:
\begin{itemize}
\item $\rma_{\ell-1}=\{ \bfs_{1,2},\bfs_{2,3},\bfs_{3,4},\ldots,\bfs_{\ell-1,\ell}\};$
\item $\rmd_{\ell-1}=\{ \bfs'_{123},\bfs_{2,3},\bfs_{3,4},\ldots,\bfs_{\ell-1,\ell}\};$
\item $\rma_1\oplus\rma_{\ell-2}=\{\bfs_{1,2}\}\sqcup
\{\bfs'_{123},\bfs_{3,4},\bfs_{4,5},\ldots,\bfs_{\ell-1,\ell}\};$
\item $\rme_k\oplus\rma_{\ell-k-1}=\{
 \bfs'_{123},\bfs_{1,2},\bfs_{2,3},\ldots,\bfs_{k-1,k}\}\sqcup
 \{\bfs_{k+1,k+2},\bfs_{k+2,k+3},\ldots,\bfs_{\ell-1,\ell}\}$ %
 with some $3\leq k\leq8$.
\end{itemize} 

\item[\sliv] Let $\omega^*$ be a symplectic form on $X$ which satisfies the period
 conditions \eqqref(period-cond-E), and  $\Gamma_\Omega(\omega^*)<\Gamma_\lor$ the group
 generated by all $\omega'$-symplectomorphisms with respect to all forms $\omega'$
 deformationally equivalent to $\omega^*$. Then $\Gamma_\Omega(\omega^*)$ is the Coxeter-Weyl
 subgroup $\sfw(\cals'(X))$ of\/ $\Gamma_\sfw$ corresponding to Coxeter system
 $\cals'(X)\subset\cals(X)$ of type $\rme_\ell$ consisting of generators
 $s_0,\ldots,s_{\ell-1}$.

\end{itemize}
}

\state Remarks.~1. In the ``non-del-Pezzo'' case $\ell>8$ the group $\Gamma_\sfw(\omega^*)$
is always strictly smaller than $\Gamma_\Omega(\omega^*)$. In contrary, in the ``del-Pezzo''
case $\ell\leq8$ one has the equality $\Gamma_\sfw(\omega^*)=\Gamma_\Omega(\omega^*)$ provided there are only
equalities in the period conditions, \ie, in the case
$\frac13\lambda^*=\mu^*_1=\cdots=\mu^*_\ell$. \\[2pt]
{\bf2.} In the case $\ell\leq9$ the group $\Gamma_\sfw{\cdot}\{\pm\id\}$ is the whole
automorphism group $\aut(\sfh_2(X,\zz))=\bfo(1,\ell,\zz)$ of the lattice
$\sfh_2(X,\zz)$. However, in the case  $\ell\geq10$ there exist automorphisms of 
$\sfh_2(X,\zz)$ which can not be realised by diffeomorphisms. See more details
in Remark on the page \pageref{remark-SW}.

\medskip%
\state Theorem 4''. {\it Let $X_0$ be the product $Y{\times}S^2$ of the sphere with a
 surface of genus $g(Y)\geq1$, $X$ a blow-up of $X_0$ at $\ell\geq2$ points,
 $\bfy,\bff;\bfe_1,\ldots,\bfe_\ell\in\sfh_2(X,\zz_2)$ the homology classes of
 respectively a horizontal section $Y{\times}\{s_0\}$, a generic vertical fibre
 $\{y_0\}{\times}S^2$ and the exceptional spheres $E_1,\ldots,E_\ell$.
\begin{itemize}
\item[\sli] The group $\Gamma_\sfw<\aut(\sfh_2(X,\zz))$ is a Coxeter-Weyl\/ group 
 $\sfw(\cals(X))$ with the Coxeter graph $\cals(X)$\\[-29pt]\
\begin{align} 
&\xymatrix@C-10pt@R-13pt{
& & s_0\ar@{-}[d]
\\
\rmbd_{\ell+1}=\wt \rmb_\ell\ (\ell\geq3): 
&  s_1 \ar@<2pt>@{-}[r] & s_2\ar@<2pt>@{-}[r] & 
\ldots \ar@<2pt>@{-}[r] & s_{\ell-2} \ar@<2pt>@{-}[r] & s_{\ell-1}\ar@<1pt>@{=}[r] & s_{\ell}
}
\\
\noalign{\ \ \;
in the case $\ell\geq3$ or respectively}
&\xymatrix@C-10pt@R-13pt{
\rml_3(4,4)=\wt \rmb_2\ (\ell=2): 
& s_1\ar@<1pt>@{=}[r] & 
s_{2}\ar@<1pt>@{=}[r] & s_{0}
}
\end{align}
\ \\[-17pt]in the case $\ell=2$.

The generators $s_i$ are represented by Dehn twists along $(-2)$-spheres in the
homo\-logy classes $\bfs_{i,i+1}:=\bfe_{i}-\bfe_{i+1}$ in the case $i=1,\ldots,\ell-1$,
by the Dehn twist along $(-2)$-sphere in the homology class
$\bfs'_{12}:=\bff-(\bfe_1+\bfe_2)$ in the case of the generator $s_0$, and the
twist along the exceptional sphere $E_\ell$ in the case of the generator $s_\ell$.

\item[\slii] Let $\omega$ be a symplectic form on $X$ positive on the class $\bff$.
 Then there exists a diffeomorphism $F\in\diff_+(X)$ which is a product of the
 generators as in Part \slip, such that the form $\omega^*:=F_*\omega$ has the following
 properties:
\begin{itemize}
\item the classes $\bff,\bfy,\bfe_1,\ldots,\bfe_\ell$ are represented by
 $\omega^*$-symplectic surfaces;
\item the periods  $\varsigma^*:=\int_\bff\omega^*, \mu^*_i:=\int_{\bfe_i}\omega^*$ of the form
 $\omega^*$ satisfy the \slsf{period conditions}
\begin{equation}\eqqno(period-cond-D)
\varsigma^*\geq\mu^*_1+\mu^*_2,\ \ \mu^*_1\geq\mu^*_2\geq\cdots\geq\mu^*_\ell>0;
\end{equation}
\item Let $\cals(\omega^*)$ be the subset of those classes
 $\bfs'_{12},\bfs_{1,2},\bfs_{2,3},\ldots,\bfs_{\ell-1,\ell}$ for which
 $\int_{\bfs_{\ldots}}\omega^*=0$. Then every class $\bfs_{\ldots}$ from $\cals(\omega^*)$ is
 represented by $\omega^*$-Lagrangian sphere $\Sigma_{\ldots}$. Furthermore, the image
 $\Gamma_\sfw(\omega^*)$ of the symplectomorphism group $\symp(X,\omega^*)$ in $\Gamma_\sfw$ is a
 finite reflection group generated by symplectic Dehn twists along Lagrangian
 spheres $\Sigma_{\ldots}$ whose homology class $[\Sigma_{\ldots}]=\bfs_{\ldots}$ lies in $\cals(\omega^*)$.
 In particular, $\Gamma_\sfw(\omega^*)$ is a reflection group $\sfw(\cals(\omega^*))$ with
 the Coxeter system $\cals(\omega^*)$.
\end{itemize} 

\item[\sliii] Let $\omega^*$ be a symplectic form on $X$ which satisfies the period
 conditions \eqqref(period-cond-D), and let $\Gamma_\Omega(\omega^*)<\Gamma_\lor$ be a the group
 generated by all $\omega$-symplectomorphism with respect to all forms $\omega$
 deformationally equivalent to $\omega^*$. Then $\Gamma_\Omega(\omega^*)$ is the Coxeter-Weyl\/
 $\sfw(\cals'(X))$ subgroup of\/ $\Gamma_\sfw$ corresponding to Coxeter system
 $\cals'(X)\subset\cals(X)$ of type $\rmd_\ell$ consisting of generators
 $s_0,\ldots,s_{\ell-1}$. In particular, $\Gamma_\sfw(\omega^*)$ is the whole group $\Gamma_\Omega(\omega^*)$
 if and only if one has equality in all period
 conditions \eqqref(period-cond-D) except the strict inequality $\mu^*_\ell>0$.
\end{itemize}
}

\medskip%
Now we describe the group $\Gamma_\sfw$ or $\wt\Gamma_\sfw$ in the remaining special
cases $X=\cp^2$, $X=Y{\times}S^2$, $X=Y{\tix}S^2$, and $X=(Y{\times}S^2)\#\barr\cp^2$
where $Y$ is some closed oriented surface of genus $g\geq0$.  Recall (see
\refsubsection{ruled}) that the manifolds $X=Y{\times}S^2$ and $X=Y{\tix}S^2$ admit
non-singular rulings, \ie, $S^2$-bundle structure $\pr:X\to Y$. We denote by
$\bff$ the homology class of the fibre the rulings above, by $\bfy_0$ the
homology class of the horizontal section $Y{\times}\pt$ of the projection
$\pr:X=Y{\times}S^2\to Y$, so that self-intersection of\/ $\bfy_0$ is $0$, and by
$\bfy_{\pm1}$ disjoint sections of self-intersection $\pm1$ (case $X=Y{\tix}S^2$).
Notice that $\bfy_1=\bfy_{-1}+\bff$. Notice also that we can realise
$X=Y{\tix}S^2$ as projective completion of the total space of a holomorphic
line bundle $\scrl$ on $Y$ equipped with some complex structure $J_Y$ with
$c_1(\scrl)=1$. In this realisation the classes $\bfy_{\pm1}$ are represented by
the zero and infinity sections of $\scrl$. Further, recall that
$X=S^2{\tix}S^2$ can be described also as $\cp^2$ blown-up at one point, such
that $\bfy_{-1}$ and respectively $\bfy_{+1}$ become the homology classes of
the exceptional curve and of a generic line. In this case
$X=S^2{\tix}S^2$ we denote the classes by $\bfe=\bfy_{-1}$ and
$\bfl=\bfy_{+1}$. Furthermore, notice that $X=(S^2{\times}S^2)\#\barr\cp^2$ is the
$\cp^2$ blown up twice, and in this case we use the notation
$\bfl,\bfe_1,\bfe_2$ for the homology classes of a complex line and the
arising exceptional rational curves. Finally, in the case
$X=(Y{\times}S^2)\#\barr\cp^2$ the notation $\bfy_0$ and $\bff$ has the same meaning
as in the case $X=Y{\times}S^2$, and $\bfe_1$ denotes classes of the exceptional
sphere.

\medskip%
\state Theorem 4'{'}'. {\it \sli The group $\wt\Gamma_H(\cp^2)$ is $\zz_2$
 generated by the complex conjugation.  \\[4pt]
 \slii In the case $X=S^2{\times}S^2=\cp^1{\times}\cp^1$ the group $\wt\Gamma_H$ is
 $(\zz_2\oplus\zz_2)\rtimes\zz_2$ in which the first two factors $\zz_2$ are
 generated by the complex conjugations of the corresponding factors $\cp^1$
 and the last factor $\zz_2$ is generated by the holomorphic involution
 $(z,w)\mapsto(w,z)$ interchanging the
 factors of  $X=S^2{\times}S^2$. \\[4pt]
 \sliii In the case $X=S^2{\tix}S^2$ the group $\wt\Gamma_H$ is
 $(\zz_2\oplus\zz_2)\rtimes\zz_2$ in which the first two factors $\zz_2$ act on
 $\sfh_2(X,\zz)=\zz\lan\bfl,\bfe\ran$ by inverting the classes $\bfl$ and
 respectively $\bfe$, and the last factor $\zz_2$ is represented by an
 orientation inverting involution which interchange the classes $\bfl$ and
 $\bfe$. \\[4pt]
 \sliv In the case $X=Y{\times}S^2$ with $g(Y)\geq1$ the group $\wt\Gamma_\sfw$ is the
 group $\zz_2\oplus\zz_2$ in which factors $\zz_2$ act on
 $\sfh_2(X,\zz)=\zz\lan\bfy_0,\bff\ran$ by inverting the classes $\bfy_0$ and
 respectively $\bff$, and are represented by anti-holomorphic involutions
 (``complex conjugations'') of the single factors $Y$ and respectively $S^2$.
 \\[4pt]
 \slv In the case $X=Y{\tix}S^2$ with $g(Y)\geq1$ the group $\wt\Gamma_\sfw$ is the
 group $\zz_2\oplus\zz_2$ in which the first factors $\zz_2$ act on
 $\sfh_2(X,\zz)=\zz\lan\bfy_1,\bff\ran$ by inverting both the classes $\bfy_1$
 and $\bff$ and is represented by anti-holomorphic involution (``complex
 conjugation'') of $X$ with respect to some complex structure, and the other
 factor $\zz_2$ can be represented by an involution which inverts the class
 $\bff$ and interchanges the classes $\bfy_{+1}$ and $\bfy_{-1}$. In
 particular, the latter involution inverts the orientation of $X$.
 \\[4pt]
 \slvi In the case $X=(S^2{\times}S^2)\#\barr\cp^2$ the group $\Gamma_H=\Gamma_\sfw$ is the
 integer orthogonal group $\bfo_+(1,2,\zz)$ of the lattice
 $\sfh_2(X,\zz)=\zz\lan \bfl,\bfe_1,\bfe_2\ran$ %
 preserving the positive cone in $\sfh_2(X,\rr)$.  It is isomorphic to the
 Coxeter-Weyl group of the type $\rml_3(4,\infty)$ with the Coxeter graph
 $\xymatrix@C-10pt@R-13pt{ s_1\ar@<1pt>@{=}[r] & s_2\ar@<1pt>@{-}[r]^\infty &
  s^*_0}$ in which $s_1$, $s_2$, and $s_0^*$ are the reflections with respect
 to the classes $\bfs_{1,2}:=\bfe_1-\bfe_2$, $\bfe_2$, and
 $\bfe'_{12}:=\bfl-\bfe_1-\bfe_2$, respectively.
 \\[4pt]
 \slvii In the case $X=(Y{\times}S^2)\#\barr\cp^2$ with $g(Y)\geq1$ the group
 $\wt\Gamma_\sfw$ is $\zz_2\ltimes\zz$ realised as the Coxeter-Weyl group of type
 $\rmi_2(\infty)=\wt\rma_1$ with the Coxeter graph $\xymatrix@C-10pt@R-13pt{
  s_1\ar@<1pt>@{-}[r]^\infty & s_1^*}$ such that the generators $s_1$ and $s^*_1$
 are the reflections with respect to the classes $\bfe_1$ and
 $\bfe'_1:=\bff-\bfe_1$. }\rm


\state Remark. In each case above we have a Coxeter-Weyl group. In
particular, $\zz_2=\sfw(\rma_1)$, $(\zz_2\oplus\zz_2)\rtimes\zz_2=\sfw(\rmb_2)$,
and $\zz_2\oplus\zz_2=\sfw(\rmd_2)=\sfw(\rma_1\oplus\rma_1)$.

\newdefi{def-barMap} Let $Y$ be a closed connected oriented surface and
$\mbfy=\{y^*_1,\ldots,y^*_\ell\}$ a collection of pairwise distinct points on $Y$. In
this paper typically $Y$ will be a base of some ruling $\pr:X\to Y$ (see
\refdefi{def-ruled}), so that $\pi_1(X)\cong\pi_1(Y)$, and $\mbfy$ will be singular
values of\/ $\pr:X\to Y$.  Denote by $\diff_+(Y,\mbfy)$ the group of
diffeomorphisms of\/ $Y$ preserving the orientation and the points $\mbfy$, by
$\diff_0(Y,\mbfy)$ its component of the identity, and by
$\map(Y,\ell):=\diff_+(Y,\mbfy)/\diff_0(Y,\mbfy)$ the \slsf{mapping class group}
of the surface $Y$ with $\ell$ marked points. 

Denote the group $\map(Y,0)$ by $\map(Y)$.  Then there exists
a natural forgetful homomorphism $\map(Y,\ell)\to\map(Y)$. In the case $g(Y)\geq2$ the
kernel of\/ $\map(Y,\ell)\to\map(Y)$ is the \slsf{pure braid groups} $\pbr_\ell(Y)$ on
$\ell$ strands of the surface $Y$.  Fix a base point $y^*$ on $Y$ and paths
$\gamma_1,\ldots,\gamma_\ell$ connecting $y^*_i$ with $y^*$.  There exists a natural surjective
homomorphism $\pbr_\ell(Y)\to(\sfh_1(Y,\zz))^\ell$ (direct product of groups
$\sfh_1(Y,\zz)$). Further, consider the quotient map
$(\sfh_1(Y,\zz))^\ell\to(\sfh_1(Y,\zz))^\ell/\sfh_1(Y,\zz)$ where $\sfh_1(Y,\zz)$ is embedded in
$(\sfh_1(Y,\zz))^\ell$ in the diagonal way, and the composition
$\pbr_\ell(Y)\to(\sfh_1(Y,\zz))^\ell/\sfh_1(Y,\zz)$. The kernel of the latter homomorphism
is normal in $\map(Y,\ell)$.  Denote this kernel by $\sfk(Y,\ell)$ and define
$\barr\map(Y,\ell)$ as the quotient of\/ $\map(Y,\ell)$ by $\sfk(Y,\ell)$.
Thus $\barr\map(Y,\ell)$ is a group extension of\/ $\map(Y)$ with
the kernel $(\sfh_1(Y,\zz))^{\ell-1}$. 

In the case $g=1$ when $Y$ is the torus $T^2$ the situation changes as
follows. Here $\pi_1(T^2)=\sfh_1(T^2,\zz)$. Further, the kernel of\/
$\map(T^2,\ell)\to\map(T^2)$ is the quotient
$\pbr_\ell(T^2)\big/\pi_1(T^2)=\pbr_\ell(T^2)\big/\sfh_1(T^2,\zz)$ of the pure braid
group.  and we obtain the homomorphism
$\pbr_\ell(T^2)\big/\sfh_1(T^2,\zz)\to(\sfh_1(T^2,\zz))^\ell\big/\sfh_1(T^2,\zz)$. We define
$\barr\map(T^2,\ell)$ as the quotient of\/ $\map(T^2,\ell)$ by the kernel of the
above homomorphism. Thus $\barr\map(T^2,\ell)$ is also a group extension of\/
$\map(T^2)$ with the kernel $(\sfh_1(T^2,\zz))^\ell/\sfh_1(T^2,\zz)\cong(\sfh_1(T^2,\zz))^{\ell-1}$.
\end{defi}

\state Acknowledgements.  I want to express my gratitude to Denis Auroux who
communicated me an example of a symplectomorphism as in the part \sli of
\slsf{Theorem 1} showing a mistake in the first version of the paper. 

Also I would like to thank Stefan Nemirovski who helped me to understand the
Seiberg-Witten theory and to find the correct form of the genus inequality. 

\newsection[preli]{Definitions and preliminary results.}

\newsubsection[j-curves]{Pseudoholomorphic
 curves in $4$-manifolds} %
In this subsection we simply list some facts from
the theory of  pseudoholomorphic curves and 
refer to \cite{Iv-Sh-1,Iv-Sh-2,Iv-Sh-3,McD-Sa-1} for the proofs. 

Here and below ``smoothness'' means some  $C^{k,\alpha}$-smoothness
with  $0<\alpha<1$ and $k\in\nn$ sufficiently large.

Consider the following situation. Let $X$ be a $4$-manifold. Denote by $\scrj$
the space of all smooth almost complex structures on $X$ tamed by \emph{some}
symplectic form. Then $\scrj$ becomes a Banach manifold.

Let $\scrj_0$ be some fixed connected component of $\scrj$.  Fix also an
integer homology class $[A]\in\sfh_2(X,\zz)$ and an integer $k\geq0$. Denote by
$X^{[k]}\subset X^k=X\times\ldots\times X$ the space of $k$-tuples $\mbfx=(x_1,\ldots,x_k)$ of pairwise
distinct points on $X$, and by $\scrm=\scrm(X,[A],k)$ the space of tuples
$(J,C;\mbfx)$ where $J\in\scrj_0$, $C\subset X$ an immersed irreducible non-multiple
$J$-holomorphic curve of genus $g$, and $\mbfx=(x_1,\ldots,x_k)$ a collection of
pairwise distinct non-singular points on $C$. In the case $k=0$ we use an
abbreviated notation. The space $\scrm$ is a smooth Banach manifold and the
natural projection $\pr:\scrm\to\scrj_0\times X^{[k]}$ is a smooth Fredholm map of
$\rr$-index $\ind(\pr)=2\cdot(c_1(X)\cdot[A]+g-1-k)$ where $c_1(X)$ is the first Chern
class of structures $J\in\scrj_0$. We denote
$\scrm_{J,\mbfx}:=\pr\inv(J,\mbfx)$. In the case of a path $h:I\to\scrj$ with
$J_t:=h(t)$ we denote
\[
\scrm_{h,\mbfx}:=\{\,(C,\mbfx,t):t\in I,J_t:=h(t), 
(J_t,C,\mbfx)\in\scrm_{J_t,\mbfx} \,\}.
\]

\newlemma{lem-S2-J} \sli In the case $\ind(\pr)\leq0$ the map $\pr:\scrm\to\scrj_0\times X^{[k]}$
is a smooth immersion of codimension $-\ind(\pr)$, and an imbedding in the case
when pseudoholomorphic curves $C$ in $\scrm$ are imbedded and $[C]^2<k$.

\slii Assume that $\ind(\pr)\geq0$ and $\scrm$ is non-empty. Then there exists a
subset $\scrz\subset\scrj\times X^{[k]}$ of second Baire category
such that for each $(J,\mbfx)\in\scrz$ the set $\scrm_{J,\mbfx}$  is a
non-empty finite dimensional manifold. Moreover, for a generic path
$h:I\to\scrj$ the space $\scrm_{h,\mbfx}$ is a smooth manifolds of dimension
$\ind(\pr)+1$.

\sliii Assume that curves $C$ in $\scrm$ are rational, $C\cong S^2$, and
$\ind(\pr)\geq0$. Then the set $\scrz$ above is \emph{connected}. 
\end{lem} 

\proof These are standard facts of Gromov's theory. We make only two remarks.
Firsts, part \sliii is proved in \cite{H-L-S,Iv-Sh-3}.  Second remark is that
the local intersection index of $J$-holomorphic curves at any intersection
point is strictly positive. Thus if two $J$-curves $C,C'$ are homologous and
pass through $x_1,\ldots,x_k$, then $[C]^2\geq k$. This explains why in the second part
of \sli for each $J$ one has at most one $J$-holomorphic curve.
\qed

\newcorol{generic} \sli Assume that the class $[A]\in\sfh_2(X,\zz)$ is
represented by an embedded symplectic $(-1)$-sphere. Then for generic
$J\in\scrj_0$ there exists a unique exceptional $J$-holomorphic curve
representing $[A]$.

\slii Assume that the class $[A]\in\sfh_2(X,\zz)$ is represented by an embedded
symplectic sphere and $[A]^2=0$. Then for any $J\in\scrj_0$ and any point $x\in X$
there exists a unique $J$-holomorphic curve $C$ (possibly singular and
reducible) which represents the class $[A]$ and passes through $x$. Moreover,
there exists a unique \emph{continuous} fibration $\pr:X\to Y$ such that $C$ is
its fibre $\pi\inv(y)$.  Furthermore, if $J\in\scrj_0$ is generic, then $C$ is
either smooth or consists of two exceptional spheres meeting transversally in
a single point.
\end{corol}

\newsubsection[ruled]{Ruled symplectic $4$-manifold} %
For general information of ruled  $4$-manifolds we refer to \cite{Go-St}.
For the purpose of this paper we need a generalisation the standard definition
of rulings on $4$-manifolds admitting singular fibres.

\newdefi{def-ruled} A \slsf{holomorphic ruling} on a complex surface $X$ with
the complex structure $J$ is a proper holomorphic projection $\pr:X\to Y$ on a
complex curve $Y$ for which  there exists a holomorphic bundle
$\barr\pr:\barr{X}\to Y$ with the fibre $\cp^1$ such that $X$ is obtained by a
multiple blow-up and the projection $\pr:X\to Y$ is the composition of the
projection $\barr\pr:\barr{X}\to Y$ with the natural map $p:X\to\barr{X}$.

A \slsf{ruling} on symplectic $4$-manifold $(X,\omega)$ is
given by a proper projection $\pr:~\!\!X \to Y$ on an oriented surface $Y$ such
that 
\begin{itemize}
\item there exist finitely many \slsf{singular values} $y^*_1,\ldots,y^*_n$ such that
$\pr:X \to Y$ is a spherical fibre bundle over $Y^\circ:=Y\bs\{y^*_1,\ldots,y^*_n\}$;
\item for every singular value $y^*_i$ there exists a disc neighbourhood
 $U_i\subset Y$ of $y^*_i$ and an integrable complex structures $J$ in $\pr\inv(U_i)$
 and $J_Y$ in $U_i$ such that the restricted projection $\pr:\pr\inv(U_i)\to U_i$
 is a holomorphic ruling;
\item every fibre $X_y := \pr\inv(y)$ is a union $\cup_j S_j$ of spheres, called
 \slsf{components} of the fibre, and $\omega$ is positive on each component $S_j$
 of a fibre.
\end{itemize}

A fibre $X_y = \pr\inv(y)$ is \slsf{regular} if it has a single component and
\slsf{ordinary singular} if it consists of two exceptional spheres.

{\rm The usual definition of the ruling  does not allow singular fibres. Ruling
with only ordinary singular  fibres are a special case of \slsf{Lefschetz
 fibrations}.}

A \slsf{section} of a ruling $\pr: X\to Y$ is a surface $S \subset X$ such that the
restriction $\pr|_S: S \to Y$ is a diffeomorphism, or a map $\sigma: Y \to X$ such that
$\pr\circ \sigma =\id_Y$, with the obvious bijective correspondence between these two
notions.

We shall also consider the following more general \slsf{singular symplectic
 rulings}: the map $\pr:X\to Y$ is merely continuous, there are finitely many
points $y^*_i$ on $Y$ such that $\pr$ is smooth outside the fibres
$F^*_i:=\pr\inv(y^*_i)$, and such that there exists a \emph{homeomorphism}
$\Phi:X\to X$ such that
\begin{itemize}
\item $\Phi:X\to X$ maps every fibre (singular or not) onto itself, $\Phi$ is smooth
 in outside the singular fibres $F^*_i$ and on each component of those, and
 such that $\pr\circ\Phi$ is a smooth symplectic ruling in the above sense.
\end{itemize} 
{\rm In particular, every regular fibre $\pr\inv(y)$ and every component of a
singular fibre is $\omega$-symplectic.}

\smallskip
A symplectic $4$-manifold $(X,\omega)$ is \slsf{ruled} if it admits a symplectic
ruling, and  \slsf{rational} if $X$ is either $\cp^2$ or a ruled
$4$-manifold such that the base $Y$ is the sphere $S^2$. Non-rational ruled
$4$-manifolds are called \slsf{irrational}. A ruled $4$-manifold $(X,\omega)$ is
\slsf{minimal} if it admits a ruling without singular fibres.
{\rm Notice that $\cp^2$ is the only rational
non-ruled $4$-manifold.}
\end{defi} 

The structure of (closed) ruled symplectic $4$-manifold is well understood,
see \cite{La,La-McD,McD-2,McD-1,McD-Sa-2}. In particular, there exists a
topological characterisation of ruled symplectic $4$-manifold. 

\newprop{ruled1} Let $(X,\omega)$ be a ruled symplectic manifold.

\sli If $\omega'$ is another symplectic form, then there exists a ruling $\pr':X\to Y$
with $\omega'$-symplectic fibres.

\slii If moreover $\omega'$ is cohomologous to $\omega$, then there exists a
diffeomorphism $F:X\to X$ with $F_*\omega=\omega'$.  
\end{prop}

The crucial technique used in the proof is the Seiberg-Witten theory. It
provides the existence of symplectic sphere $F\subset(X,\omega')$ with $[F]^2=0$. The
rest follows by ``usual'' Gromov's theory: \refcorol{generic} show that for
generic $\omega'$-tamed almost complex structure $J$ there exists a unique singular
ruling $\pr_J:X\to Y$ with $J$-holomorphic fibres homologous to $F$ such that
all singular fibres are ordinary.

Unfortunately, such a ruling $\pr_J:X\to Y$ in almost never smooth in a
neighbourhood of singular fibres.  However, it can be regularised by means of
the following construction: Take a $C^0$-small perturbation $\ti J$ of $J$
which is holomorphic in a neighbourhood of singular fibres and such that
singular fibres remain $\ti J$-holomorphic. Then the Gromov compactness
theorem for $C^0$-continuous almost complex structures (see \cite{Iv-Sh-4})
imply that the new ruling $\pr_{\ti J}:X\to Y$ will be $C^0$-close to the old
one. This implies the existence of the homeomorphism $\Phi:X\to X$ in the
definition above transforming the old ruling in the new one. Now, since
$\ti{J}$ is holomorphic in a neighbourhood of the ``problem locus'' (singular
fibres), the ruling $\pr_{\ti{J}}:X\to Y$ is smooth.

\newdefi{regularise} We call the above construction \slsf{regularisation of
 the ruling}.
\end{defi}

Notice also that $\ti J$ remains generic enough. Namely, with exception of fibres,
every compact $\ti J$-holomorphic curve escapes somewhere from the set $U$
where $\ti J$ is holomorphic. Since we can deform $\ti J$ in $X\bs U$
arbitrarily, the transversality/genericity property of $\ti J$ is remained.
In particular, the results of \refcorol{generic} are valid.

\smallskip %
Finally, let us notice that every non-minimal ruled symplectic $4$-manifolds
is a multiple symplectic blow-up of some \emph{minimal} ruled complex surface
$X_0$.  Depending on the parity $w_2(X_0)$, $X_0$ is topologically either a
product $S^2{\times}Y$, or a twisted product $Y{\ti \times}S^2$ which is a topologically
non-trivial $S^2$-bundle over $Y$.

\smallskip%
The following result will be used later.

\begin{lem}\label{lagr-S2} Let $(X,\omega)$ be a ruled symplectic $4$-manifold,
 $\pr:X\to Y$ an $\omega$-symplectic ruling, and $E_1$, $E_2$ two exceptional spheres
 which are components of two disjoint singular fibres of $\pr$. Assume that
 $\int_{E_1}\omega = \int_{E_2}\omega$. Then there exists a Lagrangian sphere $S\subset(X,\omega)$
 representing the homology class $[E_1]-[E_2]$.
\end{lem} 

\proof Fix an $\omega$-tame almost complex structure $J_0$ which is generic enough
and such that $E_1,E_2$ are $J_0$-holomorphic. Then there exists a
$J_0$-holomorphic ruling $\pr_0:X\to Y$ such that the homology classes of fibres
of $\pr$ and $\pr_0$ are the same. Replacing $\pr$ by $\pr_0$ we may assume
that $\pr$ has only ordinary singular fibres. In particular, there exist
$\omega$-symplectic exceptional curves $E'_1,E'_2$ such that $E_1\cup E'_1$ and
$E_2\cup E'_2$ are two such singular fibres.  Observe that the spheres $E_1,E'_1$
are disjoint from $E_2,E'_2$. Let $p_1:=E_1\cap E'_1$ and $p_2:=E_2\cap E'_2$ be the
intersection point. Applying to $E'_1$ and to $E'_2$ appropriate
symplectomorphisms supported in small neighbourhoods of the point $p_1$ and
respectively $p_2$ we can make the tangent plane $T_{p_1}E'_1$ $\omega$-orthogonal
to $T_{p_1}E_1$ and $T_{p_2}E'_2$ $\omega$-orthogonal to $T_{p_1}E_2$. Notice that
$E_1,E'_1$ remain disjoint from $E_2,E'_2$. Applying Moser's technique we can
shows that in a neighbourhood $U_1$ of $p_1$ there exist Darboux coordinates
$\xi_1,\eta_1,\xi'_1,\eta'_1$ such that the symplectic form $\omega$ is given by $\omega=d\xi_1\land
d\eta_1+d\xi'_1\land d\eta'_1$ %
and such that $E_1\cap U_1$ and $E'_1\cap U_1$ are discs given by $\xi'_1=\eta'_1=0$ and
respectively by $\xi_1=\eta_1=0$. In other words, locally near $p_1$ the pieces
$E_1\cap U_1$ and $E'_1\cap U_1$ are coordinate discs for the complex coordinate
$\zeta_1:=\xi_1+\cpi \eta_1$ and $\zeta'_1:=\xi'_1+\cpi\eta'_1$. Let us make a similar
construction in a neighbourhood $U_2$ of the point $p_2$ getting complex
coordinates $\zeta_2:=\xi_2+\cpi \eta_2$ and $\zeta'_2:=\xi'_2+\cpi\eta'_2$ such that $\omega=d\xi_2\land
d\eta_2+d\xi'_2\land d\eta'_2$. %
Applying Moser's technique again, in a neighbourhood $V_1$ of $E_1\cup E'_1$ and a
neighbourhood $V_2$ of $E_2\cup E'_2$ we can construct local rulings $\pr:V_1\to\Delta_1$
and $\pr:V_2\to\Delta_2$ with the following properties:
\begin{itemize}
\item Both local rulings  are $\omega$-symplectic;
\item  $E_1\cup E'_1$ and $E_2\cup E'_2$ are the only singular fibres;
\item In the neighbourhoods $U_1$ and $U_2$ they are given by
 $\pr:(\zeta_i,\zeta'_i)\mapsto\zeta_i^2+\zeta'_i{}^2$ with respect to complex coordinates $z_1$ in
 $\Delta_1$ and respectively $z_2$ in $\Delta_2$.
\end{itemize} 
It follows that in the neighbourhoods $V_1,V_2$ there exist an $\omega$-tame almost
complex structure $J$ which is generic enough, such that both local rulings
$\pr:V_1\to\Delta_1$ and $\pr:V_2\to\Delta_2$ are $J$-holomorphic, and such that in the
neighbourhoods $U_1$ and $U_2$ the coordinate $\zeta_i,\zeta'_i$ are $J$-complex. Now
we can extend the local rulings $\pr:V_1\to\Delta_1$ and $\pr:V_2\to\Delta_2$ to a global
$\omega$-symplectic ruling $\pr_1:X\to Y$. In particular, $\Delta_1$ and $\Delta_2$ become
disjoint embedded discs on $Y$. Replacing $\pr$ by $\pr_1$ we may assume
that $\pr$ has the above properties.

Set $y_1:=\pr(p_1)$ and $y_2:=\pr(p_1)$. Fix an embedded path $\alpha:I\to Y$ with the
following properties: 
\begin{itemize}
\item $\alpha$ connects the points $y_1$ and $y_2$ and avoids other critical values
 of the ruling $\pr$;
\item In each disc  $\Delta_1$ and $\Delta_2$ the path $\alpha$ is the radial interval in the
 positive real direction.
\end{itemize}
Denote by $W$ the pre-image $\pr\inv(\alpha)$ without the singular fibres $E_1\cup
E'_1$ 
and $E_2\cup E'_2$. Let $h$ be a local defining function of $\alpha$ on $Y$ such
that in the discs $\Delta_1$ and $\Delta_2$ one has $h|_{\Delta_1}=\Im(z_1)$ and
$h|_{\Delta_2}=-\Im(z_2)$. Denote by $H$ the composition $H:=h\circ\pr$ and consider it
as a local Hamiltonian function on $(X,\omega)$. Denote by $\Phi^H$ the Hamiltonian
flow of $H$ on $X$, and let $\mbfv^H$ be the Hamilton vector field. Since $H$
vanishes on $W$, $\mbfv^H$ is tangent to $W$, and the restriction $\mbfv^H|_W$
is an isotropic vector field of the restriction $\omega|_W$. Further, let us
notice that the map $\pr:W\to\alpha$ is a $S^2$-bundle and that $\omega$ is non-degenerate
on each fibre of $\pr:W\to\alpha$. It follows that $\mbfv^H$ is transversal to each
fibre $\pr\inv(y)$, $y\in\alpha$.

Let $L_1$ be the disc in the neighbourhood $U_1$ given by the equations
$\eta_1=\eta'_1=0$, \ie, the disc in the coordinate plane given by the coordinates
$\xi_1,\xi'_1$. Then $L_1$ is $\omega$-Lagrangian and the projection $\pr:L_1\to Y$
takes values in $\alpha$ and there exists a parameter $\tau$ on the path $\alpha$
such that the projection $\pr:L_1\to Y$ is given by the formula
$\pr:(\xi_1,\xi'_1)\mapsto\tau=\xi_1^2+\xi'_1{}^2$. Moreover, the Hamiltonian vector field is
given by $\mbfv^H=2(\xi\partial_\xi+\xi'\partial_{\xi'})$ and is tangent to $L_1$.

Fix a point $y_0$ on $\alpha$. Consider flow lines $\ell$ of the field $\mbfv^H$
passing through points on $L_1$. As we have seen, all these lines $\ell$ are
tangent to $L_1$ and hence lie on $L_1$. Moreover, these flow lines are
radial. This means that introducing polar coordinates $(\rho_1,\varphi_1)$ on $L_1$ we
can parametrise the set of flow lines $\ell$ on $L_1$ by the angle coordinate
$\varphi_1$.  Further, as it was shown above, every flow line $\ell$ is transversal the
fibres of the ruling $\pr:X\to Y$. Hence every projection $\pr:\ell\to\alpha$ is a local
diffeomorphism and we can use the parameter $\tau$ on $\alpha$ as a parameter on $\ell$.
Fix a point $y_0$ on $\alpha$. Extend each flow line $\ell_{\varphi_1}$ on $L_1$ to the flow
line $\ell^+_{\varphi_1}$ over the interval $(y_1,y_0]$ on $\alpha$ and denote by $L^+_1$ the
union of all flow lines $\ell^+_{\varphi_1}$, $\varphi_1\in[0,2\pi]$, including the point
$p_1$. Then $L^+_1$ is a smooth disc with the centre $p_1$ such that the
punctured disc $L^+_1\bs\{p_1\}$ is embedded in $W$ and the projection
$\pr:L^+_1\to\alpha$ is a Morse function with the unique critical point $p_1$. Since
every flow line $\ell$ of $\mbfv^H$ is isotropic with respect to $\omega|_W$, the
restriction $\omega|_{L^+_1}$ vanishes. This implies that  $L^+_1$ is
$\omega$-Lagrangian. Denote by $\gamma_1:=\partial L^+_1$ the boundary circle. By the
construction, $\gamma_1$ is an embedded circle on the fibre $F_0:=\pr\in(y_0)$.

Repeat the above construction for the point $p_2$, the disc $L_2$ and the
interval $[y_0,y_2]$ on $\alpha$. We obtain an $\omega$-Lagrangian disc $L^+_2$ with the
projection $\pr:L^+_2\to[y_0,y_2]$ and the boundary circle $\gamma_2=\partial L^+_2$ lying on
$F_0$. The circle $\gamma_1$ (resp.\ $\gamma_2$) divides the sphere $F_0\cong S^2$ into two
discs which we denote by $D_1,D'_1$ (resp.\ $D_2,D'_2$). We claim that
after reindexation if needed the areas of these discs are
\begin{equation}\eqqno(Dw=Ew)\textstyle 
\int_{D_1}\omega=\int_{E_1}\omega,\quad
\int_{D'_1}\omega=\int_{E'_1}\omega,\qquad
\int_{D_2}\omega=\int_{E_2}\omega,\quad
\int_{D'_2}\omega=\int_{E'_2}\omega
\end{equation}
Indeed, for every $\tau\in(y_1,y_0]\subset\alpha$ the fibre $\pr\inv(\tau)$ is cut by $L^+_1$
into two discs.  Denote them by $D_1(\tau)$ and respectively $D'_1(\tau)$ so that
both collections become continuous families. Since $L^+_1$ is $\omega$-Lagrangian,
both functions $\int_{D_1(\tau)}\omega$ and  $\int_{D'_1(\tau)}\omega$ must be
constant. Further, as $t\in(y_1,y_0]$ tends to $y_1$ the discs $D_1(\tau)$ and
respectively $D'_1(\tau)$ converge to the spheres $E_1$ and respectively
$E'_1$ which gives the first two equalities in \eqqref(Dw=Ew). Two other
follow in the same way.

We want to adjust $L^+_2$ to $L^+_1$ and glue them together yielding the
desired Lagrangian sphere $S$. Our construction is as follows. Notice that
some neighbourhood $V_0$ of $F_0$ is symplectomorphic to the product
$(\Delta,\omega_\Delta){\times}(F_0,\omega|_{F_0})$. Applying an appropriate diffeomorphism we can
suppose that in the neighbourhood $V_0$ the ruling $\pr:X\to Y$ coincides with
the composition of the projection $\Delta{\times}F_0\to\Delta$ with an appropriate embedding
$\Delta\hookrightarrow Y$.  Let us denote the image of the latter embedding by $\Delta_0$ and identify
$\Delta$ with its image $\Delta_0\subset Y$. Further, we may assume that $\alpha_0:=\alpha\cap\Delta_0$ is an
interval. Observe that locally in $V_0$ the submanifold $W$ also is a product:
$W\cap V_0=\alpha_0{\times}F_0$. It follow that the Hamiltonian vector field $\mbfv^H$ is
horizontal in $V_0$, and every line flow of $\mbfv^H$ in $V_0$ has the form
$\alpha_0{\times}p$ for an appropriate point $p\in F_0$. Consequently, the discs
$L^+_1,L^+_2$ are also product in $V_0$: $L^+_i=\alpha_{0,i}{\times}\gamma_i$ where $\alpha_{0,1}$
and $\alpha_{0,2}$ are the parts of $\alpha_0$ lying to the left and to the right from
$y_0$.

The ``period condition'' $\int_{E_1}\omega=\int_{E_2}\omega$ transforms now into
$\int_{D_1}\omega=\int_{D_2}\omega$.  It follows that there exists a Hamiltonian
symplectomorphism $\Psi:(F_0,\omega)\to(F_0,\omega)$ which maps $D_1$ onto $D_2$ and thus
$\gamma_1$ onto $\gamma_2$. Let $\psi_t$ denote a generating (time depending) Hamiltonian
function. Extend $\psi_t$ to a function $\ti\psi_t$ in $V_0$ with the following
properties:
\begin{itemize}
\item The restriction $\ti\psi_t|_{F_0}$ coincides with $\psi_t$;
\item In some smaller neighbourhood $V^-_0$ of $F_0$ the function $\ti\psi_t$
 depends only on the vertical coordinate $p\in F_0$ (and the time coordinate
 $t$), but not on the horizontal coordinate $y\in Y$;
\item The function  $\ti\psi_t$ has compact support in $V_0$.
\end{itemize}
Using the last property we extend $\ti\psi_t$ by $0$ to the whole manifold
$X$. Let $\wt\Psi:(X,\omega)\to(X,\omega)$ be the symplectomorphism generated by the
Hamiltonian function $\ti\psi_t$. Apply $\wt\Psi$ to $L^+_2$ and denote the result
by $\wt L^+_2:=\wt\Psi(L^+_2)$. Then $\wt L^+_2$ is a Lagrangian disc which fits
to $L^+_1$: in the neighbourhood $V^-_0$ it has a product form %
$\wt L^+_2=\alpha^-_{0,2}{\times}\gamma_1$ for interval $\alpha^-_{0,2}$ on $\alpha$ to the right from
$y_0$. It follows that $S:=L^+_1\cup\wt L^+_2$ is a smooth Lagrangian sphere in
the desired homology class. 
\qed

\newsubsection[infla]{Inflation, deflation, and symplectic blow-up.} %
The standard symplectic form in $\rr^4$ is $\omega_\st:=dx_1\land dy_1+dx_2\land dy_2$. We
shall use the following results about symplectic blow-ups proved
\cite{McD-1,McD-2}:

\newlemma{lem-blow-up} \sli Let $(X,\omega)$ be a symplectic $4$-manifold and $B\subset(X,\omega)$
a symplectic embedding of the closed ball $B(R)\subset(\rr^4,\omega_\st)$ of radius $R$.
Then there exists a symplectic manifold $(X',\omega')$ and a symplectic exceptional
sphere $\Sigma\subset(X',\omega')$ with the following properties:
\begin{itemize}
\item $(X'\bs\Sigma,\omega')$ and $(X\bs B,\omega)$ are naturally symplectomorphic;
\item The area of\/ $\Sigma$ is $\int_\Sigma\omega'=\pi R^2$, the new volume is
 $\int_{X'}{\omega'{}^2}=\int_{X}{\omega^2}-(\pi R^2)^2$.
\end{itemize} 

\slii Let $(X,\omega_\mu)$ be a symplectic $4$-manifold and $\Sigma\subset(X,\omega_\mu)$ an
exceptional symplectic sphere with $\int_{\Sigma}\omega_\mu=\mu$. Then there exists a
symplectic manifolds $(X_0,\omega_0)$ and a symplectic embedding $B(R)\subset(X_0,\omega_0)$
of the radius $R=\sqrt{\mu/\pi}$ with the following properties:
\begin{itemize}
\item $X\bs\Sigma$ is naturally symplectomorphic to $X_0\bs B(R)$;
\item For any $0<\tau<\mu$, $X$ admits a symplectic form $\omega_\tau$ which coincides with
 $\omega_\mu$ outside a given neighbourhood of\/ $\Sigma$ and such that
 $\int_{\Sigma_\tau}\omega_\tau=\tau$. Moreover $(X\bs\Sigma,\omega_\tau)$ is naturally symplectomorphic to
 $X_0\bs B(r)$ with $r:=\sqrt{\tau/\pi}$.
\end{itemize} 

\sliii For any ruled symplectic $4$-manifold $(X,\omega)$ there exist a collection
$\{E_i\}$ of disjoint exceptional symplectic spheres whose contraction yields a
\emph{minimal} ruled symplectic $4$-manifold $(X_0,\omega_0)$.
\end{lem} 

\newdefi{deflate} We call the deformation $(X,\omega_\mu)\Rightarrow(X,\omega_\tau)$ the
\slsf{deflation} of\/ $\Sigma$, and the symplectic surgery $(X,\omega_\mu)\Rightarrow(X_0,\omega_0)$ the
\slsf{(symplectic) contraction} or \slsf{blow-down} of\/ $\Sigma$. The inverse
procedures $(X,\omega_\tau)\Rightarrow(X,\omega_\mu)$ and $(X,\omega,B)\Rightarrow(X',\omega',\Sigma)$ are called
\slsf{inflation} and \slsf{symplectic blow-up}.

\rm Notice that our notion of inflation(deflation) resembles but does not
coincide with that used in \cite{La-McD}. In both cases one changes the
cohomology class $[\omega]$ by an explicit construction which increases the volume,
but resulting deformations of the class $[\omega]$ are different.
\end{defi} 

The deflation construction can be applied to any symplectically embedded
surface $\Sigma$ in a symplectic $4$-manifold $(X,\omega)$. We give a brief description
of the construction. Let $\Sigma$ be an oriented closed surface $\pi:\bfp\to \Sigma$ a
principle $\SO(2)$-bundle. Equip $\bfp$ with a connection, let $\Theta$ be the
connection form and $\calf$ the curvature of the connection. Recall that
$\calf$ is a $2$-form on $\Sigma$ and $\Theta$ a $1$-form on $\bfp$ related by
$\pi^*\calf=-d\Theta$. Let $(r,\theta)$ be the polar coordinates on the disc $\Delta(R)$ of
radius $R$, so that the volume form is $rdr\land d\theta$ and the standard rotation
action of $\SO(2)$ is given by the formula $(r,\theta)\mapsto(r,\theta+\varphi)$.

Define $U_R:=\bfp\times_{\SO(2)}\Delta(R)$ as the $\Delta(R)$-bundle over $\Sigma$. Fix a
symplectic form $\omega_\Sigma$ on $\Sigma$. Let us use the same notation $\rho,\Theta,\calf,\omega_\Sigma$ for
the pull-backs of these objects on $U_R$ (notice that they are
well-defined). Finally, set
\begin{equation}\eqqno(omega-U-R)
\textstyle
\omega:=\omega_\Sigma+rdr\land\Theta-\frac{r^2}{2}\,\calf.
\end{equation}
A direct computation shows that $\omega$ is a symplectic form on $U_R$ provided the
curvature $\omega_\Sigma-s\calf$ is remains positive for every $s\in[0,R^2/2]$, \eg, for
negative $\calf$.  Notice that the restriction of $\omega$ on each fibre of the
natural projection $\pr:U_R\to \Sigma$ is the standard form $rdr\land d\theta$. %
Introduce a new coordinate $\rho:=\frac{r^2}{2}$, it varies in $\rho\in[0,R^2/2]$. Now
fix some $\tau>0$ and let $\rho'$, $r'$ be another functions defined by the
relations $\rho'=\rho+\tau$ and $\rho'=\frac{r'{}^2}{2}$, and such that $\rho'$ varies in
$[0,R^2/2+\tau]$. Substituting of $r'=r'(\tau)$ in \eqqref(omega-U-R) we obtain
\[\textstyle
\omega_\tau:=\omega_\Sigma+\tau\calf+r'dr'\land\Theta-\frac{r'{}^2}{2}\,\calf.
\]
As above, $\omega_\tau$ is a symplectic form on $U_{R'}=\bfp\times_{\SO(2)}\Delta(R')$ with
$R':=\sqrt{R^2+2\tau}$ provided $\omega_\Sigma-s\calf$ is remains positive for every
$s\in[-\tau,R^2/2]$. In this case we can consider $\omega_\tau$ as a deformation on
$\omega_0=\omega$. Observe also that by the very construction $(U_R\bs\Sigma,\omega)$ is
\emph{symplectomorphic} to $(U_{R'}\bs U_{\sqrt{2\tau}}, \omega_\tau)$.

Now assume that $\Sigma$ is symplectically embedded surface in some symplectic
$4$-manifold $(X,\omega)$. Let $\pi:\bfp\to \Sigma$ a principle $\SO(2)$-bundle associated
with the normal bundle to $\Sigma$. Fix a connection on $\bfp$ whose curvature
$\calf$ is proportional to $\omega|_\Sigma$, $\calf=a{\cdot}\omega|_\Sigma$ with the constant
$a=\frac{[\Sigma]^2}{\int_\Sigma\omega}$. Then for some $R>0$ the disc bundle $U_R$ constructed
above is symplectomorphic to some neighbourhood $U\subset X$ of $\Sigma$. Using the
symplectomorphisms $(U_R\bs\Sigma,\omega)\cong(U_{R'}\bs U_{\sqrt{2\tau}}, \omega_\tau)$ %
above as a transition function we cut out the neighbourhood $U$ isomorphic to
$U_R$ and glue in the neighbourhood $U_{R'(\tau)}$ with
$R'(\tau):=\sqrt{R^2+2\tau}$. Fixing appropriate diffeomorphisms $U_R\cong U_{R'(\tau)}$
identical near the boundary and near $\Sigma$ we obtain a family of closed $2$-form
$\omega_\tau$ on $X$ such that the restriction of $\omega_\tau$ on $\Sigma$ is
$\omega_\tau|_\Sigma=(1+a\tau){\cdot}\omega|_\Sigma$.  Moreover, $\omega_\tau$ is a symplectic form iff $1+a\tau>0$, in
particular, for every $\tau>0$ in the case $[\Sigma]^2\geq0$. Further, the cohomology class
of $\omega_\tau$ is $[\omega]+\tau[\Sigma]$ (here $[\Sigma]$ denotes the Poincar\'e dual to the homology
class of $\Sigma$) and thus the $\omega_\tau$-volume of $X$ increases by
$2\tau\int_\Sigma\omega+\tau^2[\Sigma]^2=\tau(2+a\tau)\int_\Sigma\omega$.

\newdefi{deflate-Sig} We call the constructed above deformation
$(X,\omega_\tau)$ the \slsf{deflation} of\/ $\Sigma$.
\end{defi} 

\newlemma{defla(-2)} \sli Let $(X,\omega)$ be a symplectic $4$-manifold, $\Sigma\subset(X,\omega)$
a closed oriented Lagrangian surface of genus $g(\Sigma)\neq1$, $U\subset X$ a neighbourhood
of $\Sigma$, and $\omega_\Sigma$ a symplectic form on $\Sigma$.  Then for some $\varepsilon>0$ there exists
a smooth family $\omega_t$ of symplectic structure on $X$ parametrised by
$t\in[-\varepsilon,\varepsilon]$ whose restriction on $\Sigma$ is $\omega_t|_\Sigma=t{\cdot}\omega_\Sigma$ and such that $\omega-\omega_t$
is supported in $U$. 

\slii Let $(X,\omega_\mu)$ be a symplectic $4$-manifold, $\Sigma\subset(X,\omega_\mu)$ a symplectic
$(-2)$-sphere with $\int_{\Sigma}\omega_\mu=\mu$, and $U\subset X$ a neighbourhood of $\Sigma$.  Then there
exists a smooth family $\omega_t$ of symplectic structure on $X$ parametrised by
$t\in[-\mu,\mu]$ such that the restriction satisfies $\omega_t|_\Sigma=\frac{t}{\mu}{\cdot}\omega_\mu|_\Sigma$
and such that $\omega_\mu-\omega_t$ is supported in $U$. In particular, $\Sigma$ is a
Lagrangian sphere with respect to $\omega_0$.
\end{lem}

\proof \sli Since the normal bundle $N\Sigma$ of a Lagrangian surface is isomorphic
to the tangent bundle $T\Sigma$ and $g(\Sigma)\neq1$, $N\Sigma$ is non-trivial. Find a
connection in $N\Sigma$ whose curvature $\calf$ is proportional to $\omega_\Sigma$,
$\omega_\Sigma=a\calf$. It follows that there exists a closed $2$-form $\psi$ with compact
support in $U$ whose restriction to $\Sigma$ is $\psi|_\Sigma=\calf$ (see
\cite{Bo-Tu}). Consequently, $a\psi$ gives an extension of $\omega_\Sigma$ to a closed
$2$-form $\psi$ with compact support in $U$. The formula $\omega_t:=\omega+ta\psi$ yields a
family with the desired properties.

\medskip\noindent%
\slii Let $\omega_t$ be the family of symplectic forms on $X$ constructed above by
the deflation of $\Sigma$ in which the parameter $t$ is chosen by the condition
$\int_\Sigma\omega_t=t$. In particular, $t$ varies in $(0,\mu]$. Further, we assume that
making the deflation of $\Sigma$ we start with a tubular neighbourhood $U_{R_0}\subset X$
in which the original form is given by \eqqref(omega-U-R) and the coordinate
$r$ varies in $[0,R_0]$ with sufficiently small $R_0>0$ to be chosen later.
After the full deflation $U_{R_0}$ is transformed into
$U_R=\bfp{\times_{\SO(2)}}\Delta(R)$ with $R=\sqrt{R_0^2+\mu/2}$ and the form $\omega$ is
deformed into a closed $2$-form $\omega_0$ given by
$\omega_0=rdr\land\Theta-\frac{r^2}{2}\,\calf$.

Consider the sphere $S^2$ equipped of radius $a>R/2$ with the coordinates
$\psi\in[0,\pi]$, $\theta\in[0,2\pi]$, and with the metric $g_a=a^2(d\psi^2+\sin^2\psi\,d\theta^2)$, so
that up to constants the coordinates $\psi$ and $\theta$ are the latitude and
respectively the longitude, and the standard rotation action of $\SO(2)$ is
given by the formula $(\psi,\theta)\mapsto(\psi,\theta+\varphi)$. Let $\omega_a:=a^2\sin\psi\,d\psi\land d\theta$ be volume
form. Since $2a>R$, the substitution $r=2a\,\sin\frac{\psi}{2}$ defines a
symplectic $\SO(2)$-equivariant embedding $\Delta(R)\subset(S^2,\omega_a)$. Set
$Q:=\bfp{\times_{\SO(2)}}(S^2,\omega_a)$ and equip it with the form
$\omega_0:=\sin\psi\,d\psi\land\Theta-2\sin^2\frac{\psi}{2}\,\calf$. Then $Q$ is a ruled
$4$-manifold, and the natural projection
$\pr_Q:Q=\bfp{\times_{\SO(2)}}(S^2,\omega_a)\to Q=\bfp{\times_{\SO(2)}}\pt=\Sigma$ is a
ruling. Further, the embedding $\Delta(R)\subset S^2$ induces the embedding $U_R\subset Q$
compatible with the bundle structures $\pr_Q:Q\to\Sigma$ and $U_R\to\Sigma$. Moreover, the
form $\omega_0$ on $Q$ is an extension of the form $\omega_0$ on $U_R$. Notice also that
the equations $\psi=0$ and respectively $\psi=\pi$ define two sections of the ruling
$\pr_Q:Q\to\Sigma$, such that the first one $\{\psi=0\}$ is the embedding $\Sigma\subset U_R\subset Q$. Let
us denote by $\Sigma_1$ the second section $\{\psi=\pi\}$ and by $V$ the complement
$V:=Q\bs U_R$. Then $\Sigma_1$ has self-intersection $+2$ and $V$ is a tubular
neighbourhood of $\Sigma_1$. 

Observe that applying appropriate diffeomorphisms we can consider the forms
$\omega_t$ in $U_R$ with $t\in(0,\mu]$ as deformations of the form $\omega_0$ constant near
the boundary $\partial U_R$. This allows us to extend each $\omega_t$ to a form $\omega_t$ on
$Q$, still denoted by $\omega_t$, such that each $\omega_t$ coincides with $\omega_0$ in $V$.
In particular, for each  $t\in(0,\mu]$ the form  $\omega_t$ is symplectic and $\Sigma$ is
$\omega_t$-symplectic. 

Now consider the real projective space $\rp^3$ with coordinates
$[x_0:x_1:x_2:x_3]$ and its complexification $\cp^3$ with coordinates
$[z_0:z_1:z_2:z_3]$. Denote by $\wt Q$ the quadric given by
$z_0^2=z_1^2+z_2^2+z_3^2$ and by $\wt\Sigma:=\wt Q\cap\rp^3$ its real form. Then
$\wt\Sigma$ is the standard sphere in $\rr^3\subset\rp^3$. Moreover, being a real form of
$\wt Q$, $\wt\Sigma$ is Lagrangian with respect to the Fubini-Study form
$\omega_{FS}$. Define $\ti\omega_0$ as the restriction $a\omega_{FS}|_{\wt Q}$ where the
constant $a$ is chosen from the condition $\int_{\wt
 Q}\ti\omega_0^2=\int_Q\omega_0^2$. %
The quadric $\wt Q$ is diffeomorphic to the product $\cp^1{\times}\cp^1=S^2{\times}S^2$
and $\wt\Sigma$ has homology class $\bff_1-\bff_2$ where $\bff_1,\bff_2$ are the
horizontal and respectively vertical homology classes. Further, the projection
$\wt\pr:\wt Q\to S^2$ on the first factor is a holomorphic ruling.  Consider also
the complex curve $C$ on $\wt Q$ given by the equation
$2z_0=z_1+z_2+z_3$. Then $C$ is a holomorphic section of $\wt\pr:\wt Q\to S^2$ in
the homology class $\bff_1+\bff_2$, which is disjoint from $\wt\Sigma$.  

For $t$ small enough, let $\ti\omega_t$ be a family of symplectic deformations of
the form $\ti\omega_0$ constructed in the part \sli of the lemma, normalised by the
condition $\int_{\wt\Sigma}\ti\omega_t=t$. Moreover, we may assume that $\ti\omega_t=\ti\omega_0$
in some neighbourhood $\wt V$ of $C$.  The classification of minimal ruled
symplectic $4$-manifolds with a section (see \cite{La-McD,McD-2}) implies that
for some sufficiently small values of $t_0>0$ above there exists a
symplectomorphism $\Phi:(Q,\omega_{t_0})\to(\wt Q,\ti\omega_{t_0})$ which maps $\Sigma$ onto
$\wt\Sigma$. Further, the moduli space of deformations of $C$ has complex dimension
$3$, and moreover, there exists a unique holomorphic curve in the homology
class $[C]=\bff_1+\bff_2$ passing through three given generic points.  It
follows that the image $\Phi(\Sigma_1)$ is symplectically isotopic to $C$, and
moreover, there exists an ambient symplectic isotopy $\wt Q\to\wt Q$ which
preserves $\wt\Sigma$ and moves $\Phi(\Sigma_1)$ into $C$. Thus we may assume that our
symplectomorphism $\Phi:(Q,\omega_{t_0})\to(\wt Q,\ti\omega_{t_0})$ maps $\Sigma_1$ onto
$C$. Furthermore, choosing the constant $R_0>0$ small enough we may
additionally assume that $\Phi$ maps the neighbourhood $V$ of $\Sigma_1$ onto $\wt V$.
Now, transferring the problem from $X$ to $\wt Q$ by means of the
symplectomorphism $\Phi$, we extend the constructed above family $\omega_t$ of
symplectic forms on $X$ from the interval $t\in(0,\mu]$ to an interval $t\in(-\delta,\mu]$
for some $\delta>0$. 

Repeating the argument above we show that for every sufficiently small $t_0>0$
there exists a symplectomorphism $\wt\Psi:(\wt Q,\ti\omega_{-t_0})\to(\wt Q,\ti\omega_{t_0})$
which is identical in $\wt V$ and maps $\wt\Sigma$ onto itself inverting the
orientation. By means of $\Phi$ we can transfer $\wt\Psi$ onto $X$ and obtain a
symplectomorphism $\Psi:(X,\omega_{-t_0})\to(X,\omega_{t_0})$ which is identical outside
$U_R$ and maps $\wt\Sigma$ onto itself inverting the orientation. Applying $\Psi$ to
the family $\omega_t$ with $t\in[t_0,\mu]$ we obtain the last piece
$(\omega_t)_{t\in[-\mu,-t_0]}$ of the family of symplectic structures with the desired
properties. %
\qed


\newsubsection[def-omega]{Deformation classes of symplectic structures.}
Recall that $\iom(X,\omega)$ denotes the space of symplectic forms $\omega'$
representing the class $[\omega]$.  Let $\iom_0(X,\omega)$ be the component containing
the form $\omega$. We start with the following

\newlemma{pi0-dfrm} \sli The group $\diff_0$ acts transitively on
$\iom_0(X,\omega)$.

\slii Let $\omega'_t$ and $\omega''_t$ be two families of symplectic forms on $X$
depending smoothly on $t\in[0,1]$ such that $\omega'_0=\omega''_0$ and such that
$[\omega'_t]=[\omega''_t]\in\sfh^2(X,\rr)$. Then there exists a family of diffeomorphisms
$F_t$ depending smoothly on $t\in[0,1]$ such that $F_0=\id_X$ and such that
$\omega''_t=F^*_t\omega'_t$ for each $t$. In particular, the forms $\omega'_t,\omega''_t$ lie in
the same component of $\iom(X,\omega'_t)$.

\sliii The set $\pi_0\iom(X,\omega)$ parametrises the deformation classes of
symplectic structures on $X$ cohomologous to $\omega$. Moreover, if $(X',\omega')$ is
obtained from $(X,\omega)$ by symplectic blow-up or contraction of exceptional
symplectic spheres, then $\pi_0\iom(X,\omega)$ and $\pi_0\iom(X',\omega')$ are naturally
identified.

\sliv Let $D_\omega(X)\subset\sfh^2(X,\rr)$ be the set of the $2$-cohomology classes
represented by symplectic forms defining a given orientation and a given
positive cone in $\sfh^2(X,\rr)$. Then every connected component of $D_\omega(X)$
is contractible.
\end{lem} 

\state Remark. In \cite{McD-4} McDuff constructed a manifold admitting pairs of
symplectic forms $\omega',\omega''$ which are cohomologous and deformationally
equivalent but nevertheless not isomorphic.

\proof \slip--\sliip. These facts follow easily by Moser's
technique.

\smallskip\noindent %
\sliii Two previous assertions of the lemma and \propo{ruled1}, \sliip, ensure
a global stability in the ``monodromic'' form: For any deformation path $\omega_t$
of symplectic structures no components of $\iom(X,\omega_t)$ appear or disappear,
\ie, the sets $\pi_0\iom(X,\omega_t)$ form a non-ramified covering of the parameter
path $t\in I=[0,1]$.  The assertion \sliv provides the triviality of the
monodromy along closed paths.

Observe any two symplectic balls in $(X,\omega)$ of the same sufficiently small
radius are symplectically isotopic. It follows that two isotopic
symplectic forms remain isotopic after a blow-up, provided the involved
symplectic balls have the same small enough radius. The uniqueness of
symplectic blow-down proved in \cite{McD-2}, \slsf{Theorem 1.1}, shows that
isotopic symplectic structure remain isotopic after contraction homologous
exceptional spheres.

\smallskip\noindent %
\sliv Let $\omega,\omega'$ be two cohomologous symplectic forms on $X$. Fix a basis
$\{\bfe_i\}$ of integer homology group $\sfh_2(X,\zz)$. Perturbing the forms
slightly we may assume that their periods $\int_{\bfe_i}\omega=\int_{\bfe_i}\omega'$ are
rationally independent. By \cite{La-McD,La} there exists a diffeomorphism
$F:X\to X$ with $F_*\omega=\omega'$. The condition on periods ensure that
$F_*:\sfh_2(X,\zz)\to\sfh_2(X,\zz)$ is identical.  It follows that if the class
$\bfe\in\sfh_2(X,\zz)$ is represented by an $\omega$-symplectic exceptional sphere
$E$, then $\bfe$ is represented also by an $\omega'$-symplectic exceptional sphere
$E'$.

\smallskip%
Assume that $X$ is not minimal. Let $\omega,\omega'$ be two symplectic forms on
$X$. Then there exist collections $\{E_1,\ldots,E_\ell\}$ and $\{E'_1,\ldots,E'_\ell\}$ of
exceptional spheres symplectic with respect to $\omega$ and $\omega'$ respectively,
minimal ruled symplectic manifolds $(\barr X,\bar\omega)$ and $(\barr X',\bar\omega')$,
and maps $\pi:X\to\barr X$, $\pi':X\to\barr X'$ which are symplectic contractions
contraction of the collections $\{E_1,\ldots,E_\ell\}$ and respectively $\{E'_1,\ldots,E'_\ell\}$.
Each of minimal ruled manifolds $\barr X,\barr X'$ is either a trivial bundle
$Y{\times}S^2$, or a non-trivial bundle $Y{\ti \times}S^2$. In the latter case we can
change the contracted spheres $E_i$ (resp. $E'_i$) such that the contraction
of the spheres from the new collection gives the trivial bundle
$Y{\times}S^2$. Indeed, let us fix on $X\cong Y{\ti \times}S^2$ an almost complex structure
$J$ compatible with the form $\bar\omega$. Let %
$\barr\pr:\barr X\to Y$ be the $J$-holomorphic ruling on $\barr X$. Deforming
$J$ appropriately we may assume that $J$ is integrable near some fibre of the
ruling.

\smallskip%
Now let $\eta_s$ be a family of cohomological classes parametrised by a sphere
$s\in S^k$, $k\geq1$, such that each class is represented by a symplectic form.
Find a sufficiently fine covering $U_\alpha$ of $S^k$ such that there exist
families of symplectic forms $\omega_{\alpha,s}$, each parametrised by $s\in U_\alpha$ and
representing the class $\eta_s$. For a fixed value $s_0\in S^k$, find a maximal
collection of $\omega_{s_0}$-symplectic exceptional spheres $E_1(s_0),\ldots,E_\ell(s_0)$.
Refining the covering if needed, we may assume that there exists a families
$J_{\alpha,s}$ of $\omega_{\alpha,s}$-tamed almost complex structures with the following
property: For each $s\in U_\alpha$ there exists an exceptional $J_{\alpha,s}$-holomorphic
sphere $E_{i,\alpha,s}$ homologous to $E_i$. Now deflate each sphere $E_{i,\alpha,s}$
forming a collection of families of symplectic forms $\omega_{\alpha,s,t}$, such that
the areas $\lan[\omega_{\alpha,s,t}], [E_i]\ran$ remain independent of $\alpha$. This is
possible in view of \lemma{lem-blow-up}. The above construction shows that
each element $\psi\in\pi_k(D_\omega(X))$ is represented by a family $\eta_s$ such that all
periods $\lan\eta_s, [E_i]\ran$ are sufficiently small. It follows that the
groups $\pi_k(D_\omega(X))$ are isomorphic to the group $\pi_k(D_\omega(X_0))$ where $X_0$
is the minimal ruled symplectic manifolds ($X_0=\cp^2$ in the rational case)
obtained by contraction all symplectic spheres $E_i$ above. 

The set $D_\omega(\cp^2)$ is the positive ray. In the irrational case the set
$D_\omega(X_0)$ was described in \cite{McD-2,La}, it is also contractible. 
This implies the assertion \slivp.
\qed

\newsubsection[dehn]{Symplectic Dehn twists in dimension $4$} We refer to
\cite{Sei1,Sei2} for a detailed description of symplectic Dehn twists.  In
this subsection we simply give a definition and list some properties.

Let $\Sigma$ be a Lagrangian sphere in a symplectic $4$-manifold $(X,\omega)$. Then
there exists a symplectomorphism $T_\Sigma:X\to X$, called \slsf{symplectic Dehn twist
 along $\Sigma$} which has the following properties:
\begin{itemize}
\item $T_\Sigma$ is supported in a given neighbourhood of\/ $\Sigma$;
\item $\Sigma$ is invariant with respect to $T_\Sigma$ and $T_\Sigma$ acts on $\Sigma$ as the
 antipodal map. In particular, $(T_\Sigma)_*[\Sigma]=-[\Sigma]$ and the action of $T_\bfs$ on
 the homology groups $\sfh_2(X,\zz),\sfh_2(X,\rr)$ is the reflection with
 respect to the hyperplane in $\sfh_2(X,\rr)$ orthogonal to $[\Sigma]$
 \begin{equation}\eqqno(S-refl)
[A]\in\sfh_2(X,\rr)\to [A]+([A]\cdot[\Sigma])\,[\Sigma].
\end{equation}
\item The symplectic isotopy class of $T_\Sigma$ is defined uniquely.
\item  $T_\Sigma^2$ is smoothly isotopic to the identity. 
\item Assume that $\Sigma'\subset(X,\omega)$ is a symplectic $(-2)$-sphere. Then there exists
 a deformation $\omega'$ of $\omega$ supported in a given neighbourhood of\/ $\Sigma'$ such that
 $\Sigma'$ is $\omega'$-Lagrangian. In particular, $T_{\Sigma'}$ is a well-defined
 diffeomorphism.
\end{itemize}

\newlemma{s2-lagr-exists} Let $(X,\omega)$ be a symplectic $4$-manifold and
$E_1,E_2$ two disjoint symplectic exceptional spheres. Then there exists a
deformation $\omega'$ of $\omega$ and an $\omega'$-symplectic sphere $\Sigma$ representing the
homology class $[E_1]-[E_2]$.
\end{lem} 

\state Remark. Observe that the Dehn twist along $\Sigma$ interchanges the classes
$[E_1],[E_2]$.

\proof Blow-down the sphere $E_2$. Let $(X^*,\omega^*)$ be the arising manifold,
$B_2$ the arising symplectic ball, and $E_1^*$ the image of $E_1$ in
$X^*$. Move the centre of $B_2$ into a point $p_2^*$ lying on $E^*_1$ and then
blow-up back using a small symplectic ball centred at $p_2^*$. The obtained
manifold $(X',\omega')$ is naturally diffeomorphic to $X$ and the form $\omega'$ is
isotopic to the form obtained from $\omega$ by an appropriate deflation of $E_2$.
The proper preimage of $E^*_1$ in $(X',\omega')$ is the desired $(-2)$-sphere $\Sigma$.
\qed

\medskip In \slsf{Theorems 3} and \slsf{4} we shall need some results about
Dehn twists. The first one is proved in \cite{Sei1}, \slsf{ Proposition 8.4.}

\newprop{braid-rel} Let $S_1$, $S_2$ be two Lagrangian spheres in some
symplectic $4$-manifold $X$ which intersect transversally in a single
point. Then the compositions of Dehn twists $T_{S_1}T_{S_2}T_{S_1}$ and
$T_{S_2}T_{S_1}T_{S_2}$ are symplectically isotopic.
\end{prop}

\medskip Another one result is a refining of the \slsf{Picard-Lefschetz
 formula}. In its formulation we use the following construction. Let $B$ be
the unit ball in $\cc^2$ equipped with the standard symplectic form
$\omega_\st=\frac{\isl}{2}(dz_1\land d\barz_1+dz_2\land d\barz_2) = dx_1\land dy_1+dx_2\land dy_2$ and
$\Delta_x,\Delta_y\subset B$ the coordinate discs: $\Delta_x=\{(x_1,0;x_2,0): x_1^2+x_2^2\leq1\}$ and
$\Delta_y=\{(0,y_1;0,y_2): y_1^2+y_2^2\leq1\}$.  For any $\varepsilon>0$ sufficiently small the
set
\[\textstyle
V_\varepsilon := \{\; (r\,\cos\varphi,r\,\sin\varphi;\;\frac{\varepsilon}{r}\,\cos\varphi,\,\frac{\varepsilon}r\,\sin\varphi)
\;:\; r>0, r^2+\frac{\varepsilon^2}{r^2}<1, \varphi\in[0,2\pi]\;\}.
\]
It is easy to see that $V_\varepsilon$ is $\omega_\st$-Lagrangian surface and that $V_\varepsilon$ is
topologically the connected sum of discs $\Delta_x,\Delta_y$. Further, after an
appropriate Lagrangian deformation $V_\varepsilon$ will coincide near the boundary with
the discs $\Delta_x,\Delta_y$. 

\newdefi{con-sum} We call the (Lagrangian or respectively smooth) isotopy
class of $V_\varepsilon$ relative the boundary the \slsf{connected sum of Lagrangian
 submanifolds $\Delta_x,\Delta_y$} at the intersection point and denote by $\Delta_x\#\Delta_y$.
\end{defi}
Let us notice that if orient discs in the usual way by the forms $dx_1\land dx_2$
and respectively $dy_1\land dy_2$ then the intersection index will be $-1$, but
the orientations on $V_\varepsilon$ induced by $\Delta_x$ and respectively $\Delta_y$ will be
opposite. Vice versa, orienting the discs $\Delta_x,\Delta_y$ in such a way that the
intersection index is $+1$, the induced orientation on $V_\varepsilon=\Delta_x\#\Delta_y$ is
well-defined.

This construction is local and will be used in a local situation. In the smooth
case we need only to know the orientation of the ambient $4$-manifold.  Now we
state the result.

\newprop{picard-lef} Let $S$ be a (Lagrangian) $(-2)$-spheres in some
(symplectic) $4$-manifold $X$, $T_S$ the induced (symplectic) Dehn twist, and
$L\subset X$ a (Lagrangian) surface intersecting $S$ transversally in a single
point. Then the surface $T_S(L)$ is (Lagrangian) isotopic to the sum $L\#S$.
\end{prop}

The proof can be found in \cite{Sei1}, \slsf{ Proposition 8.4.} There
P.~Seidel considers the special case when $L$ itself is a Lagrangian sphere,
but his construction applies in our more general situation.

\medskip%
{\sl Starting from this point we use the language and certain well-known facts
 from complex geometry.}

Consider the family of quadrics $Z_\lambda=\{z=(z_1,z_2,z_3):z_1^2+z_2^2+z_3^2=\lambda\}$ in
$\cc^3$, parametrised by the unit disc $\Delta:=\{|\lambda|<1\}$. Then $Z_\lambda$ are smooth for
$\lambda\neq0$ and $Z_0$ is the standard quadratic cone in $\cc^3$. The vertex $0$ of
$Z_0$ is a singular point of type $\rma_1$ and its desingularisation by means
of (holomorphic) blow-up yields the holomorphic cotangent bundle
$T^*\cp^1$. Moreover, if $C$ is a holomorphic rational $(-2)$-curve on a
smooth complex surface $X$ and $U$ its neighbourhood, then one can contract
$C$ into a singular point such that the obtained complex space $U'$ is
isomorphic a neighbourhood of the vertex $0$ of $Z_0$.

Further, for $\lambda=\rho e^{\isl\theta}$ with $\rho>0$ the quadric $Z_\lambda$ contains the
Lagrangian sphere $\Sigma_\lambda$ given by
$\Sigma_\lambda:=\{e^{\isl\theta/2}(x_1,x_2,x_3):x_1^2+x_2^2+x_3^2=\rho\}$ where
$(x_1,x_2,x_3)\in\rr^3$.  It is not difficult to show that the Lefschetz
monodromy of the family $Z_\lambda$ around the origin $0\in\Delta$ is exactly the Dehn
twist along the Lagrangian sphere $\Sigma_\lambda$.

Now consider the following situation. Let $\Delta'$ be the unit disc with the
coordinate $w'$, $X'$ the product $\Delta'\times\cp^1$, and $\pr':X'\to\Delta'$ the natural
holomorphic ruling. Blow up a point on the central fibre $\{0\}\times\cp^1$. Denote
the obtained manifold by $X''$ and the obvious induced ruling by
$\pr'':X''\to\Delta'$. Then the new central fibre is the union of two exceptional
rational holomorphic curves, say $E_1$ and $E_2$. Further, in a neighbourhood
of the nodal point of the central fibre there exist local complex coordinates
$z_1,z_2$ in which $\pr''$ is given $(z_1,z_2)\mapsto w'=z_1^2+z_2^2$. Take another
unit disc $\Delta$ with the coordinate $w$ and consider the covering $f:\Delta\to\Delta'$ given
by $w\mapsto w'=w^2$. Let $f^*X''=:X'''$ be the pull-back of $X''$ and
$\pr''':X'''\to\Delta$ the induced projection. The relations $w'=w^2$ and
$w'=z_1^2+z_2^2$ imply that $X'''$ contains a singular point of type $\rma_1$,
say $p^*$, such that the functions $w,z_1,z_2$ generate the local ring
$\scro_{X''',p^*}$ and satisfy the relation $w^2-z_1^2-z_2^2=0$. In other
words, $w^2-z_1^2-z_2^2=0$ is a local equation for $X'''$ at $p^*$.
Furthermore, the obvious induced projection $\pr''':X'''\to\Delta$ is given locally
by $(z_1,z_2)\mapsto\pm\sqrt{z_1^2+z_2^2}$. Let $X_0$ be the desingularisation of
$X'''$, $\Sigma$ the arising holomorphic $(-2)$-sphere, and $\pr_0:X_0\to\Delta$ the
induced ruling. It has a unique singular fibre which is a chain of curves
$E_1,\Sigma,E_2$. On the other hand, let $X_\lambda$ be the deformations of $X'''$ given
in a neighbourhood of the singular point $p^*$ by the equation
$w^2-z_1^2-z_2^2=\lambda$, $\lambda\in\Delta$. Then $X_\lambda$ admit rulings $\pr_\lambda:X_\lambda\to\Delta$ locally
given by $(w,z_1,z_2)\mapsto w$. In particular, for $\lambda\neq0$ the ruling $\pr_\lambda$ has two
ordinary singular fibre over the singular values $w_\pm(\lambda):=\pm\sqrt{\lambda}$.

\smallskip%
\noindent {\sl Braid groups.} For the definition and basic properties of
\slsf{braid group} $\br_d$, \slsf{pure braid group} $\pbr_d$ we refer to
\cite{Bir}. For the definition of the usual $2$-dimensional Dehn twist $T_\delta$
along an embedded circle $\delta$ on a surface $Y$ we refer to \cite{F-M,Ivn}.

Below we shall make use of the following facts: (i) In one of its geometric
realisations, $\br_d$ is the fundamental group of the configuration space of
$d$ pairwise distinct points on a disc $\Delta$ or in the complex plane $\cc$. Let
$w^*_1,\ldots,w^*_d\in\Delta$ be a basic configuration. (ii) In this realisation, the
standard generators of $\br_d$ are so called \slsf{half-twists} $\tau$ exchanging
exactly two points $w^*_i, w^*_j$. (iii) The square $\tau^2$ of such a
generator $\tau$ is a Dehn twist $T_\delta$ in $\Delta$ along an embedded closed curve $\delta$
surrounding the points $w^*_i, w^*_j$ and no other points. Such Dehn twists
$\tau^2=T_\delta$ form a system of generators of the pure braid group $\pbr$.

\begin{thm}\label{twists-4-2} \sli The monodromy $\Phi$ of the above family $(X_\lambda,\pr_\lambda)$
along the path $|\lambda|=\rho>0$ act on $X_\rho$ as the Dehn twist along $\Sigma$ and on the
critical values $w_\pm:=\pm\sqrt{\rho}$ as a geometric half-twist in the braid group
$\br_2(\Delta)$.

In particular, the monodromy interchanges the fibres over $+\sqrt{\rho}$
and $-\sqrt{\rho}$.

\slii Let $\delta\subset\Delta$ be a closed curve surrounding the points $w_\pm=\pm\sqrt{\rho}$, $A$
an annual neighbourhood of $\delta$ disjoint from $w_\pm$, and
$f:\Delta\to\Delta$ its geometric realisation of the Dehn twist $T_\delta$ supported in
$A$. Trivialise the set $\pr_\rho\inv(A)\cong A\times\cp^1$ and define the map $F:X_\rho\to X_\rho$
setting $F:=f\times\id_{\cp^1}$ in   $\pr_\rho\inv(A)$ and $F:=\id$ outside. Then $F$
is smoothly isotopic to the square $\Phi^2$ of the monodromy $\Phi$ above and thus
smoothly isotopic to the identity map $\id:X_\rho\to X_\rho$.  
\end{thm} 

\proof  \sli follows directly from the construction of the family. 

\smallskip\noindent%
\slii Let us make a geometric analysis of the problem. As we show below,
changing coordinates one can realise the Dehn twist $T_\delta$ as follows: Instead
of points $w_\pm=\pm\sqrt{\rho}$ we take new $w_0:=0$, the centre of $\Delta$, and
$w_1:=\rho$, letting $w_1$ make one turn around $w_0$. The parameter expression
for this turn is is given by $w_1(t)=\rho e^{2\pi\isl t}$. The equation in
coordinates is $(w-\lambda)^2=\lambda^2$. Notice that the appearance of $\lambda^2$ instead of
$\lambda$ expresses the fact that we take the full twist which is the square of the
half-twist. Notice also that the square of the Dehn twist along $\Sigma$ is given
by a local equation $(w-\lambda)^2-z_1^2-z_2^2=\lambda^2$. However, this argument is not
correct: The Dehn twist $T_\Sigma$ is supported in a neighbourhood of the curve $\Sigma$
which is 2-dimensional, whereas the map $F$ is supported in the set $A\times\cp^1$
which is a neighbourhood of a 3-dimensional set $\delta\times\cp^1$.

The idea of the proof is to revert certain part of $F$ by means of an
appropriate isotopy such that the obtained will have substantially smaller
support. First, we rearrange our geometric objects. Recall that we changed our
basic constellation to points $w_0:=0$ which will be constant and $w_1:=\rho$
which will move. As the new curve $\delta$ we take the circle $|w|=\rho+2\varepsilon$ with
sufficiently small $\varepsilon$. Choose a function $\chi_1(r)$ such that $0\leq\chi_1(r)\leq1$,
$\chi_1(r)\equiv1$ for $r\leq\rho+\varepsilon$, and $\chi_1(r)\equiv0$ for $r\geq\rho+4\varepsilon$.  Then the Dehn twist
$T_\delta$ can be realised by a map $f_1:\Delta\to\Delta$ given in polar coordinates $(r,\theta)$ by
the formula $f_1(r,\theta)=(r, \theta+2\pi\chi_1(r))$. The geometric meaning of $f_1$ is
clear: we make the full rotation on the disc $\{|w|\leq\rho+\varepsilon\}$ leaving fixed the
annulus $\{\rho+4\varepsilon\leq|w|\leq1\}$. This description can be viewed via the family
$f_{1,t}(r,\theta):=(r, \theta+t\cdot2\pi\chi_1(r))$.

Now choose another function $0\leq\chi_2(r)\leq1$ such that $\chi_2(r)\equiv\chi_1(r)$ for
$\rho+\varepsilon\leq r\leq\rho+3\varepsilon$ and $\chi_2(r)\equiv0$ for $r\leq\rho-\varepsilon$. Then the map $f_2:\Delta\to\Delta$ given by
$f_2(r,\theta)=(r, \theta+2\pi\chi_2(r))$ is isotopic to $f_1$ relative boundary $\partial\Delta$ and the
points $w_0,w_1$. Thus $f_2$ is another geometric realisation of the Dehn twist
$T_\delta$.  The geometric meaning of $f_2$ is as follows: we make the full
rotation of the circle $\{|w|=\rho\}$ leaving fixed the annulus $\{\rho+4\varepsilon\leq|w|\leq1\}$ and
the disc $\{|w|\leq\rho-\varepsilon\}$, and as in the case of $f_1$, the family
$f_{2,t}(r,\theta):=(r, \theta+t\cdot2\pi\chi_2(r))$ illustrates the interpretation.

Next, we want to make a similar rearrangement in $\cp^1$-direction. For this
purpose let us realise $X_\rho$ as the blow-up of the product $\Delta\times\cp^1$ in two
points $p_0$ and $p_1$ lying on the fibres $C_0:=\{w_0\}\times\cp^1$ and
$C_1:=\{w_1\}\times\cp^1$.  The exact position of points plays no role. Let $z$ be a
complex projective coordinate on $\cp^1$ such that $z=0$ at $p_1$. Fix a
function $0\leq\eta(z)\leq1$ such that $\eta(z)\equiv1$ in the disc $\{|z|\leq\varepsilon\}$ and $\eta(z)\equiv0$
outside the disc $\{|z|\geq2\varepsilon\}$ where as above $\varepsilon$ is sufficiently small.
Define the family of  maps $F_{2,t}:\Delta\times\cp^1\to\Delta\times\cp^1$ setting 
\[
F_{2,t}(r,\theta;z):= ( f_{2,t\eta(z)+1-t}(r,\theta); z) = (r, \theta+(t\eta(z)+1-t)\cdot2\pi\chi_2(r); z). 
\]
Then $F_{2,t}$ is an isotopy of diffeomorphisms relative points $p_0,p_1$,
$F_{2,t}$ are identical in a neighbourhood of $p_0,p_1$, and for $t=0$ we have
$F_{2,0}=f_2\times\id_z$. The geometric structure of $F_{2,1}$ is as follows. Some
small neighbourhood of $p_1$ turns around the fibre $C_0=\{w_0\}\times\cp^1$ along the
circle $\{|w|=\rho\}\times\{p_1\}$. 

Finally, define the map $F_2:X_\rho\to X_\rho$ as the blow-up of the map
$F_{2,1}:\Delta\times\cp^1\to\Delta\times\cp^1$ at points $p_0,p_1$. The geometric picture of $F_2$
is reformulation of that for the map $F_{2,1}$. However, it has also another
description: It is the monodromy map of the family of blow-ups of $\Delta\times\cp^1$ at
two points in which the point $p_0$ lies on the central fibre
$C_0=\{w_0\}\times\cp^1$ and $p_2$ turns around $C_0$.  But this is exactly the
family $X_\lambda$ above, so $F_2:X_\rho\to X_\rho$ is smoothly isotopic to the squared
symplectic Dehn twist $\Phi^2$.
\qed

\newsubsection[cox]{Reflection- and Coxeter-Weyl groups.} %
In this paragraph we collect some basic facts about Coxeter systems and
Coxeter-Weyl groups, relying as the reference on the books \cite{Bou,Hum}.

\newdefi{def-cox}
 A \slsf{Coxeter matrix} over a set $\cals =\{s_1, \ldots, s_r\}$ is a symmetric $r\times
 r$-matrix %
 $M=M_\cals=(m_{ij})$ with entries in $\nn \cup \{\infty\}$, such that $m_{ii} =1$ and
 $m_{ij} \geq2$ for $i\neq j$.  A \slsf{Coxeter system} is a set $\cals$ equipped
 with a Coxeter matrix $M_\cals$.  If $s=s_i, s'=s_j \in \cals$, then $m_{ij}$
 is also denoted by $m_{s,s'}$. The number $r$ of elements in $\cals$ is the
 \slsf{rank} of $\cals$. 

 The \slsf{Coxeter graph} $\Delta = \Delta(\cals)$ associated with a Coxeter system
 $\cals$ has $\cals$ as the set of vertices; $s_i$ and $s_j$ are connected by
 a edge iff $m_{ij} \geq3$; if $m_{ij} \geq4$, this edge is labelled by the entry
 $m_{ij}=m_{ji}$.  In this paper we shall always have $m_{ij}\leq4$ or
 $m_{ij}=\infty$. We use a ``double line'' edge
 $\xymatrix@1@C-12pt@R-15pt{s_i\ar@<1pt>@{=}[r] & s_j}$ to denote an entry
 $m_{ij}=m_{ji}=4$ mimicking the notation used in Dynkin diagrams.

The labelled Coxeter graph $\Delta$ completely determines the associated Coxeter
matrix.  Given an appropriate graph $\Delta$ we denote by $\cals(\Delta)$ the associated
Coxeter system.

For any subset $\cals' \subset \cals$ the induced Coxeter matrix $M'$ is the
restriction of $M= M_\cals$ onto $\cals'$. The associated graph
$\Delta':=\Delta(\cals')$ is a full subgraph of $\Delta$, \ie, any two vertices
$s_1,s_2\in\cals'$ are connected and labelled in $\Delta'$ in the same way as in $\Delta$.
\end{defi}

\newdefi{def-gra-group}  
For any two letters $a,b$ and a non-negative integer $m$ we denote by $\langle ab\rangle^m$
the word $abab\ldots$ of the length $m$ consisting of alternating letters $a$ and
$b$. If $a$ and $b$ lie in a group $G$, then $\langle ab\rangle^m$ denotes the
corresponding product in $G$.

The \slsf{Coxeter--Weyl group} $\sfw(\cals)= \sfw(\cals, M)$ of a Coxeter
system $\cals$ is the group generated by $\cals$ with the \slsf{reflection
 relation} $s_j^2=1$ for each $s_j\in\cals$ and \slsf{(generalised) braid relations}
\[
\langle s_is_j\rangle^{m_{ij}} = \langle s_js_i\rangle^{m_{ij}}
\qquad
\text{for each $i\neq j$, such that $m_{ij}<\infty$}
\]
as defining relations. 
\end{defi}

The standard examples of Coxeter--Weyl groups are symmetric groups $\sym_d$
and Weyl groups $\sfw(\frg)$ of semi-simple complex Lie algebras $\frg$. The
corresponding generators of the Coxeter--Weyl groups are elementary
transpositions $\tau_{i,i+1}$ in $\sym_d$ and respectively the reflections
corresponding to \slsf{simple roots} $\alpha_i$ of the Lie algebra $\frg$, see \eg
\cite{Bou} for details.

In Lie theory one uses \slsf{Dynkin diagrams} to encode the combinatorial
structure of \slsf{systems of simple roots} of semi-simple complex Lie
algebras $\frg$. Those diagrams are graphs with simple, double-, or
triple-line edges, equipped with direction arrow like
$\xymatrix@1@C-11pt@R-15pt{s\ar@<1pt>@2{->}[r] & s'}$ or
$\xymatrix@1@C-11pt@R-15pt{s\ar@<1pt>@3{->}[r] & s'}$ in the last two
cases. The triple line edges correspond to $m_{ij}=6$ and will not appear in
this paper.  For a moment let us consider the double-line edges with
an arrow head as those without arrow it, \ie,
$\xymatrix@1@C-11pt@R-15pt{s\ar@<1pt>@2{->}[r] & s'}$ as
$\xymatrix@1@C-11pt@R-15pt{s\ar@<1pt>@2{-}[r] & s'}$. This rule allows to
transform an appropriate Dynkin diagram $\Delta$ into a Coxeter graph, and so to
construct a Coxeter-Weyl group $\sfw(\Delta)$. Notice that this group is exactly
the Weyl group $\sfw(\frg)$ of the complex Lie $\frg$ algebra with the Dynkin
diagram $\Delta$.

We maintain the notation $\rma_\ell,\rmb_\ell,\rmd_\ell$ for the Coxeter systems and
the corresponding graphs from the Lie theory. By
$\rml_\ell(m_{1},m_{2},\ldots,m_{\ell-1})$ with $m_i\geq3$ we denote the Coxeter system of
rank $\ell$ whose Coxeter matrix $M_\cals$ has by-diagonal entries
$m_{i,i+1}=m_{i+1,i}$ are equal $m_i$ and such that all remaining entries
$m_{ij}$ with $|i-j|\geq2$ are $2$. In other words,
$\rml_\ell(m_{1},m_{2},\ldots,m_{\ell-1})$ is the Coxeter system with the linear Coxeter
graph of length $\ell$ with markings $m_1,m_2,\ldots,m_\ell$ on the edges. Recall that we
set no marking if $m_{i}=m_{i,i+1}=3$ and use a double line to mark an edge
with $m_{i}=m_{i,i+1}=4$.

Every Coxeter group has a natural representation $\rho_\cals$ in certain vector
space $V_\cals$ equipped with some special inner product such that
$\rho_\cals(s_j)$ is a reflection in $V_\cals$. Namely, for a given Coxeter
system $\cals$ let $V_\cals$ be the real vector space with the set $\cals$ as
a basis. Let $\mbfe_j$ or $\mbfe_{s_j}$ denote the corresponding basis
vectors. Define an inner product in $V_\cals$ setting
$\lan\mbfe_j,\mbfe_j\ran:=-1$ and
$\lan\mbfe_i,\mbfe_j\ran:=\cos(\frac{\pi}{m_{ij}})$.  The latter condition can
be interpreted as follows: In the case $m_{ij}<\infty$ the vectors $\mbfe_i$ and
$\mbfe_j$ span a negatively defined 2-plane and the angle between them is
$\pi(1-\frac{1}{m_{ij}})$. In the case $m_{ij}=\infty$ the plane spanned by the
vectors $\mbfe_i$ and $\mbfe_j$ is degenerate with 1-dimensional isotropic
space generated by the vector  $\mbfe_i+\mbfe_j$. 

Notice that usually one uses the opposite sign
convention with $\lan\mbfe_j,\mbfe_j\ran:=+1$, however the choice made above
will give us later the correct signature of $V_\cals$.

\newdefi{geom-repres} \sl The representation
 $\rho_\cals:\sfw(\cals)\to\Gl(V_\cals)$ is called \slsf{geometric}.
\end{defi}

With this definition, one can prove the
following results, see \eg\ \cite{Bou} or~\cite{Hum}:

\begin{prop}\label{refl-repres} 
 Let $\cals$ be a Coxeter system, and $(V_\cals,\lan\cdot,\cdot\ran)$ the geometric
 representation. Denote by $\sigma_j$ the reflection in $V_\cals$ with respect to
 the hyperplane orthogonal to $\mbfe_j$. Then there exists a natural
 isomorphism $\rho_\cals$ from $\sfw(\cals)$ to the group generated by the
 collection $\{\sigma_j\}$ such that $s_j$ is mapped onto $\sigma_j$.
\end{prop}

\begin{prop}\label{refl-finite} 
The Coxeter-Weyl group $\sfw(\cals)$ is finite if and only if the inner
product $\lan\cdot,\cdot\ran$ on $V_\cals$ is negative definite.\\[3pt]
\rm In this case we say that the system $\cals$ has \slsf{finite type}.
\end{prop}

\smallskip%
Let us make a few remarks about the case $m_{ij}=\infty$. In this case no relation
between the generators $s_i$ and $s_j$ is imposed. The Coxeter system on two
generators $s_1,s_2$ with $m_{1,2}=\infty$ is usually denoted by $\rmi_2(\infty)$ or by
$\wt\rma_1$. It has the Coxeter graph %
$\xymatrix@C-10pt@R-13pt{ s_1\ar@<1pt>@{-}[r]^\infty & s_2}$.  The second notation
$\wt\rma_1$ comes from the fact that the corresponding Coxeter group
$\sfw(\rmi_2(\infty))$ is the affine Weyl group of the Lie algebra
$\frg(\rma_1)=\frak{sl}(2,\cc)$. The group $\sfw(\rmi_2(\infty))$ is isomorphic to
the semi-direct product $\zz_2\ltimes\zz$ in which $\zz_2$ acts on $\zz$ by
inversion, and the generators $s:=(-\id,0)$ and $T:=(\id,1)$ of $\zz_2\ltimes\zz$
are represented by elements $s_1$ and $s_1s_2$ from $\sfw(\rmi_2(\infty))$.

\smallskip %
As one can see from \slsf{Theorems 4' {\rm and} 4''} Coxeter-Weyl groups is a
natural tool in description of the action of the diffeotopy group on the
$2$-homology $\sfh_2(X,\rr)$ of a rational and ruled $4$-manifold $X$. We want
to show that the integer homology $\sfh_2(X,\rr)$ also fits in this
description.

\newdefi{crystal}\sl
 Let $\cals$ be a Coxeter system whose off-diagonal entries $m_{ij}=m_{ij}$
 ($i\neq j$) of the Coxeter matrix are $2,3,4$ or $\infty$. 
A \slsf{crystallographic structure} on a $\cals$ is given by a
decomposition of the set $\cals$ into disjoint subsets $\cals',\cals''$ with
the following properties:
\begin{enumerate}
\item if $m_{ij}=3$ or $\infty$, then $s_i,s_j$ lie both in $\cals'$ or both in $\cals''$;
\item if $m_{ij}=4$, then $s_i,s_j$ lie in different parts of the
 decomposition $\cals=\cals'\sqcup\cals''$.
\end{enumerate}
The vertices from $\cals'$ and  $\cals''$ are called \slsf{short} 
and respectively \slsf{long}.

In term of Coxeter graphs the definition means that we have only single-,
double-, and $\infty$-line edges, the edges connecting vertices 
inside each part $\cals',\cals''$ are single- or $\infty$-lines, and double-line edges
connect vertices from different parts.

The \slsf{Dynkin diagram} of a crystallographic structure is obtained from the
Coxeter graph adding the arrow head to each double-line edge pointing on the
short vertex.

The \slsf{lattice $V^\zz_\cals$ associated to the crystallographic structure}
is generated by vectors $\mbfv_j$ such that $\mbfv_j=\mbfe_j$ in the case of a
short edge $s_j\in\cals'$ and $\mbfv_j=\sqrt2\,\mbfe_j$ in the long edges case
$s_j\in\cals''$. In particular, $\lan\mbfv_j,\mbfv_j\ran=-1$ for short edges and
$\lan\mbfv_j,\mbfv_j\ran=-2$ for long ones.
\end{defi}

It is not difficult to see that the Dynkin diagram of a crystallographic
Coxeter system $\cals$ determines the system if and only if it has at least
one double line edge (arrow) in each its connected component. 

\state Remark. In a more general (and probably more correct) definition, a
Coxeter system $\cals$ is crystallographic if it admits a full rank lattice
$V^\zz_\cals$ in the geometric representation $V_\cals$ invariant with respect
to $\sfw(\cals)$.  Such systems $\cals$ can also have triple line edges with
similar properties.  Since such systems are irrelevant to the purpose of this
paper, we use a more restrictive definition.

\begin{prop}\label{cryst-E8} \sli For every crystallographic
Coxeter system $\cals$ the associated lattice $V^\zz_\cals$ is invariant with
respect to the Coxeter group $\sfw(\cals)$.

\slii Under the hypotheses of \slsf{Theorem 4'} the Coxeter system $\cals(X)$
has natural crystallographic structure given by the Dynkin diagram
\begin{align} 
&\xymatrix@C-10pt@R-13pt{
& & & s_0\ar@{-}[d]
\\
\rmbe_{\ell+1}\ (\ell\geq4): 
&  s_1 \ar@<2pt>@{-}[r] & s_2\ar@<2pt>@{-}[r] & s_3 \ar@<2pt>@{-}[r] &
\ldots \ar@<2pt>@{-}[r] & s_{\ell-2} \ar@<2pt>@{-}[r] & s_{\ell-1}\ar@<1pt>@2{->}[r] & s_{\ell}
}
\\
\noalign{\ \!\!\!\!\!\!in the case $\ell\geq4$ or respectively}
&\xymatrix@C-10pt@R-13pt{
\rml_4(3,4,4)\ (\ell=3): 
&  s_1 \ar@<2pt>@{-}[r] & s_2\ar@<1pt>@2{->}[r] & 
s_{3}\ar@<1pt>@2{<-}[r] & s_{0}
}
\end{align}
\ \\[-17pt]
in the case $\ell=3$. The associated lattice  $(V^\zz_{\cals(X)},\lan\cdot,\cdot\ran)$ is
the integer homology group $\sfh_2(X,\zz)$ with the intersection form.

\sliii Under the hypotheses of \slsf{Theorem 4''} the Coxeter system $\cals(X)$
has natural crystallographic structure given by the Dynkin diagram
\begin{align} 
&\xymatrix@C-10pt@R-13pt{
& & s_0\ar@{-}[d]
\\
\rmbd_{\ell+1}=\wt\rmb_\ell\ (\ell\geq3): 
&  s_1 \ar@<2pt>@{-}[r] & s_2\ar@<2pt>@{-}[r] & 
\ldots \ar@<2pt>@{-}[r] & s_{\ell-2} \ar@<2pt>@{-}[r] & s_{\ell-1}\ar@<1pt>@2{->}[r] & s_{\ell}
}
\\
\noalign{\ \!\!\!\!\!\!in the case $\ell\geq3$ or respectively}
&\xymatrix@C-10pt@R-13pt{
\rml_3(4,4)=\wt\rmb_2\ (\ell=2): 
&  s_1\ar@<1pt>@2{->}[r] & 
s_{2}\ar@<1pt>@2{<-}[r] & s_{0}
}
\end{align}
\ \\[-17pt]
in the case $\ell=3$. The associated lattice  $(V^\zz_{\cals(X)},\lan\cdot,\cdot\ran)$ is
the sublattice in $\sfh_2(X,\zz)$ consisting of elements orthogonal to
$\bff$. This sublattice is generated by $\bfe_1,\ldots,\bfe_\ell$ and $\bff$.

\sliv Under the hypotheses of \slsf{Theorem 4'{'}'} the Coxeter system
$\rml_3(4,\infty)$ giving the Coxeter-Weyl group $\sfw(\rml_3(4,\infty))\cong\Gamma_H(X)$ of
the manifold $X=(S^2{\times}S^2)\#\cp^2$ has natural crystallographic structure
given by the Dynkin diagram $\xymatrix@C-10pt@R-13pt{s_1\ar@<1pt>@2{->}[r] &
 s_2\ar@<1pt>@1{-}[r]^\infty&s^*_0}\!\!\!$.  The associated lattice $(V^\zz,\lan\cdot,\cdot\ran)$
is the integer homology group $\sfh_2(X,\zz)$ with the intersection form.
\end{prop}

\proof \sli One needs only to verify the invariance of $V^\zz_\cals$ with
respect to generating reflections $\sigma_i$ from \propo{refl-repres}. This can be
done explicitly.  

The parts \sliip, \sliiip, and \sliv are also verified explicitly.  \qed

\smallskip%
Finally, we discuss the meaning of the period conditions from the point of
view of Coxeter-Weyl groups. This will give us the part \slii of
\slsf{Theorems 4' {\rm and} 4''}. The following results were proved in
\cite{Vi}, see also \cite{Vi-Sh}.

\newprop{P-fund} Let $\cals=\{s_1,\ldots,s_\ell\}$ be a Coxeter system of finite rank
such that the scalar product on the geometric representation $V_\cals$ has
Lorentz signature.  Set $\calk:=\allowbreak\{v\in V_\cals:v^2>0\}$ and
$\bfp:=\{v\in\calk:\lan\mbfe_i,v\ran\geq0,i=1,\ldots,\ell\}$.

\sli  The set $\calk$ consists of two
convex cones such that one of them, called the \slsf{positive cone} and
denoted by $\calk_+$, contains $\bfp$. Moreover, $\calk_+$ is
$\sfw(\cals)$-invariant and the subset $\bfp$ is a fundamental domain for the
action of\/ $\sfw(\cals)$ on $\calk_+$. 

Further, let $H_1,\ldots,H_\ell$ be the hyperplanes orthogonal to the vectors
$\mbfe_1,\ldots,\mbfe_\ell$ respectively, and $H_\alpha$ the hyperplanes obtained from $H_1,\ldots,H_\ell$ by
means of the action of $\sfw(\cals)$. Then the union $\cup_\alpha H_\alpha$ cuts $\calk_+$
into convex cones $\bfp^\circ_{\!\!p}$ such that the closure $\bfp_{\!p}$ of each
$\bfp_{\!\!p}^\circ$ is also a fundamental domain for the
action of\/ $\sfw(\cals)$ on $\calk_+$. 

\smallskip %
\slii Let $v$ be any vector from $\bfp$. Denote by $\cals(v)$ the Coxeter
subsystem of $\cals$ consisting of generators $s_j\in\cals$ such that $H_j$
contains $v$. Then $\cals(v)$ is the Coxeter system of finite type, so that
the subgroup $\sfw(\cals(v))<\sfw(\cals)$ is finite. Moreover, if $v$ lies on
some hyperplane $H_\alpha$, then there exists $\phi\in\sfw(\cals(v))$ which maps $H_\alpha$
onto some  hyperplane $H_i$ such that the corresponding reflection $s_i$
belongs to  $\cals(v)$.

\end{prop}


\newsection[sec-ti-w2]{Secondary Stiefel-Whitney classes and
 symplectomorphisms.}

\newsubsection[defi-ti-w2]{Secondary cohomology classes.}  

Let $X$ be a smooth manifold and $f:X\to X$ a diffeomorphism which acts
trivially in $\sfh^\bullet(X,\zz_2)$. Denote $I:=[0,1]$. Define the \slsf{map torus}
$T(X,f)$ as the quotient of $X{\times}[0,1]$ under the identification
$(x,1)\sim(f(x),0)$. In other words, $T(X,f)$ can be realised as a fibre bundle
over $S^1$ with fibre $X$ and monodromy $f$. The condition on
$\sfh^\bullet(X,\zz_2)$ implies the K\"unneth formula
$\sfh^\bullet(T(X,f),\zz_2)\cong\sfh^\bullet(X,\zz_2)\otimes\sfh^\bullet(S^1,\zz_2)$. Define the class
$[S^1]\in\sfh_1(T(X,f),\zz_2)$ as the image of the fundamental class of the
circle $[S^1]$ under the K\"unneth isomorphism.

\newdefi{w2-f} In the above situation, the \slsf{secondary
 Stiefel-Whitney class} $\tiw_2(f)$ is defined as the slash-product
$w_2(T(X,f))\msmall{/}[S^1]$. This is an element of\/
$\sfh^1(X,\zz_2)$. 
\end{defi} 

Let $(X,\omega)$ be a symplectic manifold and $f:X\to X$ a symplectomorphism which
acts trivially in $\sfh^\bullet(X,\zz)$. Then the K\"unneth formula
$\sfh^\bullet(T(X,f),\zz)\cong\sfh^\bullet(X,\zz)\otimes\sfh^\bullet(S^1,\zz)$ holds over integers.
Denote by $V(X,f)$ the fibrewise tangent bundle to the projection
$T(X,f)\to S^1$. Thus $V(X,f)$ is the quotient of $TX{\times}[0,1]$ under the
identification $(x,v;1)\sim(f(x),df(v);0)$. The form $\omega$ on $X$ induces the
fibrewise symplectic structure on $V(X,f)$. In particular, the Chern classes
$c_i(V(X,f))$ are well-defined.

\newdefi{c1-f} The \slsf{secondary
 Chern class} $\tic_1(f)$ is defined as the slash-product\break 
$c_1(V(X,f))\msmall{/}[S^1]$. This is an element of\/
$\sfh^1(X,\zz)$. 
\end{defi} 

Let us consider some obvious properties of  classes $\tiw_2$ and $\tic_1$.

\begin{lem}\label{tiw2-easy} \sli For a symplectomorphism $f:X\to X$ acting
 trivially on $\sfh^\bullet(X,\zz)$ one has $\tiw_2(f)\equiv\tic_1(f)\mod2$.

\slii $\tic_1(f\circ g)=\tic_1(f)+\tic_1(g)$ ($\tiw_2(f\circ g)=\tiw_2(f)+\tiw_2(g)$)
for symplectomorphisms (reps.\ diffeomorphisms) $f,g:X\to X$ acting trivially on
$\sfh^\bullet(X,\zz)$ (resp.\ on $\sfh^\bullet(X,\zz_2)$)

\sliii {\sl (Restriction formula)} Let $Z\subset X$ be a connected
cooriented hypersurface and $f:X\to X$ a diffeomorphism stabilising $Z$,
preserving its coorientation, and acting trivially on $\sfh^\bullet(X,\zz_2)$ and
$\sfh^\bullet(Z,\zz_2)$. Then $\tiw_2(f|_Z)=\tiw_2(f)|_Z$.

\sliv Let $\gamma:=S^1$ be a circle and $X_\gamma:=\gamma{\times}S^2$ an $S^2$-bundle. Then there
exists two isotopy classes of fibrewise trivialisations
$\varphi:X_\gamma\xrar{\,\cong\,}\gamma{\times}S^2$. Moreover, fibrewise map $F:X_\gamma\to X_\gamma$ identical on
$\gamma$ and preserving the orientation of the fibre interchanges the isotopy
classes of fibrewise trivialisations if and only if $\tiw_2(F)\neq0$.
\end{lem} 

\state Remark. Let us note in this connection that $X_\gamma$ admits exactly two
non-isomorphic spin structures, and $F:X_\gamma\to X_\gamma$ as above exchanges spin
structures if and only if $\tiw_2(F)\cdot[\gamma]\neq0$, \ie, when $F$ exchanges also the
classes of trivialisations. The proof is omitted since we shall not use this
fact.

\proof \slip, \sliip, and \sliii follow easily from the definitions.

\smallskip\noindent%
\sliv The group $\diff_+(S^2)$ is homotopy equivalent to $\SO(3)$. Since
$\pi_1(\SO(3))\cong\zz_2$ there exist two isotopy classes of trivialisations
$\varphi:X_\gamma\cong\gamma{\times}S^2$, and any fibrewise diffeomorphism $F:X_\gamma\to X_\gamma$ as in the
hypotheses either preserves or interchanges these classes of trivialisations
$\varphi:X_\gamma\cong\gamma{\times}S^2$. In the first case $T(F,X_\gamma)$ is diffeomorphic to $T^2{\times}S^2$
and $w_2(T^2{\times}S^2)=0$, in the second case $T(F,X_\gamma)$ is the non-trivial
$S^2$-bundle $T^2{\tix}S^2$. This manifold is diffeomorphic to a complex
minimal ruled surface $Z$ with a ruling $\pr:Z\to E$ over a complex elliptic
curve $E$ admitting disjoint holomorphic sections $S_1,S_{-1}$ of
self-intersection $+1$ and respectively $-1$. The first Chern class of $Z$ is
Poincar\'e dual to $[S_1]+[S_{-1}]$, and its reduction modulo $2$ gives
$w_2(Z)=w_2(T^2{\tix}S^2)$. In particular, the evaluation of $w_2(Z)$ on each
section is always non-trivial. As a consequence, we conclude that $F:X_\gamma\to X_\gamma$
as above exchanges the classes of trivialisations if and only if
$\tiw_2(F)\cdot[\gamma]\neq0$. Here $[\gamma]\in\sfh_1(\gamma,\zz_2)\cong\sfh_1(X_\gamma,\zz_2)$.
\qed

\medskip%
Now we prove the first assertion of {\slsf{Theorem 1}}.

\begin{prop}\label{tic1-nontriv} Let $(X,\omega)$ be an irrational ruled symplectic
 $4$-manifold and $\pr:X\to Y$ some singular symplectic ruling.  Then there
 exist a symplectomorphism $F:X\to X$ with any prescribed value $\tic_1(F)$
 which acts trivially on $\pi_1(X)$ and $\sfh_2(X,\zz)$.
\end{prop}

\state Remark. I want to use the occasion and thank once more Denis Auroux who
pointed me out on the existence of symplectomorphisms $F$ with non-trivial
$\tiw_2(F)$.  The explicit example of such a symplectomorphism constructed in
the proof coincides essentially with the one communicated to me by Denis
Auroux.

\proof First, we construct a symplectomorphism $F:X\to X$ with a given
$\tic_1(F)$ for the case of minimal irrational ruled symplectic $4$-manifolds.
Let $Y$ be a closed oriented surface of genus $g$ and $\pi:\bfp\to Y$ a principle
$\SO(2)$-bundle. Equip $\bfp$ with a connection, let $\Theta$ be the connection
form and $\calf$ the curvature of the connection. Recall that $\calf$ is a
$2$-form on $Y$ such that $\pi^*\calf=-d\Theta$. Let $(\rho,\theta)$ be the polar coordinates on
$S^2$ such that the standard round metric is given by
$ds^2=d\rho^2+\sin^2(\rho)d\theta^2$ and the volume form is $\sin(\rho)d\rho\land d\theta$. Define the
rotation action of $\SO(2)$ on $S^2$ by the formula $(\rho,\theta)\mapsto(\rho,\theta+\varphi)$. Then $\rho$
is an $\SO(2)$-invariant function. 

Set $X_0:=\bfp\times_{\SO(2)}S^2$. Then $X_0$ is an $S^2$-bundle over $Y$ and hence
a ruled $4$-manifold. Fix a symplectic form $\omega_Y$ on $Y$. Let us use the same
notation $\rho,\Theta,\calf,\omega_Y$ for the pull-backs of these objects on $X_0$ (notice
that they are well-defined). Finally, set
\[\textstyle
\omega:=\omega_Y+\sin(\rho)d\rho\land\Theta-2\,\sin^2(\frac\rho2)\,\calf.
\]
A direct computation shows that $\omega$ is a symplectic form on $X_0$ if and only
if $\omega_Y-t\calf$ is a symplectic form on $Y$ for every $t\in[0,2]$.  Notice that
the restriction of $\omega$ on each fibre of the natural projection $\pr:X_0\to Y$ is
the standard form $\sin(\rho)d\rho\land d\theta$. In particular, this gives us an explicit
example of a symplectic ruling on $\pr:X_0\to Y$. $X_0$ is diffeomorphic to the
product $Y\times S^2$ in the case when $c_1(\bfp)$ is even and to $Y{\tix}S^2$
otherwise.

Let $\alpha\in\sfh^1(Y,\zz)\cong\sfh^1(X_0,\zz)$ be a primitive class. Then $\alpha$ is
Poincar\'e dual to an embedded non-separating curve $\delta\subset Y$.  Let $U$ be a
bi-collar neighbourhood of $\delta$ and $s$ a local coordinate on $Y$ defined in
some neighbourhood of the closure $\barr U$ such that $\delta$ is given by $\{s=0\}$,
such that for every $a\in]-1,1[$ the set $\{s=a\}$ is an embedded curve, and such
that $U$ is given by $\{-1<s<1\}$. Let $\chi(s)$ be a smooth function which is
identically $0$ for $s\leq-\half$ and identically $1$ for $s\geq\half$. Set
$X_U:=\pr\inv(U)\subset X_0$. Use the same notation $s$ for the composition $s\circ\pr$,
this is now a local (coordinate) function on $X_U$. Define
$H:=\chi(s)\cdot2\sin^2(\frac\rho2)$, consider it as a Hamiltonian function on $X_U$,
and let $\Phi_\tau:X_U\to X_U$ be the Hamiltonian flow generated by $H$ with the time
coordinate $\tau$. Then $\Phi_\tau$ is trivial (constant) in the set $\{s\leq-\half\}$ and
coincides with the Hamiltonian flow $\Psi_\tau$ of the function $2\sin^2(\frac\rho2)$
in the set $\{s\geq\half\}$. At this point we make an easy but important
observation that the Hamiltonian flow $\Psi_\tau$ of $2\sin^2(\frac\rho2)$ is simply
the rotation action of $\SO(2)$, namely, $\Psi_\tau:(y,\rho,\theta)\mapsto(y,\rho,\theta+\tau)$.

Now define the symplectomorphism $F:X_0\to X_0$ setting $F:=\id$ in $X_0\bs X_U$
and $F:=\Phi_{2\pi}$ in $X_U$. We claim that $F$ acts trivially on $\pi_1(X)$ and
$\sfh_2(X,\zz)$ . Indeed, the set $Y^*:=\{\rho=0\}$ is a section
of $\pr:X_0\to Y$ such that the maps $\Phi_\tau,\Psi_\tau,F$ are identical on $Y^*$. Thus $F$
acts trivially on the fundamental group $\pi_1(Y)$ and on the homology class
$[Y^*]$. Besides, $F$ clearly preserves the homology class of the fibre
$\pr\inv(y)$. 

Finally, let us show that $\tic_1(F)=\alpha$. Clearly, $\tic_1(F)\cdot[\gamma]=0$ for any
curve $\gamma$ disjoint from $U$. This means that $\tic_1(F)$ is a multiple of $\alpha$,
$\tic_1(F)=m{\cdot}\alpha$. Thus $\langle\tic_1(F),[\gamma]\rangle=m{\cdot}\langle\alpha,[\gamma]\rangle$ for every closed curve
$\gamma\subset Y$. Fix an embedded curve $\gamma$ intersecting $\delta$ transversally in a single
point. Let $Z:=T(X,F)$ be the map torus of $F$. Since $F$ is identical on the
section $Y^*=\{\rho=0\}\subset X$, the map torus $Z^*:=T(Y^*\!\!,F)$ is naturally
isomorphic to the product $Y^*{\times}S^1$. The restricted tangent bundle
$TX_0|_{Y^*}$ has a natural decomposition in to the sum $TY^*\oplus V_\pr$ where
$V_\pr$ is the vertical tangent bundle to the projection $\pr:X_0\to Y$.  Recall
that the bundle $V(X_0,F)$ is the vertical tangent bundle of the projection
$Z=T(X,F)\to S^1$, and the fibre of $V(X_0,F)$ in a point $(x,t)\in Z$ is the
tangent plane $T_xX_0$. It follows from the construction, that over $Y^*$ the
map $F$ preserves the decomposition $TX=TY^*\oplus V_\pr$. Hence the restricted
bundle $V(X_0,F)|_{Z^*}$ admits a similar natural decomposition into the sum
$TY^*\oplus\wt V_\pr$ in which the fibre of $\wt V_\pr$ over $(y,t)\in Z^*=Y^*{\times}S^1$
coincides with the fibre of $V_\pr$ over $y\in Y^*$. Since both subbundles
$TY^*$ and $\wt V_\pr$ are symplectic,
$c_1(V(X_0,F)|_{Z^*})=c_1(Y^*)+c_1(\wt{V}_\pr)$. In the latter formula we
identify $c_1(Y^*)\in\sfh^2(Y,\zz)$ with the corresponding summand in the
K\"unneth formula for $Z^*=Y^*{\times}S^1$. It follows that the slash product
$c_1(Y^*)\msmall{\slash}[S^1]$ vanishes.

The calculation of $c_1(\wt{V}_\pr)\msmall{\slash}[S^1]$ is done as follows: The map
$F$ induces a linear isomorphism of $V_\pr$ preserving the orientation of
$V_\pr$. The automorphism bundle of $V_\pr$ is a trivial(isable) principle
bundle with the structure group $\Gl_+(2,\rr)$ of orientation preserving
linear isomorphisms of $\rr^2$. Since $\Gl_+(2,\rr)$ is homotopically
equivalent to $S^1$, the isotopy class of the isomorphism $F:V_\pr\to V_\pr$ is
identified with the induced homotopy class $[F]$ in the group $[Y^*,S^1]$ of
homotopy classes of maps between $Y^*$ and $S^1$. The homotopy $[Y^*,S^1]$ is
canonically isomorphic to $\sfh^1(Y^*,\zz)$, and moreover, the image of
$[F]\in[Y^*,S^1]$ under this isomorphism is the class
$c_1(\wt{V}_\pr)\msmall{\slash}[S^1]$.

Let $\gamma^*$ be the lift of the curve $\gamma\subset Y$ chosen above to the section $Y^*$.
The choice of the Hamiltonian function $H$ above was made to ensure that the
flow $\Phi_\tau$ (where it is defined) acts on $V_\pr$ as the rotation with the
angle velocity $\chi(s)$. Consequently, the map $F$ has the following action on
the restricted bundle $V_\pr|_{\gamma^*}$: This action is trivial where $\chi(s)$ is
constant, and otherwise it is rotation by angle $2\pi\chi(s)$. (Here, in abuse of
notation, we denote by $\chi(s)$ its restriction on $\gamma^*$.)  Hence $F$ makes a
single twist of $V_\pr$ along $\gamma^*$ and  $\tic_1(F)=\alpha$ as
desired.

\smallskip%
The above construction proves the proposition in the special case
$X_0=Y{\times}S^1$. The general case follows from the blow-up technique and
\propo{ruled1}.
\qed

\newsubsection[barMap]{Mapping class groups and symplectomorphisms} %
The results proved in this subsection will be applied in the assertion \sliv
of \slsf{Theorem 3.}  We use some properties of the mapping class groups
$\map(Y,\ell)$ (see \refdefi{def-barMap}) for which we refer to
\cite{Bir,Ivn,F-M} and \cite{Ger}.

\begin{prop}\label{prop-MapY1} Let $(X,\omega)$ be a ruled symplectic $4$-manifold,
 $\pr:X\to Y$ a symplectic ruling, $\mbfy=\{y^*_1,\ldots,y^*_\ell\}$ the set of singular
 values of the ruling $\pr$, and $f\in\diff_+(Y,\mbfy)$ a diffeomorphism. Then
 there exists a diffeomorphism $f_1:Y\to Y$ isotopic to $f$ relative
 $\mbfy=\{y^*_1,\ldots,y^*_\ell\}$, a fibrewise diffeomorphism $F_1:X\to X$ such that
 $\pr\circ F_1=f_1\circ\pr$, and symplectomorphism $F:X\to X$ smoothly isotopic to
 $F_1$.
\end{prop} 

\proof Recall that the group $\map(Y,\mbfy)=\diff_+(Y,\mbfy)/\diff_0(Y,\mbfy)$
is generated by Dehn twists along embedded curves $\gamma\subset Y\bs\mbfy$. 
Since we are interested only on isotopy classes of diffeomorphisms, it is
sufficient to prove the proposition in the special case when $f:Y\to Y$ is such a
Dehn twist.

Let $\gamma\subset Y$ be an embedded curve avoiding the points $\mbfy$. Denote
$X_\gamma:=\pr\inv(\gamma)$.  Consider first the special case when there exists a
neighbourhood $U\subset Y$ of $\gamma$ such that the projection $\pr:\pr\inv(U)\to U$ is a
trivial $S^2$ bundle, so that $\pr\inv(U)\cong U{\times}S^2$ and that in the set
$\pr\inv(U)$ the symplectic structure $\omega$ is a product of structures $\omega_Y$ on
$U$ and $\omega_{S^2}$ on $S^2$. Then we can isotopically deform the map $f:Y\to Y$
into a map $f_1:Y\to Y$ with the support in $\pr\inv(U)$, and then set
$F_1:=f_1{\times}\id_{S^2}$ inside $\pr\inv(U)$ and $F_1:=\id_X$ outside
$\pr\inv(U)$.

The general case follows from \lemma{pi0-dfrm}. Indeed, we start with the
product $Y{\times}S^2$, equip it with the product symplectic structure, and then
make appropriate blow-ups and deformations of the symplectic structure.
\qed

\smallskip
Using previous lemma we define a homomorphism from $\map(Y,\ell)$ to $\Gamma_\lor(X)$
where $\ell$ is the maximal number of pairwise disjoint exceptional symplectic
spheres on $X$.  Let $\pr:X\to Y$ be a symplectic ruling with only ordinary
singular fibres over points $y^*_i$. Fix disc neighbourhoods $D_i\subset Y$ of points
$y^*_i$, set $Y^\circ:=Y\bs\cup_iD_i$, $X^\circ:=\pr\inv(Y^\circ)$, and fix a topological
trivialisation $X^\circ\cong Y^\circ{\times}S^2$ compatible with the projection $\pr:X^\circ\to Y^\circ$.
Represent any element $\phi\in\map(Y,\ell)$ by a diffeomorphism $f:Y\to Y$ with support
in $Y^\circ$, set $F:=f{\times}\id_{S^2}$ inside $X^\circ$ and $F:=\id$ outside $X^\circ$. Now
we define the desired homomorphism $\map(Y,\ell)\to\Gamma_\lor(X)$ associating to $\phi$ the
isotopy class of the diffeomorphism $F$ described just above.

\begin{prop}\label{ImMapYk} The image of\/ $\map(Y,\ell)$ in $\Gamma_\lor(X)$ is isomorphic
 to  $\barr\map(Y,\ell)$.
\end{prop} 

\proof We must describe the kernel of the homomorphism from
$\map(Y,\ell)\to\Gamma_\lor(X)$. Denote this kernel by $\sfk$. Recall that $\pi_1(X)\cong\pi_1(Y)$
and that the group $\map(Y)$ can is isomorphic to the \slsf{outer automorphism
 group} $\out_+(\pi_1(Y,y^*))$ of $\pi_1(Y)$, \ie, the quotient of the
automorphism group $\aut_+(\pi_1(Y,y^*))$ by the group $\pi_1(Y,y^*)$ acting on
itself by conjugation. Composing the homomorphism $\map(Y,\ell)\to\Gamma_\lor$ with the
action of $\Gamma_\lor$ on the fundamental group, we see that the kernel of
$\map(Y,\ell)\to\Gamma_\lor$ lies in the kernel of $\map(Y,\ell)\to\map(Y)$.

Consider first the case $g(Y)\geq2$. Recall that in this case the kernel of
$\map(Y,\ell)\to\map(Y)$ is the pure braid group $\pbr(Y,\ell)$. It is generated by
the following moves: a single point $y^*_i\in\mbfy$ moves along an embedded
closed curve $\gamma$ on $Y$ avoiding the rest of the points $\mbfy$ which stay
without move.  A special case of such a move is as follows: For an embedded
arc $\alpha$ connecting two points $y^*_i,y^*_j\in\mbfy$ and avoiding the rest of the
points $\mbfy$, the point $y^*_i$ approaches to $y^*_j$ along $\alpha$, turns
around it once, and then comes back also along $\alpha$. It is not difficult to see
that this move is isotopic to the square $\tau_\alpha^2$ of the half-twist $\tau_\alpha$ along
the arc $\alpha$, and so  $\tau_\alpha^2$ is isotopic to the Dehn twist $T_\gamma$ along the
closed curve $\gamma$ surrounding the arc $\alpha$. As we have shown in
\refsubsection{dehn} the image of such a squared half-twist $\tau_\alpha^2$ is trivial
in $\Gamma_\lor$.

The triviality of squared half-twists $\tau_\alpha^2$ allows us to ignore all
braidings/interlacings of strands in the pure braid group $\pbr(Y,\ell)$ and take
into account only the homotopy classes of strands of braids. This means that
all squared half-twists $\tau_\alpha^2$ generate the subgroup which is the kernel of
the natural projection $\pbr(Y,\ell)\to(\pi_1(Y,y^*))^\ell$. Further, assume that all
points $\mbfy$ lie in some small disc $D$. We can move the disc $D$ with the
points $\mbfy$ on it along arbitrary path $\gamma$ on $Y$ applying an appropriate
map $F_\gamma$ which is isotopic to the identity and hence lies in $\diff_0(Y)$.
It is not difficult to see isotopy class of $F_\gamma$ in $\diff_+(Y,\ell)$ is an
element $\varphi_\gamma\in\pbr(Y,\ell)$ whose projection to $(\pi_1(Y,y^*))^\ell$ is the diagonal
element $([\gamma],\ldots,[\gamma])\in(\pi_1(Y,y^*))^\ell$ where $[\gamma]$ is the homotopy class of $\gamma$
in $\pi_1(Y,y^*)$. Thus we have shown that the kernel
$\sfk=\ker\big(\map(Y,\ell)\to\Gamma_\lor\big)$ lies in the pure braid group $\pbr(Y,\ell)$,
contains the kernel of the projection $\pbr(Y,\ell)\to(\pi_1(Y,y^*))^\ell$, and the
image of $\sfk$ in $(\pi_1(Y,y^*))^\ell$ contains the diagonal embedding of
$\pi_1(Y,y^*)$. 

Denote by $\sfk_1$ the subgroup in $(\pi_1(Y,y^*))^\ell$ normally generated by the
diagonal embedding of $\pi_1(Y,y^*)$. We claim that $\sfk_1$ coincides with the
kernel the projection $(\pi_1(Y,y^*))^\ell\to\big(\sfh_1(Y,\zz)\big)^\ell/\sfh_1(Y,\zz)$
where $\sfh_1(Y,\zz)$ is embedded in $\big(\sfh_1(Y,\zz)\big)^\ell$ in the
diagonal way. First we notice that $\sfk_1$ contains the kernel of the
projection $(\pi_1(Y,y^*))^\ell\to\big(\sfh_1(Y,\zz)\big)^\ell$. This fact follows
directly from the property that $\sfh_1(Y,\zz)$ is the Abelianisation of
$\pi_1(Y,y^*)$. Now the claim follows from the property that the projection
$(\pi_1(Y,y^*))^\ell\to\big(\sfh_1(Y,\zz)\big)^\ell$ maps the diagonal embedding of
$\pi_1(Y,y^*)$ onto the diagonal embedding of $\sfh_1(Y,\zz)$.

The description of the subgroup $\sfk_1<(\pi_1(Y,y^*))^\ell$ made above implies
that the whole subgroup $\sfk(Y,\ell)<\pbr(Y,\ell)$ defined in \refdefi{def-barMap}
on \pageref{def-barMap} lies in the kernel of the homomorphism
$\map(Y,\ell)\to\Gamma_\lor(X)$. To show that this is the whole kernel of
$\map(Y,\ell)\to\Gamma_\lor(X)$ we consider the action of $\map(Y,\ell)$ on the second
fundamental group of $X$ $\pi_2(X,x^*)$. Describing this group we use the
following two facts: The first is that the higher homotopy groups of the
manifold $X$ and its universal covering $\wt X$ are isomorphic, and in
particular %
$\pi_2(X,x^*)\cong\pi_2(\wt X,\ti x^*)$ for any choice of the lift $\ti x^*$ of a
given base point $x^*$. The second fact is that for the simply connected
manifold $\wt X$ the natural Hurewicz homomorphism $\pi_2(\wt X,\ti
x^*)\to\sfh_2(\wt X, \zz)$ %
is an isomorphism.

Observe that the universal covering $\wt Y$ is topologically a
disc, and that the universal covering $\wt X$ admits the natural ruling %
$\wt\pr:\wt X\to\wt Y$ with the following singular fibres: for any point $y^*_i$
we fix one its lift $\ti y^*_i$ to $\wt Y$. Then all remaining lifts have the
form $\gamma{\cdot}\ti y^*_i$ with $\gamma\in\pi_1(Y,y^*)$, the whole collection %
$\{\gamma{\cdot}\ti y^*_i\}$ is the set of singular values of the ruling %
$\wt\pr:\wt X\to\wt Y$ with an ordinary singular fibre over every such point.

It follows that the group $\pi_2(X,x^*)\cong\pi_2(\wt X,x^*)\cong\sfh_2(\wt X, \zz)$ is a
free Abelian group with the basis naturally indexed by a single index $F$
corresponding to the class of a generic fibre and by the set $\{\gamma{\cdot}\ti
y^*_i\}$, %
such that each point $\gamma{\cdot}\ti y^*_i$ indexes the homotopy class of an
appropriately chosen component of the singular fibre over this
point. Moreover, the natural action of $\pi_1(X,x^*)\cong\pi_1(Y,y^*)$ on $\pi_2(X,x^*)$
is given by deck transformations of $\wt X$ and can be described as follows:
An element $\delta\in\pi_1(X,x^*)$ acts trivially on $F$ and transforms the homotopy
class indexed by $\gamma{\cdot}\ti y^*_i$ in the class indexed by $\delta{\cdot}\gamma{\cdot}\ti
y^*_i$. %
In particular, this action reduces to permutation of the basis described
above. Further, the action of the pure braid group $\pbr(Y,\ell)$ on $\pi_2(X,x^*)$
reduces to shifting of homotopy classes indexed by $\gamma{\cdot}\ti y^*_i$. Thus the
braid $\beta$ whose projection to $(\pi_1(Y,y^*))^\ell$ is an element $(\delta_1,\ldots,\delta_\ell)$
acts on $\gamma{\cdot}\ti y^*_i$ as $\delta_i{\cdot}\gamma{\cdot}\ti y^*_i$. It follows that we can
recover the whole projection $(\delta_1,\ldots,\delta_\ell)$ of the pure braid $\beta$ as above
tracing the action of $\beta$ on the second homotopy group $\pi_2(X,x^*)$.

The latter consideration must be corrected taking into account that we can
move the base point $x^*$ which was considered above to be fixed. Clearly, the
shift of the base point $x^*$ along a path $\sigma$ corresponds to the
multiplication of a tuple $(\delta_1,\ldots,\delta_\ell)$ as above by the diagonal element
$(\sigma\inv,\ldots,\sigma\inv)$. This last argument finishes the proof of the proposition in
the case $g(Y)\geq2$.

The whole argumentation remains essentially the same in the case $g(Y)=1$ when
$Y$ is the torus $T^2$. The only difference comes from the fact, distinguished
already in \refdefi{def-barMap}, that the kernel of the natural projection
$\pbr(T^2,\ell)\to(\pi_1(T^2,y^*))^\ell$ is not the pure braid group $\pbr(T^2,\ell)$ but
its quotient $\pbr(T^2,\ell)/\sfh_1(T^2,\zz)$. We leave the details to the interested
reader.
\qed

\smallskip \newsubsection[proof1]{Proof of {\sl \textsf{Theorem 1}}.}  %
In this paragraph we assume that the hypotheses \slii of 
the theorem are fulfilled.

\smallskip \noindent
 {\sl Step 1. Deforming a symplectomorphism into a fibrewise map.} %
Let $J$ be a complex structure on the sphere $S^2$ and
$p_0,p_1,p_2\in S^2$ three distinct points. Then there exists a unique
$J$-holomorphic isomorphism $\psi:(S^2,J)\to\cp^1$ with $\psi(p_0)=0,\psi(p_1)=1$ and
$\psi(p_2)=\infty$.  Consequently, if $J$ is an almost complex structure on $X$,
$\pr:X\to Y$ a $J$-holomorphic ruling, and $s_0,s_1,s_2:Y\to X$ three sections of
the ruling, then the fibre $\pr\inv(y)$ is equipped a distinguished
trivialisation provided $y$ is not a critical value of the ruling and the
points $s_0(y),s_1(y),s_2(y)\in\pr\inv(y)\cong S^2$ are distinct. 

Now assume that we have two configurations $(X,\pr,J,s_i)$ and
$(X',\pr',J',s'_i)$ of type ``ruling + structure + 3~sections'' with the
following properties:
\begin{itemize}
\item sections $s_0,s_1,s_2$ and $s'_0,s'_1,s'_2$ have only positive
 transversal intersections;
\item projections $\pr:X\to Y$ and $\pr':X'\to Y$ have only ordinary singular
 fibres which do not contain intersection points of the sections $s_0,s_1,s_2$
 and $s'_0,s'_1,s'_2$;
\end{itemize}
Further, assume that they are combinatorially equal. The latter means that
there exists a bijection between the exceptional spheres $E_i$ and $E'_i$ of
$X$ and respectively $X'$, such that the corresponding intersection indices of
the sections and the exceptional spheres are equal.  Let $y^*_i,y'{}^*_i$ and
respectively $y^\times_i,y'{}^\times_j$ and be the projections of singular fibres of
$\pr$ and $\pr'$, and respectively the projections of intersection points of
the sections. Fix bijections between points  $y^*_i,y'{}^*_i$ and
respectively $y^\times_i,y'{}^\times_j$ which preserves the combinatorics: if $y^\times_i$ is
the projection of some intersection point of say $s_0$ and $s_1$, then 
$y'{}^\times_i$ must be the projection of an intersection point of $s'_0$ and $s'_1$. 
Then there exists a diffeomorphism $f:Y\to Y$ mapping $y^*_i$ onto
$y'{}^*_i$ and $y^\times_i$ onto $y'{}^\times_i$. 

We claim that in this case there exists a fibrewise diffeomorphism $\Psi:X\to X'$
covering $f$, \ie, $\pr'\circ \Psi=f\circ\pr$, which is unique up to
isotopy. Constructing $\Psi$ we first deform slightly the structures $J$ and
$J'$ such that they become integrable near singular fibres and near fibres
through intersection points of sections. Without loss of generality we may
assume that the fibres of rulings remain pseudoholomorphic. Next, we deform
sections such that they become holomorphic near their intersection
points. After this we perform the following surgeries on $X$ and $X'$:
\begin{itemize}
\item we blow up every intersection point of sections, this gives us manifolds
 $\wh X$ and $\wh X'$ with projections on $X$ and $X'$ and with the induced
 rulings and sections;
\item for every ordinary singular fibre on $\wh X$ and $\wh X'$, we contract
 that one exceptional component which meets none or only one lifted section
 $\hat s_i$ or $\hat s'_i$. 
\end{itemize} 
Let $\wt X$ and $\wt X'$ be the resulting manifolds $\wh X \to\wt X$,
$\wh{X}'\to\wt X'$ the induced projections, and $\ti s_i,\ti s'_i$ the sections.
Then $\ti s_i,\ti s'_i$ are pairwise disjoint and there exists a unique
diffeomorphism on $\wt\Psi:\wt{X}\to\wt X'$ which is holomorphic on fibres of
the induced rulings and maps the sections in the corresponding sections. We
need to correct $\wt \Psi$ since it does not maps the blow-up centres of
$\wh{X}\to\wt X$ on those of $\wh{X}'\to\wt X'$, but this can be done easily by an
appropriate isotopy. After this correction, the inverse surgery transforms
$\wt\Psi$ in the desired fibrewise diffeomorphism $\Psi:X\to X'$.

Observe that we have shown a stronger property: the data  $(X,\pr,J,s_i)$,
$(X',\pr',J',s'_i)$, and $f:Y\to Y$ determine the isotopy class of the
constructed diffeomorphism $\Psi$. 

\smallskip %
Now we can explain the main idea of this step: We construct families
$(J_t;\pr_t;s_i(t))$ of the data as above using Gromov's technique
of pseudoholomorphic curves, and such that $J_1=F_*J_0$ and $\pr_1= F_*\pr_0$.
Then using this family we deform $F$ into a diffeomorphism fibrewise with
respect to the ruling $\pr_0$.

Using \lemma{pi0-dfrm} we may assume that $X$ is not minimal and contains
exceptional symplectic spheres. Let $g:=g(X)=g(Y)$ be the genus of $X$.  Equip
$Y$ with some generic complex structure. Consider a generic holomorphic line
bundle $L$ on $Y$ of degree $g$. Let $X_0$ be the fibrewise projective
completion of the total space of $L$. Then $X_0$ is the projectivisation of
the rank-$2$ holomorphic bundle $L\oplus\scro_Y$. Let $\pr:X_0\to Y$ be the
corresponding ruling, $Y_g$ the zero section of $\pr:X_0\to Y$, and $Y_{-g}$ the
infinity section.  Then $[Y_{\pm g}]^2=\pm g$. Blowing up $X_0$ in $\ell$ generic
points $x^*_1,\ldots,x^*_\ell$ we obtain a manifold diffeomorphic to $X$. This gives
an integrable complex structure $J^*$ on $X$.  It could occur that this
structure is not compatible with any symplectic form $\omega'$ cohomologous to our
form $\omega$. However, \lemma{pi0-dfrm} ensures that if we deform $\omega$, the
deformed form $\omega^*$ admits a symplectomorphism $F^*$ isotopic to $F$ and hence
having the same class, \ie, $\tiw_2(F)=\tiw_2(F^*)$. As such a deformation of $\omega$,
we deflate almost completely the exceptional symplectic spheres $E_1,\ldots,E_\ell$
obtained above and increase the volume of the base $Y$, until we achieve the
relation $\lan[\omega^*],[Y_{-g}]\ran>0$. Now the structure $J^*$ would be tamed by
some symplectic form of the cohomology class $[\omega^*]$. Without loss of
generality we may assume that $J^*$ is tamed by $\omega$.

Now choose a generic almost complex structure $J_0$ sufficiently close to
$J^*$ and set $J_1:=F_*(J_0)$. Then $J_1$ is also an $\omega$-tamed almost complex
structure, and so there exists a path $J_t$ of $\omega$-tamed structures connecting
$J_0$ and $J_1$. We choose the path $J_t$ to be generic enough. Then every
$J_t$ induces a ruling on $X$. Without loss of generality we may assume that
they have the same base $Y$. Regularising the family $J_t$ (see
\refdefi{regularise}) we may also assume that the rulings
$\pr_t:X\to Y$ are smooth. Let $y^*_1(t),\ldots,y^*_\ell(t)$ be the images of singular
fibres of $\pr_t:X\to Y$, and $\mbfy^*(t)$ the whole set.

By the choice of $L$, $\sfh_\dbar^0(Y,L)\cong\cc$ and $\sfh^1_\dbar(Y,L)=0$.
Consequently, $X_0$ admits complex $1$-dimensional family of holomorphic
sections of the ruling $\pr:X_0\to Y$. Moreover, for a generic point $y_0\in Y$ and
a generic point $x_0\in\pr\inv(y_0)$ there exists a unique $J^*$-holomorphic
section $C_0$ of the ruling $\pr:X_0\to Y$ passing through $x_0$. This property
retains for structures $J$ close to $J^*$. Choose a generic family $y_0(t)$ of
points on $Y$, and three generic families $x_0(t),x_1(t),x_2(t)$ of points
lying on the fibre $\pt_t\inv(y_0(t))$, such that $F(y_0(0))=y_0(1)$ and
$F(x_i(0))=x_i(1)$. Consider the moduli space $\calm$ of
quadruples $(t;C_0,C_1,C_2)$ such that $t\in[0,1]$ and $C_i$ are irreducible
$J_t$-holomorphic curves in the homology class $[Y_g]$ each passing through
the corresponding point $x_i(t)$. Then $\calm$ is a smooth manifold of
expected dimension $\dim_\rr\calm=1$ equipped with the natural projection
$\pi_\calm:\calm\to I:=[0,1]$. The choice of the involved data ensures that for a
given $t$ close to $0$ or to $1$ there exists a unique triple
$(C_0(t),C_1(t),C_2(t))$ in the fibre $\pi_\calm\inv(t)$.

Now let us observe that $\calm$ is compact and $\pi_\calm:\calm\to I$ is proper.
Indeed, if one of the curves $C_i(t)$ breaks, then the expected dimension of
the arising configuration would be negative. This eventuality is
prohibited by the generic choice of the path $J_t$.  It follows that $\calm$
contains a unique component $\calm^*$ which is an interval, whereas all
remaining components are circles whose $\pi_\calm$-projections lie in some
proper subinterval $[\varepsilon,1-\varepsilon]$, $\varepsilon>0$. Choose a parameter (coordinate) $\tau$ on
$\calm^*$ varying in the interval $I=[0,1]$. Then the dependence of the old
parameter $t$ on $\tau$ is established by the projection $\pi_\calm$. We write this
dependence in the form $J_\tau$ instead of formally correct $J_{\pi_\calm(\tau)}$ and
so on.

The genericity condition implies that for each $\tau$ every pair of the curves
$C_i(\tau),C_j(\tau)$ has $g$ transversal intersection points which depend smoothly
on $\tau$. In this way we obtain $3g$ families $x^\times_j(\tau)$ of nodal points. We
denote by $y^\times_j(\tau):=\pr_\tau(x^\times_j(\tau))$ their projections on $Y$, they are
pairwise distinct. Let $\mbfx^\times(\tau)$ and $\mbfy^\times(\tau)$ be the whole
collections. Then there exists a family $f_\tau:Y\to Y$ of diffeomorphisms such that
$f_0=\id_Y$ and such that that $f_\tau(\mbfy^*(0))=\mbfy^*(\tau)$ and
$f_\tau(\mbfy^\times(0))=\mbfy^\times(\tau)$.  This gives us the desired family of structures
$(J_\tau,\pr_\tau;C_0(\tau),C_1(\tau),C_2(\tau);f_\tau)$. The construction above shows that
there exists a family $\Psi_\tau:X\to X$ of diffeomorphisms such that $\Psi_0=\id_X$,
$\Psi_\tau(C_i(0))=C_i(\tau)$, and $\pr_\tau\circ \Psi_\tau=f_\tau\circ\pr_0$. We rewrite the latter
relation in the form $f_\tau\inv\circ\pr_\tau^{\ }=\pr_0\vph\circ\Psi_\tau\inv$.

The uniqueness of $J$-holomorphic rulings for a given structure $J$ ensures
that $\pr_1\circ F=\pr_0$.  Set $F_\tau:=\Psi_\tau\inv\circ F$.  Then $F_\tau$ is an isotopic
deformation of $F=F_0$ and
$\pr_0\circ F_1=\pr_0\circ\Psi_1\inv\circ F=f_1\inv\circ\pr_1\circ F=f_1\inv\circ\pr_0$. So setting
$f:=f_1\inv:Y\to Y$ we obtain the desired relation 
$\pr_0\circ F_1=f\circ\pr_0$.

\smallskip \noindent {\sl Step 2. Dealing with braid and mapping class
 groups.} %
Consider the constructed map $f:Y\to Y$. Denote by $y^*_i$ (instead of
$y^*_i(0)$) projections of singular fibres of $\pr_0$, and by $\mbfy^*$ the
whole collection.  The triviality of the action of $F$ on $\pi_1(X)$ implies
that $f$ acts trivially on $\pi_1(Y)$. Consequently, $f$ is smoothly isotopic to
identity in the group $\diff_+(Y)$ (see \eg., \cite{F-M}, Theorem
2.3). Further, by the construction $f$ preserves the projections $y^*_i$ of
singular fibres of $\pr_0$. Moreover, since $F$ acts trivially in homology, so
does $F_1$, and hence $F_1$ does not permute singular fibres itself and their
components. Consequently, every point $y^*_i$ is fixed by $f$. It follows that
$f$ defines an element in the pure braid group $\pbr(Y,\ell)$ of the
surface $Y$ where $\ell$ denotes the number of singular fibres.

Without loss of generality we can suppose that $f:Y\to Y$ is identical in some
small disc neighbourhood $D_i$ of each point $y^*_i$. As in the definition of
the homomorphism $\map(Y,\ell)\to\Gamma_\lor(X)$ on \pageref{prop-MapY1} set
$Y^\circ:=Y\bs\cup_iD_i$, $X^\circ:=\pr_0\inv(Y^\circ)$, fix a topological trivialisation
$X^\circ\cong Y^\circ{\times}S^2$ compatible with the projection $\pr_0:X^\circ\to Y^\circ$, and define the
diffeomorphism $F_2$ setting $F_2:=f{\times}\id_{S^2}$ inside $X^\circ$ and $F_2:=\id$
outside $X^\circ$.  Then the diffeomorphism $F_3:=F_2\inv\circ F^{\,}_1$ satisfies the
relation $\pr_0\circ F_3=\pr_0$. 

Recall that in \refsubsection{barMap} (see page \pageref{barMap}) we have
defined the homomorphism $\pbr(Y,\ell)\to\Gamma_\lor$, whose image is isomorphic to
$(\pi_1(Y,y^*))^\ell/\pi_1(Y,y^*)$, and such that the image of the isotopy class
$[f]$ of $f:Y\to Y$ is isotopy class of $F_2:X\to X$. Notice that in the proof of
\propo{ImMapYk} we show that the image of $[f]$ is completely 
determined how  $F_2$ acts on $\pi_1(X,x^*)$ and on $\pi_2(X,x^*)$. However,
since $\pr_0\circ F_1=\pr_0\circ F_2$, this action is the same for $F_1$ and $F_2$, and
hence is trivial by the assumption of the theorem. It
follows that the isotopy class of $F_2$ is trivial, and hence $F_1$ is
isotopic to the diffeomorphism $F_3$.

\smallskip \noindent %
{\sl Step 3. Discovering the secondary class $\tiw_2(F)$.} %
At the previous steps first we deformed a given topologically trivial
symplectomorphism $F$ into a diffeomorphism $F_3$, which acts fibrewisely with
respect to some ruling $\pr_0:X\to Y$. Making an additional isotopy, if needed,
we may suppose that $F_3$ is identical in a neighbourhood of every singular
fibre of $\pr_0$. Let $Y^\circ,X^\circ$ and the trivialisation $X^\circ\cong Y^\circ{\times}S^2$ be as
above. Then $F_3$ indices a smooth map $\wt F_3:Y\to\diff_+(S^2)$. Since $Y$ is
an Eilenberg-MacLane space, this map is homotopic to the identity if and only
if the induced homomorphism $(\wt F_3)_*:\pi_1(Y)\to\pi_1(\diff_+(S^2))$ is
trivial. Now observe that this homomorphism is exactly the secondary
Stiefel?Whitney class $\tiw_2(F_3)$. Thus $\wt F_3:Y\to\diff_+(S^2)$ is
homotopic to the identity iff $\tiw_2(F_3)=\tiw_2(F)$ vanishes. However, the
condition of the homotopic Al triviality of $\wt F_3:Y\to\diff_+(S^2)$ is equivalent
to the fact that $F_3$ is fibrewise isotopic to the identical map. 

The theorem follows.
\qed

\section[Structure of the diffeotopy group $\Gamma(X)$ \break of rational and
 ruled $4$-manifolds.]
{Structure of the diffeotopy group $\Gamma(X)\ $\\  of rational and
 ruled $4$-manifolds.} 
\label{Gamma}

In this section we give the proof of \slsf{Theorems 3 {\rm and} 4}. We start
with construction of certain diffeomorphisms representing those generators of
the group $\Gamma_\sfw$ which can not be represented by symplectomorphisms. Recall
that this is a subgroup of $\aut(\sfh_2(X,\zz))$ which is generated by 
diffeomorphisms preserving the orientation and the positive cone in
$\sfh_2(X,\rr)$.

\newsubsection[non-symp]{Some non-symplectic diffeomorphisms.}

\begin{lem}\label{except-twist}
Let $X$ be an oriented $4$-manifold and $E\subset X$ be an embedded sphere with
self-intersection $\pm1$. Then for any neighbourhood $U\subset X$ of $E$ there exists a
diffeomorphism $F$ of $X$ supported in $U$ such that $F(S)=S$ and such that
$F|_E:E\to E$ is an orientation inverting diffeomorphism. Moreover, one can
choose $F$ acting on $E$ as the antipodal map.

Furthermore, there exists an isotopy between  $F^2$ and the identity
map supported in $U$.
\end{lem}

\proof Inverting the orientation on $X$, if needed, we may suppose that
$[E]^2=-1$. Now without loss of generality we can assume that $E$ is the
exceptional sphere in $\cc^2$ blown-up in the origin and $U$ is a neighbourhood
of $E$. Let $(z_0,z_1)$ be the complex coordinates in $\cc^2$. Then
$[z_0:z_1]$ is the projective coordinate on $E$. Notice that the complex
conjugation of $\cc^2$ with $(z_0,z_1)\mapsto(\barz_0,\barz_1)$ indices an
orientation inverting diffeomorphism $[z_0:z_1]\mapsto[\barz_0:\barz_1]$ of $E$.
Write the coordinates $z_i$ as $x_i+\cpi y_i$. Let $\Phi(\varphi)$ be the rotation in
$y$-plane by the angle $\varphi$ preserving the $x$-plane, so that
\[
\Phi(\varphi)(x_0,x_1;y_0,y_1):=(x_0,x_1; 
\cos\varphi\, y_0 -\sin\varphi\,y_1, \sin\varphi\,y_0+\cos\varphi\,y_1).
\]
Finally, let $\chi(r)$ be a smooth cut-off function depending only on
$r:=\sqrt{|z_0|^2+|z_1|^2}$ which is supported in $U$ and such that $\chi(r)=1$
for all $r\in[0,\delta]$ with some $0<\delta\ll1$. Now the desired diffeomorphism $F$ is
given by the formula
\[
F(x_0,x_1;y_0,y_1):=\Phi(\pi{\cdot}\chi(r))(x_0,x_1;  y_0,y_1)
\]
where as above $r=\sqrt{|z_0|^2+|z_1|^2}$. 

\smallskip%
Further, let us notice that for a similar rotation map 
\[
\Psi(\varphi):(z_0,z_1)\mapsto (\cos\varphi \,z_0
-\sin\varphi\,z_1, \sin\varphi\,z_0+\cos\varphi\,z_1)
\]
and any function $f(r)$ supported in $U$ the map $G$ given by
$G(z_0,z_1):=\Psi(f(r))(z_0,z_1)$ is isotopic to identity and hence $G\circ F$ is
isotopic to $F$. For the special case $G:=\Psi(\frac{\pi}{2}{\cdot}\chi(r))(z_0,z_1)$ the
composition $F_1:=G\circ F$ acts on $E$ as the antipodal map $z\mapsto-\frac{1}{\barz}$.

\smallskip%
Finally, we show that $F^2$ is isotopic to the identity. Let us observe that
by construction the map $F$ has the form
\[
F(x_0,x_1;y_0,y_1):=\Phi(r)(x_0,x_1;  y_0,y_1)
\]
for some path $\Phi(r)$ in the orthogonal group $\SO(4)$, such that $\Phi(r)=\id$
for $r>\delta$ and such that $\Phi(0)$ is $\cc$-antilinear. For the square $F^2$ we
obtain a similar form
\[
F(x_0,x_1;y_0,y_1):=\Upsilon(r)(x_0,x_1;  y_0,y_1)
\]
with some path $\Upsilon(r)$ in the orthogonal group $\SO(4)$, such that $\Upsilon(r)=\id$
for $r>\delta$ and such that $\Upsilon(0)$ is $\cc$-linear and hence unitary,
$\Upsilon(0)\in\bfu(2)$.  Consequently, $\Upsilon$ defines an element $[\Upsilon]$ in the relative
homotopy group $\pi_1(\SO(4),\bfu(2))$. From the long exact sequence
\[
\cdots\to\pi_1(\bfu(2)) \to\pi_1(\SO(4))\to\pi_1(\SO(4),\bfu(2)) \to \pi_0(\bfu(2)) \to\cdots
\]
we see that the group $\pi_1(\SO(4),\bfu(2))$ is trivial. This implies the
existence of the desired isotopy $H_t$ between $F^2$ and $\id$ also in the
form $H_t(r)(x_0,x_1; y_0,y_1)$ for some family $H_t(r)$ of paths in $\SO(4)$
with $H_t(r)=\id$ for $r>\delta$ and with $H_t(0)\in\bfu(2)$.  
\qed 

\newdefi{E-twist}
 \sl We call a diffeomorphism $F_1$ as in the previous lemma and its isotopy
 class in $\Gamma(X)$ the \slsf{twist along an exceptional sphere $E$} (or along a
 $(\pm1)$-sphere). This is the counterpart of (symplectic) Dehn twists along
 (Lagrangian) $(-2)$-spheres for the case of $(-1)$-spheres.
\end{defi}

\smallskip

\begin{lem}\label{transvec}
 Let $W$ be the product $\Delta{\times}S^2$ blown-up in the point $(0,0)$, $E$ the
 arising exceptional sphere, $E'$ the proper pre-image of the central fibre
 $0{\times}S^2$, $\pr:W\to\Delta$ the projection, and $\bfe,\bfe',\bff$ the homology
 classes of $E$, $E'$, and respectively a generic fibre $S_y=\pr\inv(y)$.
 Further, let $W_2$ be the product $\Delta{\times}S^2$ blown-up in two points, $E_1,E_2$
 arising exceptional spheres, $\pr_2:W_2\to\Delta$ the projection, and
 $\bfe_1,\bfe_2,\bff$ the homology classes of $E_1$, $E_2$, and a generic
 fibre $S_y=\pr_2\inv(y)$.

 \sli There exists a diffeomorphism $F_1:W\to W$ mapping $\bfe$ onto $\bfe+\bff$
 and $\bff$ onto $\bff$ which preserves each fibre $S_y:=\pr\inv(y)$ over $y$
 lying on the boundary circle $\partial\Delta$. Moreover, $\tiw_2(F_1|_{\partial W})=1\in\sfh^1(\partial
 W,\zz_2)\cong\zz_2$.

 \slii There exists a diffeomorphism $F_2:W\to W$ mapping $\bfe$ onto
 $\bfe+2\bff$ and $\bff$ onto $\bff$ which is identical near the boundary $\partial
 W$.
\end{lem}

\proof \sli In view of the previous lemma, it is sufficient to find a
diffeomorphism $F'_1:W\to W$ mapping the class $\bfe$ onto $\bff-\bfe$. To find
such a map we observe that the manifold $W$ equipped with the natural complex
structure has two holomorphic contractions $\pi,\pi':W\to\Delta{\times}S^2$: the first
contracting the exceptional complex curve $E$ and the second contracting $E'$.
Let $z\in\Delta$ and $w\in\cp^1=S^2$ be the complex coordinates on $\Delta{\times}S^2$.  It
is easy to see that the desired map $F'_1:W\to W$ can be given by the formula
$F'_1(z,w)=(z,\frac{z}{w})$. 

It remains to calculate the class $\tiw_2(F'_1|_{\partial W})$. Let us compose $F'_1$
with the map $G:\partial W\to\partial W$ given by $G(z,w):=(z,\frac{1}{w})$. The map $G(z,w)$
extends in the obvious way from $\partial W=\partial(\Delta{\times}S^2)$ to the whole product
$\Delta{\times}S^2$. Therefore $\tiw_2(G)=0$ and $\tiw_2(F'_1|_{\partial W})
=\tiw_2(G\circ F'_1|_{\partial W})$. The latter composition is
$H:=G\circ F'_1(z,w)=(z,\frac{w}{z})$ and we can see the following picture: as $z$
goes along the circle $\partial\Delta$ the map $H$ rotate the sphere
counterclockwise. Therefore we conclude the desired equality
$\tiw_2(F_1|_{\partial W})=\tiw_2(H)=1\in\sfh^1(\partial W,\zz_2)\cong\zz_2$.

\smallskip\noindent%
\slii Clearly, the composition $F_2:=F_1\circ F_1$ of a map $F_1:W\to W$ as in the part
\sli maps $\bfe$ onto $\bfe+2\bff$ and $\bff$ onto $\bff$. Further, by the
construction the map $F_1$ coincides with $F'_1$ near the boundary $\partial W$. But
the explicit formula shows that $F'_1\circ F'_1=\id_W$. Hence $F_2:W\to W$ is
identical near the boundary $\partial W$. 
\qed

\newsubsection[thm3]{Proof of {\mdseries \slsf{Theorem 3}}}
{\it Part \slip.} Let $\omega$ and $\omega_1$ be two symplectic forms on $X$ having equal
cohomology class, $[\omega]=[\omega_1]$. Perturbing these forms we may assume that their
cohomology class is generic enough.  Then by
\cite{La,La-McD,McD-2,McD-1,McD-Sa-2} there exists a diffeomorphism $F_1:X\to X$
with $\omega'_1=(F_1)_*\omega'$ where $\omega',\omega'_1$ are perturbed forms. Repeating the
argument from the proof of part \sliv of \lemma{pi0-dfrm} we conclude that
$F_1$ acts trivially on $\sfh_2(X,\zz)$. Deforming the perturbed forms back to
the given ones and using \lemma{pi0-dfrm}, \sliip, we conclude the existence
of a diffeomorphism $F_2:X\to X$ with the same properties: $\omega_1=(F_2)_*\omega$ and
$(F_2)_*=\id:\sfh_2(X,\zz)\to\sfh_2(X,\zz)$.

It follows from \slsf{Propositions \ref{ImMapYk} {\rm and} \ref{prop-MapY1}}
that any diffeomorphism $F_2:X\to X$ with trivial action on $\sfh_2(X,\zz)$ and
any symplectic form $\omega$ there exists an $\omega$-symplectomorphism $F_3:X\to X$ such
which has the same action on $\sfh_2(X,\zz)$, $\pi_1(X)$, and $\pi_2(X)$ as
$F_1$. Then the composition $F_4:=F_2\circ F_3\inv$ lies in $\diff_\bullet(X)$ and
transforms $\omega$ into $\omega_1$. Finally, part \sli of \slsf{Theorem 1} allows to
nullify $\tiw_2(F_4)$ and to construct $F\in\diff_\bsq(X)$ with
$\omega_1=F_*\omega$. This shows the transitivity of the action of $\Gamma_\bsq$ on
$\pi_0(\Omega(X,\omega))$. The injectivity of this action follows from part \slii of
\slsf{Theorem 1}.

\medskip\noindent%
{\it Part \sliip} is essentially equivalent to \slsf{Theorem 1}. Indeed, if we
denote by $\Gamma_\bullet(\omega)<\Gamma_\bullet$ the group of isotopy classes of homotopically trivial
$\omega$-symplectomorphisms, then the parts \sli and \slii of \slsf{Theorem~1}
assert the surjectivity and respectively the injectivity of the homomorphism
$\tiw_2:\Gamma_\bullet(\omega)\to\sfh^1(X,\zz_2)$. In particular, $\Gamma_\bullet(\omega)$ is the desired
splitting.

\medskip\noindent%
{\it Part \sliiip} follows from the results of \refsubsection{barMap}. Notice that
the assertion is void in the rational case $g(Y)=g(X)=0$. So we assume that
$g(Y)\geq1$. Let $J$ be an $\omega$-tamed almost complex structure and $\pr:X\to Y$ the
$J$-holomorphic ruling. Denote by $\mbfy^*=\{y^*_1,\ldots,y^*_\ell\}$ the images of
singular fibres of $\pr:X\to Y$, and let $\map(Y,\mbfy^*)$ be the corresponding
mapping class group. Further, let $\Delta_1,\ldots,\Delta_\ell\subset Y$ be sufficiently small disc
neighbourhoods of the points $y^*_i$, respectively.  In the case $\ell=0$ when $X$
is a product $X=Y{\tix}S^2$ or a skew-product $X=Y{\tix}S^2$ we fix any
non-empty disc $\Delta_0\subset Y$. Set $Y^\circ:=Y\bs(\cup_i\Delta_i)$ and $X^\circ:=\pr\inv(Y^\circ)$.
Then there exists a deformation $\omega'$ of the symplectic structure $\omega$ such that
$(X^\circ,\omega')$ is the product $(Y^\circ,\omega_Y){\times}(S^2,\omega_{S^2})$. Every element of
$\map(Y,\mbfy^*)$ can be realised by a diffeomorphism $f:Y\to Y$ which is
supported in $Y^\circ$ and preserves the form $\omega_Y$. For such $f:Y\to Y$ we define
the diffeomorphism $F_f:X\to X$ setting $F_f|_{X^\circ}:=(f|_{Y^\circ}){\times}\id_{S^2}$ and
extending it to the whole $X$ by identity. Then every such map $F_f$ is an
$\omega'$-symplectomorphism. By \lemma{pi0-dfrm}, every such map is isotopic to an
$\omega$-symplectomorphism. Define the subgroup $\Gamma_\sfm<\Gamma_\lor$ as the group of
isotopy classes of maps $F_f$ for all diffeomorphisms $f:Y\to Y$ supported in
$Y^\circ$.

It follows from the construction that $\Gamma_\sfm$ acts trivially on
$\sfh_2(X,\zz)$.  Fix a base point $y_0$ lying in some disc $\Delta_i$ and
different from $y^*_i$ and fix a base point $x_0$ lying on the fibre
$\pr\inv(y_0)$. As it was shown in the proof of \propo{ImMapYk}, the group
$\pi_2(X,x_0)$ is isomorphic to $\zz\oplus\big(\oplus_{i=1}^\ell\zz[\pi_1(X,x_0)]\big)$, where
$\zz[\pi_1(X,x_0)]$ is group ring of $\pi_1(X,x_0)$ taken in $\ell$ copies. The
natural action of $\pi_1(X,x_0)$ on $\zz[\pi_1(X,x_0)]$ given transport of the
base point $x_0$ along paths representing elements $\gamma$ of $\pi_1(X,x_0)$ is the
left diagonal multiplication:
$\hbox{$\gamma*(a_0;a_1,\ldots,a_\ell)$}=(a_0;\gamma{\cdot}a_1,\ldots,\gamma{\cdot}a_\ell)$. Let us denote by
$\bar\pi_2(X)$ the group of \slsf{coinvariants} of $\pi_2(X,x_0)$ with respect to
the action of $\pi_1(X,x_0)$, \ie, the quotient of $\pi_2(X,x_0)$ by the subgroup
generated by elements of the form $\gamma*a-a$ with $a\in\pi_2(X,x_0)$ and
$\gamma\in\pi_1(X,x_0)$. Then $\bar\pi_2(X)$ is independent of the choice of the base
point $x_0$ and the group $\Gamma_\lor(X)$ acts naturally on $\bar\pi_2(X)$. Further,
recall there exists the natural isomorphism $\pi_1(X,x_0)\cong\pi_1(Y,y_0)$ and that
the group $\map(Y)$ can be defined as the \slsf{group $\out(\pi_1(Y))$ of outer
 automorphisms}. The latter group is the quotient
$\aut(\pi_1(Y,y_0))/\inn(\pi_1(Y,y_0))$ of the whole automorphism group
$\aut(\pi_1(Y,y_0))$ by the subgroup of inner automorphisms which are
conjugations by elements of the group $\pi_1(Y,y_0)$. As we have shown in the
proof of \propo{ImMapYk}, the image of $\Gamma_\lor(X)$ in
$\map(Y){\times}\aut(\bar\pi_2(X))$ is isomorphic to $\barr\map(Y,\ell)$ and the
homomorphism $\Gamma_\lor(X)\to\map(Y){\times}\aut(\bar\pi_2(X))$ maps the subgroup
$\Gamma_\sfm<\Gamma_\lor$ isomorphically onto the image of $\Gamma_\lor(X)$ in $\aut(\bar\pi_2(X))$.

This shows part \sliiip.

\medskip\noindent%
{\it Part \slivp.} Recall that $\Gamma_\sfw$ denotes the image of $\Gamma_\lor$ in
$\aut(\sfh_2(X,\zz))$.  \slsf{Theorem 4'--4'{'}'} contains a description of
the group $\Gamma_\sfw(X)$ for different types of $X$. In particular, in case
a set $\cals(X)$ of generators of $\Gamma_\sfw(X)$ is given explicitly. Every such
generator is geometric in the sense that it is represented by an explicitly
given diffeomorphism unique up to isotopy. We claim that the subgroup
$\Gamma'_\sfw<\Gamma_\lor$ generated by those geometric generators is mapped by the
homomorphism $\Gamma_\lor\to\aut(\sfh_2(X,\zz))$ isomorphically onto $\Gamma_\sfw$.  Clearly,
this claim implies the assertion \slivp.

Consider first the case treated in \slsf{Theorem 4'--4'{'}} where the
number $\ell$ of blow-ups is big enough. In this case the group $\Gamma_\sfw$ is a
Coxeter group and the set of generators $\cals(X)$ is Coxeter system. In
particular, the set of Coxeter relations is defining: together with the set of
generators $\cals(X)$ it gives a presentation of the group
$\Gamma_\sfw$. Consequently, it is sufficient to show that the geometric
realisations of the generators $s_i\in\cals(X)$ satisfy the same relations.

The first group of Coxeter relation is the involutivity of each generator:
$s_i^2=1$. Recall that each generator $s_i$ except the last one $s_\ell$ is a
Dehn twist along some $(-2)$-sphere $S_i$. It is well-known that the square of
a Dehn twist along a $(-2)$-sphere is smoothly isotopic to identity, see
\refthm{twists-4-2}. In the case of exceptional spheres this result is proved
in \lemma{except-twist}.

In the second group of relation we have the generalised braid relations
between generators $s_i\neq s_j$. According to the corresponding number $m_{ij}$
in the Coxeter matrix we can have the commutativity $s_is_j=s_js_i$ in the
case $m_{ij}=2$, the usual braid relation $s_is_js_i=s_js_is_j$ in the case
$m_{ij}=3$, and the ``long braid relation'' $s_is_js_is_j=s_js_is_js_i$ in the
case $m_{ij}=4$.

As a preparation step in the proof of these relations we show that for an
appropriate integrable complex structure $J$ every homology class yielding a
generator $s_i\in\cals(X)$ is represented by a holomorphic curve. In certain
sense, this is a strong version of \lemma{s2-lagr-exists} in the case of
rational and ruled $4$-manifolds. Clearly, for such a construction it is
sufficient to start from some integrable structure $J_0$ on the minimal
manifold $X_0$ and show the blow-up centres. 

In the rational case we proceed as follows. We start with an arbitrary line
$L$ on $\cp^2$ and blow up any point $p_1$ on it. As the result, we obtain the
first exceptional rational curve $E_1$ and the curve $F_1$ as the proper
pre-image. (Notation $F_1$ symbolise that this curve is a fibre of a
$J$-holomorphic ruling.) As the next blow-up centre $p_2$ we take the
intersection point of $E_1$ and $F_1$. As the result this time we obtain the
second exceptional rational curve $E_2$, the curve $S_{1,2}$ as the proper
pre-image of $E_1$ in the homology class $\bfs_{1,2}=\bfe_1-\bfe_2$, and the
curve $E'_{1,2}$ as the proper pre-image of $F_1$ in the homology class
$\bfe'_{1,2}=\bfl-\bfe_1-\bfe_2$. The next blow-up centre $p_3$ is the
intersection point of $E_2$ and $E'_{1,2}$. This time we obtain the next
exceptional rational curve $E_3$, the curve $S_{2,3}$ as the proper pre-image
of $E_2$ in the homology class $\bfs_{2,3}=\bfe_2-\bfe_3$, and the curve
$S'_{123}$ as the proper pre-image of $E'_{1,2}$ in the homology class
$\bfs'_{123}=\bfl-\bfe_1-\bfe_2-\bfe_3$. Starting from $k=3$, the successive
blow-up centre $p_{k+1}$ is an arbitrary point on the last obtained
exceptional curve $E_k$ apart from the proper pre-image of previous
$(-2)$-curve(s) $S_{k-1,k}$ (and $S'_{123}$ in the case $k=3$). The result of
this blow-up is the new exceptional curve $E_{k+1}$ and the $(-2)$-curve
$S_{k,k+1}$ as the proper pre-image of $E_k$ in the homology class
$\bfs_{k,k+1}=\bfe_k-\bfe_{k+1}$. The outline of the blow-up sequence is
depicted on \reffig{rat-blow}.

\smallskip

{
\begin{center}
\begin{figure}[h!]  
\includegraphics[width=6.3in]{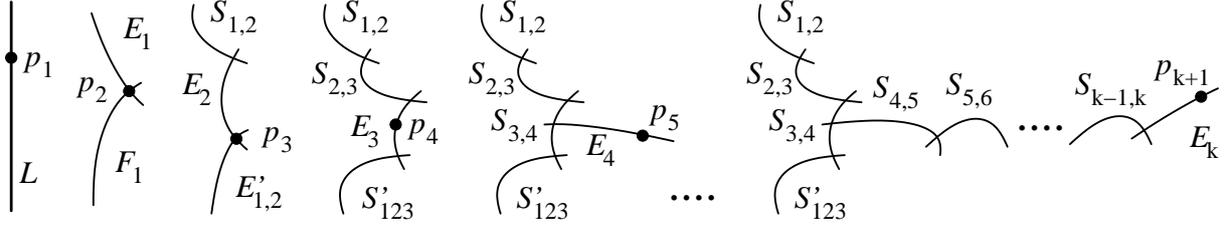}
  \caption{The blow-up constellation in the rational case.}
\label{rat-blow}
\end{figure}
\end{center}}

\vspace{-15pt}
In the irrational the procedure is completely analogous, an we simply indicate
the blow-up centres. We start with an arbitrary minimal irrational ruled
complex algebraic surface $X_0$, and the first blow-up centre is an arbitrary
point $p_1$ on an arbitrary fibre $F$. The second blow-up centre $p_2$ is the
intersection point of the exceptional curve $E_{1}$ and proper pre-image
$E'_1$ of the fibre $F$. Starting from $k=2$, the successive blow-up centre
$p_{k+1}$ is an arbitrary point on the last obtained exceptional curve $E_k$
apart from the proper pre-image of previous $(-2)$-curve(s) $S_{k-1,k}$ (and
$S'_{1,2}$ in the case $k=2$). The outline of the blow-up sequence in the
irrational case is depicted on \reffig{irrat-blow}.

\smallskip

{
\begin{center}
\begin{figure}[h!]  
\includegraphics[width=6.3in]{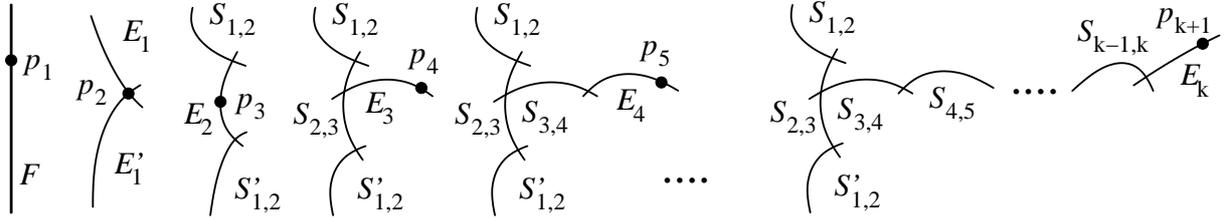}
  \caption{The blow-up constellation in the irrational case.}
\label{irrat-blow}
\end{figure}
\end{center}}

\vspace{-20pt}
Let us notice that in both rational and the irrational cases the incidence
graph of the obtained curves $S_{1,2},S_{2,3},\ldots,S_{k-1,k},E_k$ and $S'_{123}$
(or resp.\ $S'_{1,2}$ in the irrational case) is exactly the Coxeter graph
without markings. In particular, there exists a natural 1-1-correspondence
between these curves and the generators $s_0,s_1,\ldots\in\cals(X)$.

\smallskip
Now we are ready to prove the relations. The commutativity relations
$s_is_j=s_js_i$ follows from the fact that the generators
$s_i,s_j$ are (Dehn) twists along corresponding curves, and the curves
corresponding to $s_i,s_j$ are disjoint.

In the case of the usual braid relation $s_is_js_i=s_js_is_j$ we observe that
the generators $s_i,s_j$ can be realised as Dehn twists along corresponding
rational curves, and the curves
corresponding to $s_i,s_j$ are $(-2)$-spheres meeting transversally in one
point. So the relation is proved in \propo{braid-rel}.

\smallskip %
The last one is the ``long braid relation'' $s_is_js_is_j=s_js_is_js_i$. In
this case one of the generators, say $s_i$, is a Dehn twist along some
$(-2)$-sphere $S$ and another one $s_j$ is the twist along an exceptional
sphere $E$. Recall that in smooth category we have $s_i=s_i\inv$, so the
relation can be written as the commutation of $s_j$ with its conjugate:
$s_is_js_i\inv\cdot s_j=s_j\cdot s_is_js_i\inv$. Further, since $s_j$ is the twist $T_E$
along $E$, the conjugate diffeomorphism $T_S\vph T_E\vph T_S\inv$ is the twist
$T_{T_S(E)}$ along the image $T_S(E)$. Furthermore, since the question is
local it is sufficient to consider any ``standard'' local model. So we may
suppose that $S$ is a rational holomorphic $(-2)$-curve on some complex
surface, $E$ is an exceptional holomorphic curve with the intersection index
$[S]\cdot[E]=1$, and $U$ some their neighbourhood. After a generic deformation of
the complex structure the curve $E$ deforms into a new rational curve $E'$ in
the same homology class, the curve $C$ disappears, and the union $S+E$
deforms into an irreducible holomorphic curve $E_1$ in the homology class
$[S]+[E]$. Considering local holomorphic equations of the curves before and
after the deformations we see that $E_1$ is (isotopic to) the smooth connected
sum $S\#E$ in the sense of \refdefi{con-sum}. However, since the the
intersection index of $E$ and $E_1$ is zero, the curves $E'$ and $E_1$ must be
disjoint. This proves the desired property that the twists $T_E$ and
$T_{E_1}\simeq T_{T_S(E)}$ commute in the diffeotopy group $\Gamma(X)$.

\medskip%
Now we consider the ``small'' manifolds treated in \slsf{Theorem 4$'{'}{'}$}.
Recall that $\wt\Gamma_H$ denotes the image of $\diff(X)$ in $[X,X]$ and that $\Gamma_H$
consists of those elements of $\wt\Gamma_H$ which preserve the orientation and the
positive cone. In the rational case we prove a stronger result:
$\Gamma=\Gamma_\bsq\rtimes\wt\Gamma_H$. Recall also that $\barr\map(Y,\ell)$ and
$\sfh_1(X,\zz_2)$ are trivial in the rational case.

In the case \sli when $X$ is $\cp^2$ the isomorphism $\Gamma=\Gamma_\bsq\rtimes\wt\Gamma_H$
follows from \lemma{except-twist}.  Indeed, every complex line $L$ is a
$(+1)$-sphere and the twist along $L$ generates in $\Gamma$ the subgroup $\zz_2$
which defines the splitting $\wt\Gamma_H\to\Gamma$ of the projection $\Gamma\to\wt\Gamma_H\cong\zz_2$.

The case \slii when $X$ is $\cp^1{\times}\cp^1$ is also obvious: Indeed, the
complex conjugations of the factors and the involution
$(z,w)\in\cp^1{\times}\cp^1\mapsto(w,z)\in\cp^1{\times}\cp^1$ generate in $\Gamma$ the subgroup
$\wt\Gamma'_H$ which projects isomorphically onto $\wt\Gamma_H$.

\smallskip %
In the case \sliii when $X$ is the $\cp^2$ blown up once we proceed as
follows. We realise $X$ as the blow-up of $\cp^2$ in the origin
$0\in\cc^2\subset\cp^2$. Let $F:\cc^2\bs\{0\}\to\cc^2\bs\{0\}$ be the inversion given by
$F:(z,w)\mapsto(\frac{z}{|z|^2+|w|^2}, \frac{w}{|z|^2+|w|^2})$. Since $F$ preserves
the complex lines in $\cc^2$ through the origin, it extends to a map $F:X\to
X$. %
It is easy to see that the map $F$ is smooth, involutive ($F^2=\id_X$), and
interchanges the infinity line $L_\infty=\cp^2\bs\cc^2$ with the exceptional curve
$E\subset X$. Moreover, the map $F:E\to L_\infty$ is holomorphic. Now we see that similarly
to the case $X=\cp^1{\times}\cp^1$ the map $F$ and the twists along $E$ and $L_\infty$
generate in $\Gamma$ the subgroup $\wt\Gamma'_H$ which projects isomorphically onto
$\wt\Gamma_H$.

\smallskip %
In the case \sliv when $X$ is the product $Y{\times}S^2$ and the genus of $Y$ is
$g(Y)\geq1$ we fix maps $I_Y:Y\to Y$ and $I_{S^2}:S^2\to S^2$ which are complex
conjugation with respect to some complex structures $J_Y$ and $J_{S^2}$. Then
the involutions $I_Y{\times}\id_{S^2}$ and $\id_Y{\times}I_{S^2}$ generate the subgroup
$\Gamma'<\diff(X)$ isomorphic $\zz_2{\times}\zz_2$ whose projection onto $\Gamma_H$ is an
isomorphism.

\smallskip %
The case \slv of \slsf{Theorem 4$'{'}{'}$} when $X$ is the skew-product
$Y{\tix}S^2$ and the genus of $Y$ is $g(Y )\geq1$ is treated in the proof of
\slsf{Theorem 4$'{'}{'}$} below on page \pageref{Gamma-W-small}.

\smallskip %
Consider the case \slvi when $X$ is the $(S^2{\times}S^2)\#\cp^2$. First, we prove
the isomorphism $\wt\Gamma_H\cong\sfw(\rml_3(4,\infty))$. It is easy to see that the vectors
$\mbfv_0:=\bff_1-\bff_2$, $\mbfv_1:=\mbfe_1:=\bff_2-\bfe$, and
$\mbfv:=\mbfe:=\bfe$ form a basis of the lattice
$\sfh_2(X,\zz)\cong\zz\lan\bff_1,\bff_2,\bfe\ran$. Moreover, one can explicitly
verify the following two facts:
\begin{enumerate}
\item The basis $\{\mbfe_0,\mbfe_1,\mbfe\}$ with
 $\mbfe_0:=\frac{\mbfv_0}{\sqrt2}$ induces the isomorphism
 $\sfh_2(X,\rr)\cong(V_{\cals(X)},\lan\cdot,\cdot\ran)$ where $\cals(X)$ is the system
 $\rml_3(4,\infty)$. In particular, the reflections $\sigma_0,\sigma_1,\sigma$ in
$\sfh_2(X,\rr)\cong\rr\lan\bff_1,\bff_2,\bfe\ran$ with respect to hyperplanes
orthogonal to the vectors $\mbfv_0$, $\mbfv_1$, and respectively $\mbfv$
satisfy the corresponding Coxeter relations;
\item The basis $\{\mbfv_0,\mbfv_1,\mbfv\}$ with induces the isomorphism
 $\sfh_2(X,\zz)\cong(V^\zz_{\cals(X)},\lan\cdot,\cdot\ran)$ and hence defines a
 crystallographic structure in $\rml_3(4,\infty)$ with the Dynkin diagram %
 $\xymatrix@C-10pt@R-13pt{s_0\ar@<1pt>@2{->}[r] &
  s_1\ar@<1pt>@1{-}[r]^\infty&s}\!\!$.
\end{enumerate}

Our next step is to prove the isomorphism $\wt\Gamma_H\cong\sfw(\rml_3(4,\infty))$. For this
purpose is sufficient to show that the reflections $\sigma_0,\sigma_1,\sigma$ can be
represented by diffeomorphisms whose isotopy classes satisfy the same Coxeter
relations. To do this, let us realise $S^2{\times}S^2$ as the Hirzebruch surface
$\boldsymbol{\Sigma}_2$. Let $\pr:\bfSig_2\to\cp^1$ be the holomorphic ruling. Fix
some K\"ahler form $\omega$ on $\bfSig_2$. Let $Y_{\!-2}$ be the holomorphic
exceptional section of the ruling $\pr:\bfSig_2\to\cp^1$, $Y_2$ any other
holomorphic section, and $F_2,F_2'$ any two distinct fibres.  Denote by
$\bfy_{\!-2},\bfy_2,\bff_2\in\sfh_2(\bfSig_2,\zz)$ the corresponding homology
classes.  Further, smoothing appropriately the unique intersection point of
$Y_{\!-2}$ and $F'_2$ we obtain an $\omega$-symplectic sphere $F_1$ representing
the homology class $\bfy_{\!-2}+\bff_2$ and having a unique intersection point
with $F_2$. Let $\bff_1$ be the homology class of $F_1$. Then
$\bfy_{\!-2}=\bff_1-\bff_2$. Observe there exists a diffeomorphism
$\Phi:\bfSig_2\xrar{\;\cong\;}S^2{\times}S^2$ which maps the spheres $F_1$ and $F_2$ in
some horizontal and respectively vertical fibre of $S^2{\times}S^2$.  To construct
such a map $\Phi:\bfSig_2\to S^2{\times}S^2$, it suffices to take a generic enough
$\omega$-tame almost complex structure $J$ such that $F_1,F_2$ are $J$-holomorphic
and apply \refcorol{generic}, \sliip, to $J$-holomorphic deformations of
$F_1,F_2$. Now let us blow-up $\bfSig_2$ in any point lying on $F_2$ apart
from the intersection point with with $Y_{-2}$. Denote by $X_1$ the obtained
complex surface, by $E$ the arising exceptional curve, by $\bfe$ its homology
class, by $E_1$ the proper pre-image of the fibre $F_2$, and by $\bfe_1$ the
homology class of $E_1$. Identifying $X_1$ with $X=(S^2{\times}S^2)\#\cp^2$ we
obtain the equalities $\mbfv_0=\bfy_{-2}$, $\mbfv_1=\bfe_1$, and $\mbfe=\bfe$
of the homology classes. Consequently, the (Dehn) twists along $Y_{-2}$, $E_1$
and $E$ define the lifts of the reflections $\sigma_0,\sigma_1,\sigma\in\aut(\sfh_2(X,\zz))$
to the group $\Gamma(X)$. Denote these lifts by $\mbfs_0,\mbfs_1,\mbfs$. Since
$Y_{-2}$ is disjoint from $E$, $\mbfs_0$ commutes with $\mbfs$. As we have
shown above, $\mbfs_0$ and $\mbfs_1$ satisfy the relation
$(\mbfs_0\mbfs_1)^2=(\mbfs_1\mbfs_0)^2$. Taking into account the relations
$\mbfs_0^2=\mbfs_1^2=\mbfs^2=\id$ we conclude that the subgroup $\wt\Gamma'_H(X)$
of $\Gamma(X)$ generated by $\mbfs_0,\mbfs_1,\mbfs$ is a quotient of the group
$\sfw(\rml_3(4,\infty))$. On the other hand, the image of $\wt\Gamma'_H(X)$ in
$\aut(\sfh_2(X,\zz))$ is isomorphic to $\sfw(\rml_3(4,\infty))$. This implies that
$\wt\Gamma'_H(X)$ is isomorphic to $\sfw(\rml_3(4,\infty))$ and that the Coxeter
relations between $\mbfs_0,\mbfs_1,\mbfs$ corresponding to the Coxeter graph
$\xymatrix@C-10pt@R-13pt{s_0\ar@<1pt>@2{-}[r] & s_1\ar@<1pt>@1{-}[r]^\infty&s}$
form a defining set of relations. Moreover, since $\wt\Gamma'_H(X)$ is a subgroup
of $\Gamma(X)$, we obtain a splitting $\Gamma(X)=K\rtimes\wt\Gamma_H(X)$ where $K$ is the
kernel of the homomorphism $\Gamma(X)\to\aut(\sfh_2(X,\zz))$. By \slsf{Theorem 1} and
since $\barr\map(S^2,\ell)$ and $\sfh^1(S^2,\zz_2)$ are trivial, the group $K$ is
isomorphic to $\Gamma_\bsq(X)$.

This proves the assertion \sliv of \slsf{Theorem 3} in the case
$X=(S^2{\times}S^2)\#\cp^2$. 

\medskip%
It remains to show the assertion \sliv of \slsf{Theorem 3} in the case
$X=(Y{\times}S^2)\#\cp^2$ when $Y$ is a real surface of genus $g(Y)\geq1$. Denote by
$X_0$ the product $Y{\times}S^2$ and equip it with a product complex structure
$J=J_Y{\times}J_{S^2}$. Identify $X$ with the blow-up of $X_0$ at some point lying
on a given vertical fibre $F=y_0{\times}S^2$. Let $E\subset X$ be the exceptional rational
curve and $E_1\subset X$ the proper pre-image of $F$. Then $E_1$ is also an
exceptional rational curve. Denote by $\bfe,\bfe_1\in\sfh_2(X,\zz)$ the
homology classes of $E$ and respectively $E_0$, and by $\mbfs,\mbfs_0\in\Gamma(X)$
the twists along spheres $E,E_1$, respectively. Repeating the argumentation
used in the case $X=(S^2{\times}S^2)\#\cp^2$ we show that the subgroup $\wt\Gamma'_H(X)$
of $\Gamma(X)$ generated by $\mbfs,\mbfs_0$ is the group $\sfw(\rmi_2(\infty))$ given by
the Coxeter graph $\xymatrix@C-10pt@R-13pt{s_0\ar@<1pt>@1{-}[r]^\infty&s}\!\!$, and
more over, under the homomorphism $\Gamma(X)\to\aut(\sfh_2(X,\zz))$ the group
$\wt\Gamma'_H(X)$ is mapped isomorphically onto the group $\wt\Gamma_H(X)$. This gives
us the desired splitting $\Gamma(X)=K\rtimes\wt\Gamma_H(X)$ where $K$ is the kernel of
the homomorphism $\Gamma(X)\to\aut(\sfh_2(X,\zz))$ and finishes the proof of the
assertion \sliv of \slsf{Theorem 3}.

\medskip\noindent%
{\it Part \slvp.} Consider first a special case when $X_0$ is a product
$Y{\times}S^2$. Equip $X_0$ with the product complex structure
$J_0=J_Y{\times}J_{S^2}$. Assume that $(Y,J_Y)$ admits a \slsf{real structure},
\ie, a $J_Y$-antiholomorphic involution $\sigma_Y:Y\to Y$. Fix also a
$J_{S^2}$-antiholomorphic involution $\sigma_{S^2}:S^2\to S^2$. Then $F_c:=\sigma_Y{\times}
\sigma_{S^2}$ has the desired properties. The general case is obtained applying
blow-ups and contraction which are compatible with the real structure $F_c$
above in the following sense: A blow-up centre can be only a real point, \ie,
a point fixed by the antiholomorphic involution. An exceptional curve can be
contracted only if it is invariant under the antiholomorphic involution.

\medskip\noindent%
{\it Part \slvip.} The intersection form on $\sfh_2(X,\rr)$ of every ruled
$4$-manifold has the Lorentz signature. Consequently, there exists an
orientation inverting diffeomorphism if and only if the rank of
$\sfh_2(X,\rr)$ is $2$. This occurs only in the cases $X=Y{\times}S^2$ and
$X=Y{\tix}S^2$. 

An explicit construction of such a diffeomorphism is as follows: In the case
$X=Y{\times}S^2$ one can take a map of the product form $F=f_Y{\times}\id_{S^2}$ or
$F=\id_Y{\times}f_{S^2}$ where $f_Y:Y\to Y$ or respectively $f_{S^2}:S^2\to S^2$ is an
orientation inverting diffeomorphism. In the case $X=Y{\tix}S^2$ we fix a
ruling $\pr:X\to Y$ and a Riemannian metric on $X$ whose restriction on every
fibre $F_y=\pr\inv(y)$ is isometric to the standard round metric on
$S^2$. Then the fibrewise antipodal map $F:X\to X$ gives the desired inversion of
the orientation.

\smallskip
This finishes the proof of \slsf{Theorem 3}. \qed

\newsubsection[Gamma-W-irr]{Action on second homology. The irrational case.}
In this subsection we give a proof of \slsf{Theorem 4''}. Therefore we assume
that $X$ is an $\ell$-fold blow-up of the product $Y{\times}S^2$ where $Y$ is a closed
oriented surface of genus $g\geq1$. Recall that $\bff$ is the homology class a
fibre of some ruling $\pr:X\to Y$.

\newlemma{bfF-preserved} Let $F:X\to X$ be any homotopy equivalence. Then
$F_*(\bff)=\pm\bff$. 
\end{lem}

\proof The result follows from the homotopic theory point of view on rulings
and the class $\bff$. Namely, recall that the surface $Y$ is an
Eilenberg-MacLane space of the type $K(1, \pi_1(X))$ and any ruling indices an
isomorphism of fundamental groups. Since any automorphism of $\pi_1(Y,y_0)$ can be
realised by a diffeomorphism $f:(Y,y_0)\to(Y,y_0)$, a continuous map $p:X\to Y$ is
homotopic to some ruling $\pr:X\to Y$ if and only if it induces an
isomorphism of fundamental groups. 
 
Let $[\omega_Y]$ be the generator of $\sfh^2(Y,\zz)$ inducing the given
orientation on $Y$, and $p:X\to Y$ any continuous map inducing an
isomorphism of fundamental groups. Then the class $p^*[\omega_Y]$ is Poincar\'e dual to
$\pm\bff$. This is an intrinsic homotopic characterisation of $\pm\bff$ which must
be preserved by any  homotopy equivalence.
\qed

\newprop{E-irrat} Under the hypotheses of \slsf{Theorem 4''}, let $f$ be any
automorphism of the $\sfh_2(X,\zz)$ such that $f_*\bff=\bff$. Then 
$f\in\sfw(\cals(X))$.
\end{prop} 

\proof Let us notice that the classes $\bfs'_{12}=\bff-(\bfe_1+\bfe_2)$ and
$\bfs_{i,i+1}=\bfe_i-\bfe_{i+1}$ are represented by embedded spheres. In the
case $\bfs'_{12}$ we put both blow-up centres $p_1,p_2$ on the same fibre
of the canonical ruling $\pr:Y{\times}S^2\to Y$, and then the proper pre-image of this
fibre in $X$ is the desired $(-2)$-sphere in the class $\bfs'_{12}$. In the
case $\bfs_{i,i+1}$ we blow up $Y{\times}S^2$ in pairwise distinct points
$p_1,\ldots,p_i,p_{i+2},\ldots,p_\ell$ creating exceptional spheres
$E_1,\ldots,E_i,E_{i+2},\ldots,E_\ell$, and then blow up in some point $p_{i+1}$ lying on
$E_i$. Now the desired $(-2)$-sphere in the class $\bfs_{i,i+1}$ is the proper
pre-image of $E_i$.

\smallskip%
Denote $\bfe'_i:=f_*\bfe_i$. We shall try to bring the classes $\bfe'_i$ back
into $\bfe_i$ applying elements from $\sfw(\cals(X))$.

Write $\bfe'_i=k_i\bff + m_i\bfy +\sum_jn_{ij}\bfe_j$. Then %
$m_i=\bff{\cdot}\bfe'_i=\bff{\cdot}\bfe_i=0$. From the equality
$\bfe'_i{}^2=\bfe_i^2=-\sum_jn^2_{ij}=-1$ we conclude that for each $i$ there
exists exactly one index $j=j(i)$ such that $n_{i,j(i)}=\pm1$ and the remaining
coefficients vanish, $n_{ij}$ for $j\neq j(i)$. Since
$\bfe'_i{\cdot}\bfe'_j=\bfe_i{\cdot}\bfe_j=0$. We conclude that the sequence
$j(1),j(2),\ldots,j(\ell)$ is a permutation of the sequence $1,2,\ldots,\ell$.  Further,
observe that the Dehn twist along the embedded sphere $S_{i,i+1}$ in the
homology class $\bfs_{i,i+1}=\bfe_i-\bfe_{i+1}$ permutes the classes $\bfe_i$
and $\bfe_{i+1}$. Hence we can bring the class $\bfe'_i$ to the form
$\epsilon_i\bfe_i+k_i\bff$ with $\epsilon_i=\pm1$. Applying \lemma{except-twist} we can make
$\epsilon_i=+1$, so that $\bfe'_i= \bfe_i+k_i\bff$.

Denote $\bfy':=f_*\bfy$ and write $\bfy'=a\bff + c\bfy +\sum_jb_j\bfe_j$.
Then $c=\bff{\cdot}\bfy'=\bff{\cdot}\bfy=1$.  From $\bfy'{\cdot}\bfe'_i=\bfy{\cdot}\bfe_i=0$
we conclude $b_i=k_i$. Finally, the equality
$\bfy'{}^2=\bfy^2=0$ yields $2a=\sum_jk_j^2$.

At this point we observe that the the Dehn twist along a sphere in the
homology class $\bfs'_{12}=\bff-(\bfe_1+\bfe_2)$ maps the class $\bfe_1$ onto
$\bff-\bfe_2$ and $\bfe_2$ onto $\bff-\bfe_1$ preserving the remaining classes
$\bff,\bfy,\bfe_j$, $j\neq1,2$.  The composition of this Dehn twist with the Dehn
twist along a sphere in the class $\bfs_{12}=\bfe_1-\bfe_2$ maps $\bfe_1$ onto
$\bff-\bfe_1$ and $\bfe_2$ onto $\bff-\bfe_2$ preserving the remaining
classes. Observing that classes $\bff-\bfe_i$ are represented by exceptional
spheres we conclude that for given $\epsilon_1=\pm1$, $\epsilon_2=\pm1$ there exists an element
$h_{\epsilon_1,\epsilon_2}\in\sfw(\cals(X))$ with $(h_{\epsilon_1,\epsilon_2})_*\bfe_1=\bfe_1+\epsilon_1\bff$ and
$(h_{\epsilon_1,\epsilon_2})_*\bfe_2=\bfe_2+\epsilon_2\bff$ which preserves the remaining classes
$\bff,\bfy,\bfe_j$, $j\neq1,2$.  Applying such an element $h$ and Dehn twists
permuting the classes $\bfe_i$ we can reduce the general situation
$\bfe'_i=\bfe_i+k_i\bff$ to the special case when $k_2,\ldots,k_\ell$ vanish. Further,
the composition $h_{+1,+1}$ with $h_{+1,-1}$ maps $\bfe_1$ onto $\bfe_1+2\bff$
preserving the remaining classes $\bff,\bfy,\bfe_j$, $j\neq1$. It follows that we
can make $k_1$ either $1$ or $0$. However, the relation $2a=\sum_jk_j^2$
forbids the possibility $k_1=1,k_2=\cdots=k_\ell=0$. This shows that we have made vanish
all $k_i$ and brought the classes $\bfe'_i$ back
into $\bfe_i$, applying only transformations from $\sfw(\cals(X))$.
\qed

\smallskip \statep Proof. of \slsf{Theorem 4''}. The previous proposition
provides {\sl Part} \sli and gives an algebraic description of the group
$\Gamma_\sfw<\aut(\sfh_2(X,\zz))$ as the stabiliser of the element $\bff$.

\medskip\noindent%
{\sl Part \sliip.}  Let us identify $\sfh^2(X,\zz)$ with $\sfh_2(X,\zz)$ using
Poincar\'e duality. In particular, $\bfy,\bff;\bfe_1,\ldots,\bfe_\ell$ become a basis
of $\sfh^2(X,\rr)$ and we can write any class $[\omega]\in\sfh^2(X,\rr)$ as
$[\omega]=\varsigma\bfy+\nu\bff-\sum\mu_j\bfe_j$ where $\varsigma=\int_\bff\omega$, $\nu=\int_\bfy\omega$ and
$\mu_i=\int_{\bfe_i}\omega$ are the periods of $[\omega]$. Then $[\omega]$ lies in the positive
cone if and only if $\varsigma>0$ and $\varsigma{\cdot}\nu>\sum\mu^2_i$. 


Let $\bfp$ be the domain in the positive cone given by the period conditions
\eqqref(period-cond-D), namely, by $\varsigma\geq\mu_1+\mu_2,\;\mu_1\geq\mu_2\geq\cdots\geq\mu_\ell\geq0$.  It follows
from \propo{P-fund} 
that $\bfp$ is a fundamental domain for the action of $\Gamma_\sfw$ on the positive
cone and implies  the second assertion of {\sl Part \sliip.}


\smallskip%
Let $\omega^*$ be a symplectic form with the cohomology class $[\omega^*]$ in
$\bfp$. Then there exist $\ell$ pairwise disjoint exceptional $\omega^*$-symplectic
spheres $E^*_1,\ldots,E^*_\ell$. As we have shown above, the homology classes of these
spheres are $\bfe^*_i= \pm\bfe_j+n_j\bff$. Since $\bff$ is the homology class of
a generic fibre of some $\omega^*$-symplectic ruling,
$0<\int_{E^*_i}\omega^*<\int_{\bff}\omega^*$. It follows that up to permutation we must
have $\bfe^*_i=\bfe_j$ or $\bfe^*_i=\bff-\bfe_j$, and then
$\bfe_j=\bff-\bfe^*_i$. In both cases the classes $\bfe_j$ are represented by
$\omega^*$-symplectic exceptional spheres.  This implies the first assertion of
{\sl Part \sliip.} Moreover, after reindexing and replacing $\bfe^*_i$ by
$\bff-\bfe^*_i$ if needed, we obtain the equalities $\bfe^*_i=\bfe_i$.

\smallskip%
Now let $F:X\to X$ be a diffeomorphism which preserves the class
$[\omega^*]\in\sfh^2(X,\rr)$. Then it preserves the period inequalities
$\varsigma^*\geq\mu^*_1+\mu^*_2$, $\mu^*_1\geq\mu^*_2\geq\cdots\geq\mu^*_\ell>0$ where $\varsigma^*=\int_\bff\omega^*$ and
$\mu^*_i=\int_{\bfe_i}\omega^*$ and permute the classes $\bfe_i$ or $\bff-\bfe_i$
preserving $\omega^*$-periods. Observe that the period conditions can be
reformulated as increasing of $\omega^*$-integrals in the sequence
$\bfe_\ell,\bfe_{\ell-1},\ldots,\bfe_2,\bfe_1,\bff-\bfe_2,\bff-\bfe_3,\ldots,\bff-\bfe_\ell$. 
This gives us the following possibilities for permutations:
\begin{itemize}
\item[(A)] arbitrary permutation of classes $\bfe_i,\bfe_{i+1},\ldots,\bfe_j$
 satisfying $\mu^*_{i-1}>\mu^*_i=\mu^*_{i+1}=\cdots=\mu^*_j>\mu^*_{j+1}$ with $j>i$. In the
 case $i=1$ the condition $\mu^*_{i-1}>\mu^*_i$ should read $\varsigma^*-\mu^*_2>\mu^*_1$, in
 the case $j=\ell$ the condition $\mu^*_j>\mu^*_{j+1}$ is void.
\item[(D)] permutations of pairs $(\bfe_1,\bff-\bfe_1),(\bfe_2,\bff-\bfe_2),
 \ldots,(\bfe_k,\bff-\bfe_k)$ %
 and switching the elements $\bfe_i\leftrightarrow\bff-\bfe_i$ in \emph{even} number of
 pairs in the case $\frac12\varsigma^*=\mu^*_1=\mu^*_{2}=\cdots=\mu^*_k>\mu^*_{k+1}$ with $k\geq2$.
\end{itemize} 
The condition of evenness of the number of switched pairs comes from the fact
that there is no diffeomorphism interchanging the classes $\bfe_{i_0}$ and
$\bff-\bfe_{i_0}$ for some fixed $i_0$ and preserving the classes $\bff,\bfy$
and $\bfe_i$ with $i\neq i_0$, see the last argument in the proof of
\propo{E-irrat}. The corresponding permutation groups are the symmetric group
$\sym_{k+1}$ with $k=j-i$ and the Weyl group $\sfw(\rmd_k)$.  Both are
Coxeter-Weyl groups of type $\rma_k$ and $\rmd_k$, respectively, and the
corresponding systems of generators are reflections with respect to
hyperplanes orthogonal to $\bfs_{i,i+1},\bfs_{i+1,i+2},\ldots,\bfs_{j-1,j}$ and
respectively or $\bfs'_{12}=\bff-\bfe_1-\bfe_2,\bfs_{1,2},\ldots,\bfs_{k-1,k}$.

\smallskip%
To finish the proof of \slsf{Theorem 4''} it remains to show that the classes
$\bfs'_{12}=\bff-\bfe_1-\bfe_2$ and $\bfs_{i,i+1}=\bfe_i-\bfe_{i+1}$ are
represented by $\omega^*$-Lagrangian spheres provided they are orthogonal to $\omega^*$.
In the case $\bfs_{i,i+1}=\bfe_i-\bfe_{i+1}$ this fact is proved in
\lemma{lagr-S2}. In the case $\bfs'_{12}=(\bff-\bfe_1)-\bfe_2$ we observe that
the class $\bfe'_1:=\bff-\bfe_1$ is also represented by an exceptional
symplectic sphere.
\qed

\newsubsection[Gamma-W-rat]{Action on the second homology in the rational case
 and SW-theory.}  %
Recall that in the irrational case the group
$\sfw(\cals))<\aut(\sfh_2(X,\zz))$ is distinguished by by purely algebraic
condition, namely, as the stabiliser of the class $\bff$. In turn, the class
$\bff$ admits (up to sign) a simple homotopic characterisation.  The main
difference in the rational case is that the group
$\sfw(\cals))<\aut(\sfh_2(X,\zz))$ is now distinguished by means of
Seiberg-Witten invariants, \ie, by the \emph{smooth topology}.

\begin{prop}\label{E-SW}  Under the hypotheses of \slsf{Theorem 4'}, let $S$
 be a smooth sphere of self-intersection $-1$, $-2$, $-3$, or $-4$ satisfying
 the conditions $[S]{\cdot}\bfe_i=:m_i\geq0$ and $m_1\geq m_2\geq\cdots\geq m_\ell$. Then
 $[S]{\cdot}\bfl<m_1+m_2+m_3$ unless $S=3m\bfe-m(\bfe_1+\cdots+\bfe_9)-2\bfe_{10}$ for
 some $m\geq2$ and then $[S]^2=-4$.
\end{prop}


\proof Let us assume the contrary. Set $q:=-[S]^2$ and $k:=[S]{\cdot}\bfl$, so
that $k\geq m_1+m_2+m_3\geq0$, and $[S]=k{\cdot}\bfl-\sum_im_i\bfe_i$. In particular,
$[S]^2=k^2-\sum_im_i^2=-q$. In the case, $k=1$ we must have $m_2=\cdots=m_\ell=0$ and
$m_1\leq1$, and hence $[S]^2\geq0$ in contradiction with the hypothesis. Thus we may
additionally suppose that $k\geq2$.

In this situation {\bf Theorem 1.2} from \cite{Fi-St} (see also \cite{Oz-Sz})
gives the inequality $k(k-3)<\sum_im_i(m_i-1)$. This inequality is based on
SW-theory and is one of the most important applications of this theory in
symplectic and smooth topology. Let us also notice that the paper \cite{Fi-St}
contains the proof only for the special case ({\bf Theorem 1.1}) when
$m_i=0,1$ or $2$. However, this restriction on $m_i$ is never used so that the
proof goes through in general situation. Substituting 
$k^2=\sum_im_i^2-q$, we obtain 
\begin{equation}\label{SW-1}
\textstyle
\sum_im_i\leq3k+q-1.
\end{equation}
Now let us increase the sum $\sum_im_i^2$ preserving the sum $\sum_im_i$ and
the conditions $k\geq m_1+m_2+m_3$, $m_1\geq m_2\geq\cdots\geq m_\ell\geq0$. Namely, if
 $k> m_1+m_2+m_3$ and $m_l>m_{l+1}=0$, then we replace $m_1$ by $m'_1:=m_1+1$
and  $m_l$ by $m'_l:=m_l-1$, and in the case
$m_i>m_{i+1}=\ldots=m_j>m_{j+1}\geq0$ for some $3\leq i<j$ we replace $m_{i+1}$ by
$m'_{i+1}:=m_{i+1}+1$ and $m_j$ by $m'_j:=m_j-1$, remaining the rest $m_i$
unchanged. Besides in the case $m_1\geq m_2>m_3$ we replace $m_1$ by $m'_1:=m_1+1$
and $m_2$ by $m'_2:=m_2-1$. Then $\sum_im_i^2$ will increase, and we can repeat
doing so until we come to the case when  $m:=m_3=m_2=m_4=\cdots m_l>m_{l+1}\geq0$ 
and $m_1=k-2m$. Setting $m':=m_{l+1}$ we obtain
$\sum_im_i=k+(l-3)m+m'$ and 
\begin{equation}\eqqno(SW-2)\textstyle
k^2=\sum_im_i^2-q\leq m_1^2+(l-1){\cdot}m^2+m'{}^2-q
=k^2-4km+(l+3){\cdot}m^2+m'{}^2-q.
\end{equation}
This yields $4km\leq(l+3)m^2+m'{}^2-q$ or
$4k\leq(l+3)m+m'{\cdot}\frac{m'}{m}-\frac{q}{m}$. On the other hand, %
$\sum_im_i\leq k+(l-3)m+m'\leq3k+q-1$ which means 
\begin{equation}\eqqno(SW-3)
(l-3)m+m'\leq2k+q-1
\end{equation}
Combining these inequalities, we obtain
\begin{equation}\eqqno(SW-4)\textstyle
2k\leq6m+\frac{m'(m'-m)}{m}+(q-1)-\frac{q}{m}\leq6m+(q-1)-\frac{q}{m},
\end{equation}
where the latter inequality follows from $0\leq m'<0$. So we can conclude the
strict inequality $k<3m$ in the case $q=1$ and non-strict one $k\leq3m$ in the
cases $q=2$ and  $q=3$. On the other hand, the condition $m_1\geq m_2$ is equivalent to
$k\geq3m$. So in the case  $q=1$ we obtain a contradiction which proves the
proposition. 

In the case $q=2$ or $q=3$ we conclude the equality $k=3m$. Now the inequality
\eqqref(SW-2) leads to $(l-9)m^2+m'{}^2\geq q$, and hence $l>9$ or $l=9$ and
$m'\geq2$.  On the other hand, in the case $k=3m$ the estimation \eqqref(SW-3)
transforms into $(l-9)m+m'\leq q-1$, which is possible only in the following
cases:\; (a) $l<9$;\;(b) $l=9$, $m'\leq q-1$; or\; (c) $l=10$, $m=2$, $m'=0$ and
$q=3$.  In the case $q=2$ this again leads to a contradiction. Thus the only
remaining case is $l=9$, $q=3$ and $m'=2$.  Here we recall that the estimate
$(l-9)m^2+m'{}^2\geq q$ was obtained from the relation $k^2=\sum_im_i^2-q$ after
making certain transformations of the sequence $(m_1,\ldots,m_\ell)$ which increases
the sum $\sum_im_i^2$. Notice that every such transformation makes at least one
of the inequalities $m_1\geq m_2\geq\ldots\geq m_l$ sharper. On the other hand, in the
resulting sequence we have only equalities $m_1=\ldots m_9\geq m_{10}=2$. This means
that none of the transformations was applied and our coefficients are
$m_1=\ldots=m_9=m\geq m_{10}=2$, $k=3m$, (and $m_{11}=\ldots=m_\ell=0$ if $\ell>10$). But in this
case we must have $[S]^2=-4$.  Since we are now in the case $q=3$, this is a
contradiction which proves the proposition.

The final case is $q=4$. In this case we can colclude from \eqqref(SW-4) a
weaker estimate $k\leq3m+1$. The subcase $k\leq3m$ can be treated by the same
argument as the case $q=3$ and $k\leq3m$, and we can conclude that the only
possibility in this situation is $m_1=\ldots=m_9=m\geq m_{10}=2$, $k=3m$, (and
$m_{11}=\ldots=m_\ell=0$ if $\ell>10$). Let us consider the case $q=4$ and $k=3m+1$. In
this case \eqqref(SW-3) transforms into $(l-9)m+m'\leq5$, whereas \eqqref(SW-2)
yields $(l-9)m^2+m'{}^2\geq4(m+1)$. The case $l<9$ is clearly impossible. In the
case $l=9$ we obtain $m'\leq5$ and $m'{}^2\geq4(m+1)\geq4((m'+1)+1)=4m'+8$ which is
equivalent to $(m'-2)^2\geq12$. This contradicts to $m'\leq5$ and so is also
impossible. In the case $l=10$ we obtain $m+m'\leq5$ and $m^2+m'{}^2\geq4(m+1)$. The
only possibility to fulfill the both inequality is $m=5$ and $m'=0$. In this
case we should have $k=16$ and $m_1=6$, $m_2=\cdots=m_{10}=5$, but for the value
$q$ would be $\sum_im_i^2-k^2=36+9{\cdot}25-256=5$ and not $4$. As above, in this
situation we argue as follows: The $m$-sequence $m_1=6$, $m_2=\cdots=m_{10}=5$
should apprear from some original sequence by means of certain transformation
described above, but the only possibilities to make an inverse transformation
are: (a) to replace $m_{10}=5$ and $m_{11}=0$ by $m^*_{10}=4$ and $m^*_{11}=1$
or (b) $m_{1}=16$ and $m_{11}=0$ by $m^*_{1}=15$ and $m^*_{11}=1$, and in the
both cases the value $q$ would drop by $>2$ to some $q^*<4$. The final case to
consider is $l\geq11$. Then $(l-9)m+m'\leq5$ and $m>m'\geq0$ can be satisfied only if
$l=11$, $m=2$, and $m'\leq1$, but un this case we can not get
$(l-9)m^2+m'{}^2\geq4(m+1)$. This conclusion finishes the proof.
\qed

\state Remark. 
It is worth to notice that for every $m\geq2$ the class
\[
[S]:=3m\bfe-m(\bfe_1+\cdots+\bfe_9)-2\bfe_{10}
\]
with $[S]^2=-4$ can be realised by an embedded rational algebraic curve $C$ on
$10$-fold blown-up $\cp^2$.  Namely, it was shown by Dolgachev \cite{Dol} that
for every $m\geq2$ there exists an elliptic pencil $C_\lambda$ on $\cp^2$ of degree
$3m$ with $9$ base points $p_1,\ldots,p_9$ each being an ordinary point of
multipicity $m$. Every such pencil has $12$ singular curves
$C_{\lambda_1},\ldots,C_{\lambda_{12}}$ whose singularities are the above base points
$p_1,\ldots,p_9$, each of multipicity $m$, and one more nodal point
$p_{10,i}\in C_{\lambda_i}$. Now, for each $i=1,\ldots,12$, the proper pre-image of the
curves $C_{\lambda_i}$ on $\cp^2$ blown-up in ten points $p_1,\ldots,p_9$ and $p_{10,i}$
gives us the desired curve $C$ in the homology class $[S]$.

\newprop{G-W-rat} Under the hypotheses of \slsf{Theorem 4'}, the group
$\Gamma_\lor(X)$ coincides with $\sfw(\cals(X))$.
\end{prop}

\proof First, let us observe that every generator of $\sfw(\cals(X))$ is
represented by some diffeomorphism. \lemma{except-twist} shows this for the
generator $s_\ell$ inverting the class $\bfe_\ell$. In the case of the remaining
classes $\bfs_{i,i+1}=\bfe_i-\bfe_{i+1}$ and
$\bfs'_{123}=\bfl-(\bfe_1+\bfe_2+\bfe_3)$ we need embedded $(-2)$-spheres
representing these classes. Such spheres can be constructed by appropriate
blow-ups. For example, in the case $\bfs'_{123}=\bfl-(\bfe_1+\bfe_2+\bfe_3)$
we blow up $\cp^2$ in 3 points $p_1,p_2,p_3$ lying on some complex line $L$,
and then the proper pre-image of $L$ is the desired rational $(-2)$-curve.

Now let $F:X\to X$ be any diffeomorphism. Denote by $\bfe$ the image $F_*\bfe_\ell$
and write $\bfe=k\bfl -\sum_im_i\bfe_i$. We are going to apply to $\pm\bfe$
generators of the group $\sfw(\cals(X))$ trying to make $|k|$ as small as
possible.  Our goal is $k=0$, since in this case we must have
$\sum_im_i^2=-\bfe^2=1$ and hence $\bfe=\pm\bfe_j$ for some $j$.

Since we consider the class $\bfe\/$ up sign inversion, we may assume that
$k>0$. Observe that the group $\sfw(\cals(X))$ contains involutions inverting
the vectors $\bfe_i$. Thus we may assume that $m_i\geq0$. Also, the group
$\sfw(\cals(X))$ contains elements realising any given permutation of classes
$\bfe_i$. So we may additionally suppose the condition $m_1\geq m_2\geq\cdots\geq m_\ell\geq0$. In
this situation \propo{E-SW} gives us the inequality $k<m_1+m_2+m_3$.  Now let
us apply the generator $s_0$ which is the reflection with respect to the class
$\bfs'_{123}=\bfl-(\bfe_1+\bfe_2+\bfe_3)$. Then $s_0(\bfe)=
\bfe+(k-(m_1+m_2+m_3)){\cdot}\bfs'_{123}$, and in particular, %
the new coefficient $k$ is $k'=k+(k-(m_1+m_2+m_3))<k$. Notice also that
$k'>-k$ since otherwise we would have $m_1+m_2+m_3\geq3k$ and hence $m_1\geq k$. But
then we must have $m_2^2+\cdots+m_\ell^2=1+k^2-m_1^2\leq1$ which would imply $m_1=k\geq1$,
$m_2=1$, $m_2=\cdots=m_\ell=0$, and $m_1+m_2+m_3\geq3k$ would be impossible.  It follows
that in the case $k\neq0$ we can always decrease $|k|$ using the procedure above
until we come to the case $k=0$, and then $\bfe=\pm\bfe_j$ for some $j$. Then
permuting the classes $\bfe_i$ and inverting the sign if needed we can make
$\bfe=\bfe_\ell$.

Summing up we conclude that in the case $\ell\geq3$ we can replace the
diffeomorphism $F:X\to X$ by another diffeomorphism $F_1:X\to X$ which fixes the
class $\bfe_\ell$. Observe that we can repeat the same construction with $F_1$
and the class $(F_1)_*\bfe_{\ell-1}$ letting $\bfe_\ell$ fixed, obtaining a new
diffeomorphism $F_2:X\to X$ which fixes the classes $\bfe_\ell$ and
$\bfe_{\ell-1}$. Continuing this procedure, we finally come the case of a
diffeomorphism $F:X\to X$ which preserves classes $\bfe_3,\bfe_4,\ldots,\bfe_\ell$.
This means that $F$ acts on $\sfh_2(X,\zz)$ by automorphisms of the sublattice
$\zz\lan\bfl,\bfe_1,\bfe_2\ran\subset\sfh_2(X,\zz)$. 

It appears that the automorphism group of the sublattice
$\zz\lan\bfl,\bfe_1,\bfe_2\ran$ can not be generated by reflections conjugated
to the reflections $s_0,s_1$ and the reflection $s^*_2$ with respect to the
hyperplane orthogonal to $\bfe_2$, see part \slvi of \slsf{Theorem 4'{'}'}
where the precise description of the situation is given. However, enlarging
the lattice to $\zz\lan\bfl,\bfe_1,\bfe_2,\bfe_3\ran$ makes room for needed
reflections. So let us assume that $\varphi$ is an automorphism of
$\zz\lan\bfl,\bfe_1,\bfe_2,\bfe_3\ran$ preserving the positive cone.  Denote
$\bfl':=\varphi(\bfl)$ and write $\bfl':= k\bfl-\sum_{i=1}^3m_i\bfe_i$. The condition
on positive cone means that $k>0$. Applying reflections inverting the vectors
$\bfe_i$, if needed, we may assume that $m_i\geq0$. Besides, we can permute the
vectors $\bfe_i$ making $m_1\geq m_2\geq m_3$.

We claim that $k<m_1+m_2+m_3$ unless $k=1$ and $m_1=m_2=m_3=0$. Indeed,
otherwise we must have $k\geq m_1+m_2+m_3$ and then 
\[\textstyle
1=k^2-\sum_im_i^2\geq(m_1+m_2+m_3)^2-\sum_im_i^2\geq2m_1m_2,
\]
which would imply $m_2=m_3=0$ leading to a contradiction. Now, the same
argument as above shows that the reflection with respect to
$\bfs'_{123}=\bfl-(\bfe_1+\bfe_2+\bfe_3)$ transforms $\bfl'$ into
$\bfl'+(k-(m_1+m_2+m_3)){\cdot}\bfs'_{123}=:\bfl^*$, and new coefficient $k$ is
$k^*=k+(k-(m_1+m_2+m_3))<k$. Since $\bfl^*$ stays in the positive cone, we
must have $k^*>0$.  Thus repeating the procedure we come to the case $k=1$. 
This means that the obtained automorphism $\varphi$ fixes the vector $\bfl$. 

It follows that every automorphism $\varphi$ of the lattice
$\zz\lan\bfl,\bfe_1,\bfe_2,\bfe_3\ran$ which preserves the positive cone lies
in $\sfw(\cals(X))$. Indeed, applying the procedure above we write such an
automorphism $\varphi$ as a product $\varphi_1\varphi_2$ where $\varphi_1$ lies in $\sfw(\cals(X))$
and $\varphi_2$ preserves the vector $\bfl$. In particular, $\varphi_2$ maps each class
$\bfe_j$ onto some vector $\bfe'_j$ with square $(\bfe'_j)^2=-1$ still lying
on the lattice $\zz\lan\bfl,\bfe_1,\bfe_2,\bfe_3\ran$ and orthogonal to
$\bfl$. Such vector $\bfe'_j$ must be one of $\pm\bfe_1,\pm\bfe_2,\pm\bfe_3$, and
hence $\varphi_2$ permutes the classes $\bfe_1,\bfe_2,\bfe_3$ and inverts their
signs. Consequently, $\varphi_2$ also lies in $\sfw(\cals(X))$.

This finishes the proof of the proposition.
\qed

\smallskip%
\state Remark. \label{remark-SW} It is interesting to observe that for $\ell\leq9$
the assertion of \propo{E-SW} can be proved purely algebraically without
application of SW-theory. To show this we consider the following problem: {\it
 for which $\ell$ there exists integers $k$ and $m_1\geq m_2\geq\cdots\geq m_\ell\geq0$ satisfying
 $m_1+m_2+m_3\leq k$ and $\sum_im_i^2>k^2$?} Applying the argumentation similar to
one use in \propo{E-SW} one shows that the maximum of the quadratic form
$\sum_im_i^2$ under the conditions $m_1\geq m_2\geq\cdots\geq m_\ell\geq0$ and $m_1+m_2+m_3\leq k$ is
achieved for $m_1=k-2m$ and $m_2=\cdots=m_\ell=m$, and the value is
$(k-2m)^2+(\ell-1){\cdot}m^2=k^2+(\ell+3){\cdot}m^2-4km$. This value is $>k^2$ if and only
if $(\ell+3)m>4k$. Since $m_1=k-2m\geq m_2=m$ is equivalent to $k\geq3m$, we come to
the condition $(\ell+3){\cdot}m>3{\cdot}4{\cdot}m=12m$ which means $\ell>9$. And indeed, for
$\ell=10$ the values $k=3$ and $m_1=\cdots=m_{10}=1$ provide the simplest
counterexample to the algebraic counterpart of \propo{E-SW}.

The meaning of the above observation is that in the case when $X$ is the
$\cp^2$ blown up in $\ell\leq9$ points every element of $\aut(\sfh_2(X,\zz))$ 
can be represented by a diffeomorphism, whereas for $\ell\geq10$ this is not the
case.

\medskip%
\statep Proof. of \slsf{Theorem 4'.} {\it Part \slip} is simply reformulation
of the previous proposition.

\smallskip%
{\it Part \sliip.} 
Here we simply repeat the arguments used in the proof of
\slsf{Theorem 4''}, \sliip, making appropriate modifications. In particular,
we show that each class $\bfs'_{123}=\bfl-(\bfe_1+\bfe_2+\bfe_3)$ and
$\bfs_{i,i+1}=\bfe_i-\bfe_{i+1}$ can be realised by a Lagrangian sphere for an
appropriate symplectic form $\omega$ deformationally equivalent to the form $\omega^*$
as in the hypotheses of \slsf{Theorem 4'}. As it was shown at the end of the
proof of \slsf{Theorem 4''}, this fact follows from \lemma{lagr-S2}.


\smallskip%
{\it Part \sliiip.} Observe that in order to have a \emph{Lagrangian} sphere
representing a class $\bfs'_{123},\bfs_{1,2},\ldots$ the corresponding period
inequality must be an equality: $\lambda=\mu_1+\mu_2+\mu_3,\mu_1=\mu_2$ and so on. Moreover,
it follows from  \lemma{defla(-2)} that the equality is also a sufficient
condition for the existence of such a Lagrangian sphere. Thus the problem can
be reformulated as follows: how many equalities and at which places can we
have equalities in the period inequalities. Further, observe that the system
$\rme_\ell$ corresponds to the case when there are no strict inequalities,
whereas each maximal system listed in the case $\ell\geq9$ is obtained from the
system $\rme_\ell$ by removing one element and changing the corresponding period
equality into a strict inequality.

Further, let us notice that the anticanonical divisor $-K_X$ of the $\ell$-fold
blown-up $\cp^2$ has the homology class $3\bfl-\sum_{i=1}^\ell\bfe_i$, and so the
``periods'' of $-K_X$ are $\bfl{\cdot}(-K_X)=3$, $\bfe_i{\cdot}(-K_X)=1$ giving the
equality in each period inequality. Since $(-K_X)^2=9-\ell$, exists no symplectic
form $\omega$ on a rational symplectic $4$-manifold satisfying the period
conditions $\frac\lambda3=\mu_1=\ldots=\mu_9$. On the other hand, $-K_X$ is \emph{ample} for
$\ell\leq8$. Therefore its Poincar\'e dual is represented by a K\"ahler form $\omega_{-K}$,
which gives us the existence of a symplectic form admitting the
$\rme_\ell$-configuration of Lagrangian spheres for $\ell\leq8$.

It remains to show the existence of a symplectic form on the $\ell$-fold blown-up
$\cp^2$ with $\ell\geq9$ whose conditions periods correspond to one of the maximal
systems listed in the theorem. The case $\rmd_{\ell-1}$ is in fact already
contained in the proof of \slsf{Theorem 4''}: let us simply notice that in the
case $\ell\geq2$ the $\ell$-fold blow-up of $\cp^2$ is the $(\ell-1)$-fold blow-up of
$S^2{\times}S^2$, and the proof given for the case of a blow-up of $Y{\times}S^2$
applies here without changes. The case $\rma_{\ell-1}$ is easy: we just need to
find $\ell$ small symplectic balls $B_i$ in $\cp^2$ of the same radius $r$ and
blow them up.  The case $\rme_k\oplus\rma_{\ell-k-1}$ with $3\leq k\leq8$ is essentially the
combination of the cases $\rme_k$ and $ \rma_{\ell-k-1}$: Start with some symplectic $k$-fold
blow-up of $\cp^2$ having equalities in all period inequalities, and then 
blow up $\ell-k-1$ small disjoint symplectic balls.

The final case to consider is $\rma_1\oplus\rma_{\ell-2}$. Here the period conditions
are $\mu_1=\mu_2>\lambda-2\mu_1=\mu_3=\ldots=\mu_\ell$ and we proceed as follows: Fix some
sufficiently small $\mu>0$ and consider the product $X_0:=S^2{\times}S^2$ equipped
with the product symplectic form $\omega_0$ such that the volumes of fibres are
$\int_{\bff_1}\omega_0=\int_{\bff_2}\omega_0=\frac{1+\mu}{2}$. Fix $\ell-2$ disjoint
symplectic balls on $X_0$ of radius $r=\sqrt{\mu/\pi}$, indexed by
$B_0,B_3,\ldots,B_\ell$, and let $(X,\omega)$ be the blow-up of $(X_0,\omega_0)$ made using
these balls. Let $E_0,E_3,\ldots,E_\ell$ be the arising exceptional symplectic spheres
and $\bfe_0,\bfe_3,\ldots,\bfe_\ell$ their homology classes. Observe that the classes
$\bfe_1:=\bff_1-\bfe_0$ and $\bfe_2:=\bff_2-\bfe_0$ are also represented by
exceptional symplectic spheres. Similarly, the class
$\bfl:=\bfe_0+\bfe_1+\bfe_2$ is represented by
a symplectic $(+1)$-sphere. By the construction,  the periods of $(X,\omega)$ 
with respect to the basis $\bfl;\bfe_1,\ldots,\bfe_\ell$ have the type
$\rma_1\oplus\rma_{\ell-2}$ as desired.

\smallskip%
{\it Part \slivp}  follows from the previous assertions of the theorem.
\qed

\subsection{The case of small manifolds.}\label{Gamma-W-small} Here we give 

\statep Proof. of \slsf{Theorem 4'{'}'.} {\it Parts \slip--\slvp} become
obvious after taking into account \lemma{bfF-preserved}, except the existence
of an orientation inverting diffeomorphism $F$ of $Y{\tix}S^2$ in {\it Parts
 \slvp} and \sliii (here $Y=S^2$). To construct $F$, we consider a principle
$\SO(2)$-bundle $\pi:\bfp\to Y$ with $c_1(\bfp)=1$ and realise $X=Y{\tix}S^2$ the
the associated fibre bundle $\bfp{\times_{\SO(2)}}S^2$. Let $\theta\in[0,\pi]$ and
$\varphi\in[0,2\pi]$ be the standard spherical coordinates on $S^2$. Then the desired
map $F:X\to X$ can be defined as the fibrewise map given by $(\theta,\varphi)\mapsto(\pi-\theta,\varphi)$.

\smallskip%
{\it Part \slvip.} It is clear that the reflections
$\mbfs_1,\mbfs_2,\mbfs^*_0\in\aut(\sfh_2(X,\zz))$ are represented by
diffeomorphisms. It was shown in \lemma{cryst-E8}, \slivp, that these
reflections generate a subgroup of $\aut(\sfh_2(X,\zz))\cong\bfo(1,2,\zz)$
isomorphic to $\sfw(\rml_3(4,\infty))$. Further, since the class $\bfl$ and its
images $\mbfs_1(\bfl),\mbfs_2(\bfl),\mbfs^*_0(\bfl)$ lie in the positive cone
of $\sfh_2(X,\rr)$, $\sfw(\rml_3(4,\infty))$ is a subgroup of $\bfo_+(1,2,\zz)$
preserving the positive cone. So we only need to show that $\sfw(\rml_3(4,\infty))$
is is the whole group $\bfo_+(1,2,\zz)$. Take any $A\in\bfo_+(1,2,\zz)$ and set
$\wt\bfl:=A(\bfl)$. Write $\wt\bfl=l\bfl-m_1\bfe_1-m_2\bfe_2$. Applying the
reflections $\mbfs_1,\mbfs_2$ we can make $m_1,m_2$ non-negative. Since
$\wt\bfl{}^2=l^2-m_1^2-m_2^2=1$, in the case $l\geq2$ we must have $l<m_1+m_2$.
So applying $\mbfs^*_0$ we can decrease the coefficient $l$. This means that
the $\sfw(\rml_3(4,\infty))$-orbit of $\wt\bfl$ contains a vector
$\bfl'=\bfl-m'_1\bfe_1-m'_2\bfe_2$.  But in this case $m_1=m_2=0$ and hence
$\bfl'=\bfl$. Consequently, the case of a general element  $A\in\bfo_+(1,2,\zz)$
reduces to the case when $A$ preserves the class $\bfl$. In this case $A$ can
only exchange the classes  $\bfe_1,\bfe_2$ and change their sings. Hence such
a special element $A$ also lies in  $\sfw(\rml_3(4,\infty))$, and so we conclude the
desired isomorphism $\wt\Gamma_H\cong\sfw(\rml_3(4,\infty)))\cong\bfo_+(1,2,\zz)$.

\smallskip%
{\it Part \slviip.} Since the classes $\bfe_1$ and $\bfe'_1=\bff-\bfe_1$ are
represented by $(-1)$ spheres, the reflections $s_1$ and $s^*_1$ are
represented by diffeomorphisms. Further, a class
$a\bff+b\bfy+c\bfe_1\in\sfh_2(X,\rr)$ has positive square iff $2ab>c^2$, and
lies in the positive cone iff, in addition, $a>0$. So by
\lemma{bfF-preserved}, a diffeomorphism of $X=(Y{\times}S^2)\#\barr\cp^2$ preserves
the positive cone if and only if it preserve the class $\bff$. The subspace
$(\bff)^\perp$ of $\sfh_2(X,\rr)$ orthogonal to $\bff$ admits the basis
$\bfe_1,\bfe'_1$. It is not difficult to see that this basis induces the
natural isomorphism $(\bff)^\perp\cong V_{\cals(X)}$ with the geometric representation
of the group $\sfw(\cals(X))$ with the Coxeter system $\cals(X):=\{s_1,s^*_1\}$,
and that the lattice $\zz\lan\bfe_1,\bfe'_1\ran=:V^\zz_{\cals(X)}$ induces the
natural crystallographic structure in $(\bff)^\perp\cong V_{\cals(X)}$.

\smallskip%
This finishes the proof of \slsf{Theorem 4'{'}'.}
\qed



\end{document}